\documentclass[10pt]{amsart}
\usepackage{amsmath}
\usepackage{amsxtra}
\usepackage{amscd}
\usepackage{amsthm}
\usepackage{amsfonts}
\usepackage{amssymb}
\usepackage{eucal}

\theoremstyle{definition}

\theoremstyle{remark}

\newcommand{\thmref}[1]{Theorem~\ref{#1}}
\newcommand{\secref}[1]{\S~\ref{#1}}
\newcommand{\lemref}[1]{Lemma~\ref{#1}}
\newcommand{\defref}[1]{Definition~\ref{#1}}
\newcommand{\propref}[1]{Proposition~\ref{#1}}
\newcommand{\corref}[1]{Corollary~\ref{#1}}

\newcommand{\remref}[1]{Remark~\ref{#1}}
\newcommand{\exref}[1]{Example~\ref{#1}}

\newcommand{\nc}{\newcommand}
\nc{\renc}{\renewcommand}
\nc{\ssec}{\subsection}
\nc{\sssec}{\subsubsection}
\nc{\on}{\operatorname}

\nc\ol{\overline}
\nc\wt{\widetilde}
\nc\wh{\widehat}
\nc\tboxtimes{\wt{\boxtimes}}

\nc{\Aa}{{\mathbb{A}}}
 \nc{\Gg}{{\mathbb{G}}}  
\nc{\Hh}{{\mathbb{H}}}
 \nc{\Nn}{{\mathbb{N}}}
\nc{\Pp}{{\mathbb{P}}}
\nc{\Rr}{{\mathbb{R}}}
\nc{\BV}{{\mathbb{V}}}
\nc{\BW}{{\mathbb{W}}}
\nc{\Zz}{{\mathbb{Z}}}
\nc{\Qq}{{\mathbb{Q}}}
\nc{\Ss}{{\mathbb{S}}}
\nc{\Cc}{{\mathbb{C}}}

\nc{\CA}{{\mathcal{A}}}
\nc{\CB}{{\mathcal{B}}}

\nc{\CE}{{\mathcal{E}}}
\nc{\CF}{{\mathcal{F}}}
\nc{\CG}{{\mathcal{G}}}
\nc{\CL}{{\mathcal{L}}}
\nc{\CC}{{\mathcal{C}}}
\nc{\CM}{{\mathcal{M}}}
\def\Mm{\CM}
\nc{\CN}{{\mathcal{N}}}
\nc{\Oo}{{\mathcal{O}}}
\nc{\CP}{{\mathcal{P}}}
\nc{\CQ}{{\mathcal{Q}}}
\nc{\CR}{{\mathcal{R}}}
\nc{\CS}{{\mathcal{S}}}
\renc{\CD}{{\mathcal{D}}}
\nc{\CT}{{\mathcal{T}}}
\nc{\CU}{{\mathcal{U}}}
\nc{\CV}{{\mathcal{V}}}
\nc{\CK}{{\mathcal{K}}}
\nc{\CW}{{\mathcal{W}}}
\nc{\CZ}{{\mathcal{Z}}}

\nc{\fa}{{\mathfrak{a}}}
\nc{\fb}{{\mathfrak{b}}}
\nc{\fg}{{\mathfrak{g}}}
\nc{\fgl}{{\mathfrak{gl}}}
\nc{\fh}{{\mathfrak{h}}}
\nc{\fj}{{\mathfrak{j}}}
\nc{\fm}{{\mathfrak{m}}}
\nc{\fn}{{\mathfrak{n}}}
\nc{\fu}{{\mathfrak{u}}}
\nc{\fp}{{\mathfrak{p}}}
\nc{\fr}{{\mathfrak{r}}}
\nc{\fs}{{\mathfrak{s}}}
\nc{\fsl}{{\mathfrak{sl}}}
\nc{\hsl}{{\widehat{\mathfrak{sl}}}}
\nc{\hgl}{{\widehat{\mathfrak{gl}}}}
\nc{\hg}{{\widehat{\mathfrak{g}}}}
\nc{\chg}{{\widehat{\mathfrak{g}}}{}^\vee}
\nc{\hn}{{\widehat{\mathfrak{n}}}}
\nc{\chn}{{\widehat{\mathfrak{n}}}{}^\vee}

\nc{\fA}{{\mathfrak{A}}}
\nc{\fB}{{\mathfrak{B}}}
\nc{\fD}{{\mathfrak{D}}}
\nc{\fE}{{\mathfrak{E}}}
\nc{\fF}{{\mathfrak{F}}}
\nc{\fG}{{\mathfrak{G}}}
\nc{\fK}{{\mathfrak{K}}}
\nc{\fL}{{\mathfrak{L}}}
\nc{\fM}{{\mathfrak{M}}}
\nc{\fN}{{\mathfrak{N}}}
\nc{\fP}{{\mathfrak{P}}}
\nc{\fU}{{\mathfrak{U}}}
\nc{\fV}{{\mathfrak{V}}}
\nc{\fZ}{{\mathfrak{Z}}}

\nc{\bb}{{\mathbf{b}}}
\nc{\bc}{{\mathbf{c}}}
\nc{\bd}{{\mathbf{d}}}
\nc{\be}{{\mathbf{e}}}
\nc{\bj}{{\mathbf{j}}}
\nc{\bn}{{\mathbf{n}}}
\nc{\bp}{{\mathbf{p}}}
\nc{\bq}{{\mathbf{q}}}
\nc{\bF}{{\mathbf{F}}}
\nc{\bu}{{\mathbf{u}}}
\nc{\bv}{{\mathbf{v}}}
\nc{\bx}{{\mathbf{x}}}
\nc{\bs}{{\mathbf{s}}}
\nc{\by}{{\mathbf{y}}}
\nc{\bw}{{\mathbf{w}}}
\nc{\bA}{{\mathbf{A}}}
\nc{\bK}{{\mathbf{K}}}
\nc{\bI}{{\mathbf{I}}}
\nc{\bB}{{\mathbf{B}}}
\nc{\bG}{{\mathbf{G}}}
\nc{\bC}{{\mathbf{C}}}
\nc{\bD}{{\mathbf{D}}}
\nc{\bP}{{\mathbf{P}}}
\nc{\bH}{{\mathbf{H}}}
\nc{\bM}{{\mathbf{M}}}
\nc{\bN}{{\mathbf{N}}}
\nc{\bV}{{\mathbf{V}}}
\nc{\bU}{{\mathbf{U}}}
\nc{\bL}{{\mathbf{L}}}
\nc{\bT}{{\mathbf{T}}}
\nc{\bW}{{\mathbf{W}}}
\nc{\bX}{{\mathbf{X}}}
\nc{\bY}{{\mathbf{Y}}}
\nc{\bZ}{{\mathbf{Z}}}
\nc{\bS}{{\mathbf{S}}}
\nc{\bR}{{\mathbf{R}}}

\nc{\sA}{{\mathsf{A}}}
\nc{\sB}{{\mathsf{B}}}
\nc{\sC}{{\mathsf{C}}}
\nc{\sD}{{\mathsf{D}}}
\nc{\sF}{{\mathsf{F}}}
\nc{\sG}{{\mathsf{G}}}
\nc{\sK}{{\mathsf{K}}}
\nc{\sM}{{\mathsf{M}}}
\nc{\sO}{{\mathsf{O}}}
\nc{\sQ}{{\mathsf{Q}}}
\nc{\sP}{{\mathsf{P}}}
\nc{\sZ}{{\mathsf{Z}}}
\nc{\sfp}{{\mathsf{p}}}
\nc{\sr}{{\mathsf{r}}}
\nc{\sg}{{\mathsf{g}}}
\nc{\sff}{{\mathsf{f}}}
\nc{\sfb}{{\mathsf{b}}}
\nc{\sfc}{{\mathsf{c}}}
\nc{\sd}{{\ltimes}}

\nc{\tA}{{\widetilde{\mathbf{A}}}}
\nc{\tB}{{\widetilde{\mathcal{B}}}}
\nc{\tg}{{\widetilde{\mathfrak{g}}}}
\nc{\tG}{{\widetilde{G}}}
\nc{\TM}{{\widetilde{\mathbb{M}}}{}}
\nc{\tO}{{\widetilde{\mathsf{O}}}{}}
\nc{\tU}{{\widetilde{U}}}
\nc{\TZ}{{\tilde{Z}}}
\nc{\tx}{{\tilde{x}}}
\nc{\tq}{{\tilde{q}}}

\nc{\tfP}{{\widetilde{\mathfrak{P}}}{}}
\nc{\tz}{{\tilde{\zeta}}}
\nc{\tmu}{{\tilde{\mu}}}

  \nc{\Ob}{{\mathop{\operatorname{\rm Ob}}}}
  \nc{\Sym}{{\mathop{\operatorname{\rm Sym}}}}
   \nc{\Aut}{{\mathop{\operatorname{\rm Aut}}}}
 \nc{\Spec}{{\mathop{\operatorname{\rm Spec}}}}
  \nc{\spec}{{\mathop{\operatorname{\rm Spec}}}}
\nc{\Ker}{{\mathop{\operatorname{\rm Ker}}}}
 \nc{\dom}{{\mathop{\operatorname{\rm dom}}}}
\nc{\End}{{\mathop{\operatorname{\rm End}}}}
 \nc{\Hom}{\on{\Hom}}
 \nc{\GL}{{\mathop{\operatorname{\rm GL}}}}
 \nc{\Id}{{\mathop{\operatorname{\rm Id}}}}
 \nc{\rk}{{\mathop{\operatorname{\rm rk}}}} 
 \nc{\length}{{\mathop{\operatorname{\rm length}}}}
\nc{\supp}{{\mathop{\operatorname{\rm supp}}}}
\nc{\val}{{\rm val}}
\nc{\res}{{\mathop{\operatorname{\rm res}}}}

\def\tensor{{\otimes}}
\def\meet{\cap}
\def\union{\cup}

\def\g{\gamma}
\def\G{\Gamma}

\def\<{\begin}
 \def\>{\end}
\def\m{\setminus}

\nc{\seq}[1]{\stackrel{#1}{\sim}}

\def\inv{^{-1}}
\def\claim#1{{\noindent \bf Claim #1.\ }}
\def\case#1{{\noindent \bf Case #1.\ }}
\def\Claim{{\noindent \bf Claim.\ }}
\def\beq#1{{\begin{equation} \label{#1}}  }

\def\Uu{\mathbb U}

\def\prf{\begin{proof}}
\def\pv{\end{proof} }
 \def\eprf{\end{proof} }

\def\acl{\mathop{\rm acl}\nolimits}
 \def\dcl{\mathop{\rm dcl}\nolimits}

\def\liminv{\underset{\longleftarrow}{lim}\,}

\def\lbl#1{ {{\rm [#1]}}    \label{#1}  }
\def\a{\alpha}

\def\ba{\bar{a}}
\def\k{{\rm k}}

\def\cU{{\mathcal U}}
\def\Ff{{\mathcal F}}
 \nc{\K}{{\mathop{\operatorname{ { K} }}}}

\def\T{{\bf T}}
 \def\RV{{\rm RV}}
\def\rv{{\rm rv}}
\def\RES{{\rm RES}}

\nc{\valr}{{\rm val_{rv}}}
\def\VF{{\rm VF}}
\def\tX{{\tilde{X}}}

\def\eq#1#2{ \begin{equation} \label{#1} #2  \end{equation}}
 \def\Ll{{\mathbb L}}

\def\Ke{K_{exp}}

\def\bs{\backslash}

 \def\e{\epsilon}
\def\L{\Lambda}
\def\f{\mathfrak{f}}
\def\cc{{\frown}}

\begin{document}

\author{Ehud Hrushovski, David Kazhdan}
\address{\newline Institute of Mathematics, the Hebrew
University of Jerusalem, Givat Ram, Jerusalem, 91904, Israel.}  
\thanks{Supported by ISF grants   0397691 and 1048-078}
\email{ehud@math.huji.ac.il,kazhdan@math.huji.ac.il}
 \title{Motivic Poisson summation}

\def\mD{{\mathcal D}}
 \def\mDTc{{\bar{\mD}_T}}
 
\def\CR{\delta^K}

  \def\rR{\mathcal R}
  \def\bt{\bar{t}}  
\def\z{\zeta}
\def\tJ{\widetilde{J}}
\def\lf{\mathfrak l}
\def\tC{\widetilde{C}} 
 
\def\bO{\mathbf{O}}
 
\def\TR{(R \meet T)}
\def\bTR{(R \meet \bT)}

  \def\cct{\cc}  
\def\lbl{\label}
\<{abstract}
We develop a "motivic integration" version of the Poisson summation formula for function fields, with values in the Grothendieck ring of definable exponential sums.  We also study division algebras over the function field, and obtain relations among the motivic Fourier transforms of a test function  at different completions.
We use these to prove, in a special case, a motivic version of a theorem of [DKV].  
 \>{abstract}

\maketitle


\<{section}{Introduction}

  The first order theory of valued fields associated with number theory has received
a great deal of attention in the past half-century.    A region of mystery remains around local fields of positive characteristic, but 
by and large local fields and associated geometric structures, are decidable and accessible to model-theoretic tools; in the hands of Denef this has been useful 
in the study of $p$-adic integration, later leading to the motivic integration
of Kontsevich, Denef-Loeser and others.  By Feferman-Vaught methods (\cite{CK}, 6.2),
one can similarly understand products of local fields or of their rings of integers;
the underlying rings of the adeles are thus decidable; but without access
to the discrete global field embeded in the adeles, this    
permits rather limited contact with the global geometry.    No known decidable theory captures such a discrete embedding.  The closest approach is \cite{MvD}, that can be 
understood as the theory of the non-archimedean adeles over $\Qq^a$, with an embedded copy of $\Qq^a$; but as the authors make clear, it works precisely because of the cut-and-paste property, i.e. the absence of any global constraints.   Every global field or adelic construction whose first-order theory is understood at all,  is known to be undecidable.

 Based on this evidence, one might guess that   the line of decidability, for
fields associated with number theory, coincides with the local/global distinction.  
The history of number theory, however,  shows no such line at all; adelic methods are no less geometric than than local ones, and for two hundred years have consistently decided relevant problems.  We are not able 
to resolve the tension between these different conclusions in the present paper, but we try to reduce it a little.   We study function fields and their associated adeles.
We embed the function field only piecewise, as an Ind-definable object, and do not permit quantification over it.  But in this setting we are able to interpret the Poisson summation formula motivically, leading to connections between Denef-Loeser motivic integrals over distinct local fields.  

The term `motivic' in this paper is used in its sense in the context motivic integration, to say that numbers are replaced by elements of the Grothendieck ring of varieties,
and closely associated rings.  
  We discuss such rings in \S 2 and \S 3; in particular we define the ring of exponential
  sums $\K_e$, and a localization $\K_e[Gr \inv]$,   allowing division by
 certain classes of  group varieties.
  in \S 4 we recall
the motivic Fourier transform, in the very simple context of test functions that we need.  

  In \S 5, we define motivic global test functions,
and the motivic analogue $\CR$ of the functional summing a test function over rational points of a function field.   \S 6 sets up a first order context useful for ``everywhere-local" definitions.
In particular we will be interested in {\em integral conjugacy classes} in a division algebra 
$D$ over a function field $\f(t)$.  For each place $v$ we define a subring $R_v$
of $D$; we say two elements are locally integrally conjugate at $v$ if they are
conjugate by some element of $R_v^*$ (possibly after base change), and 
integrally conjugate if they are locally integrally conjugate at every place $v$.  Such sets are conveniently defined in our setting; their
$\f(t)$-rational points forms a constructible set over $\f$.

In \S 8 we compute explicitly the constructible set of 
rational points in an integral conjugacy class.  We restrict attention to division algebras of prime degree $n$.   Let $c \in D$ and let $E$ be the subfield of $D$ consisting
of elements commuting with $c$.  Then $E=\f(C)$ for some curve $C$.
Let $O_c$ be the integral conjugacy class of $c$.  The 
 value of $\CR(O_c)$ closely related to  a certain Rosenlicht generalized Jacobian of the curve $C$, with ramification data  connected with the integral structure.  

More precisely, the value we obtain is connected with (a generalized) $Pic^0(C)$.  Now 
when $C(\f)$ has no rational point, the functor $Pic^0(C)$ is not represented by the Jacobian; in fact the functor $\k \mapsto Pic^0(C \times_\f \k)$ is not representable
by a variety at all (\cite{lichtenbaum}).  Nevetheless after some discussion of 
quotients in 2.2, we manage to associate
to $Pic^0(C)$ a class in a Grothendieck ring, treating it directly as an adelic
quotient $T(\Oo) \bs T(\Aa) / T(\k(t))$.  The class becomes equal to the class of the Jacobian over any field extension with a point of $C$.
At all events, $\CR(O_c)$ is expressed in terms mentioning only a commutative subalgebra of $D$.
 
We test our method on a problem involving division rings.  We consider
two division algebras $D, \dot{D}$  over $\f(t)$ associated with two distinct elements
of a given cyclic Galois group of prime order $n$ over $\f$.  (\S 7).  
We work with a quotient $\K$ of $\K_e[Gr \inv]$ appropriate for studying $D$ or $\dot{D}$ as a division algebra; namely, we factor out the class
$[\e_L]$ of a certain zero-dimensional variety $\e_L$ that has 
no rational points in any field $\f'$ such that $D$ is a division ring over $\f'$.  
 Consider
local motivic test functions $\phi$ on $D$ over $\f((t))$ that are invariant
under conjugation by $D^*$.  We will explain below how to match such
functions on $D^*$ with their homologues on  $\dot{D}$ (\defref{match}); the conjugacy classes can be identified.  

\<{thm} \lbl{A}   Let $\phi,\dot{\phi}$ be   matching
 local motivic test functions   on $D, \dot{D}$ over $\f((t))$,
with values in $\K_e[Gr \inv]$.   Then their 
Fourier tranforms $\fF \phi, \fF \dot{\phi}$   also   match. \>{thm}

This is closely related (when specialized to finite $\f={\mathbb F}_q$ and to numerically valued test functions) to   results of \cite{DKV}.
In \cite{DKV} the irreducible representations of each division algebra are shown to correspond to certain irreducible representations of $GL_n$; as a consequence
they correspond bijectively between two cyclic division algebras of the same dimension. 
Equivalently, the multiplicative convolution algebras of conjugation invariant
test functions on the two algebras are canonically isomorphic.  
   Such results
are purely local, but appear to be very difficult to prove using local methods in high dimensions (in degree two this is done by Jacquet-Langlands.)  
See Appendix 3 for a relation between the multiplicative and additive convolution 
algebras.

The proof of \thmref{A} proceeds by expressing the local Fourier transform of $\phi$ at the place $0$
in terms of the local Fourier transform at the place $1$, and some global terms, measuring rational
points on integral conjugacy classes 
(\secref{express}).    
  The global term was shown to depend only
on commutative subalgebras of $D,\dot{D}$; these are canonically isomorphic. 
The matching of Fourier transforms over $\f((t))$ is thus reduced to the same question over $\f((t-1))$,
where it is evident, since   over $\f((t-1))$ the two algebras are isomorphic.

\>{section}

\<{section}{Some  Grothendieck ring operations}

Here $T$ may be any theory.   We say ``definable" for ``definable by a quantifier-free formula
in the language of $T$."   This 
shorthand is acceptable notationally since our main application is to $T=ACF$, a theory with quantifier-elimination.
In reality the quantifier elimination provides little help, since it is not reflected in the Grothendieck ring,
and much   of our effort is directed at staying with quantifier-free formulas; see \remref{proof1}.  Similarly we will assume that any definable function $f$ is
piecewise given by terms.
 
 Let $T_\forall$ be
the universal part of $T$; so that $A \models T_\forall$ iff $A$ embeds into
some $M \models T$.  

  Let $\K(T)$ be the Grothendieck ring of definable sets, and let 
$\K$ be any $\K(T)$-algebra.  

 For any $A   \models T_\forall$, we have a Grothendieck ring $\K_A:= \K(T_A)$.
 So the Grothendieck ring is really a functor from models of $T_\forall$
 and embeddings among them, to the category of rings.  If necessary
 we will denote the class of a definable set $X$ in $\K_A$ by $[X]_A$.
 But usually we  write $[X]$ for $[X]_A$, and write $[X]=[Y] \in \K_A$
 to express: $[X]_A = [Y]_A$.  We also write $\K_b$ for $\K_A$ if $b$ generates $A$.
 
 This point of view
 is essential in discussing definable functions into $\K$ (cf. \cite{HK}).
 Let $f: X \to Y$ be a definable function, $X,Y$ definable sets.  
 We view the map
 $$ y \mapsto [f \inv(y)]$$
 as a function on $Y$, and let $Fn(X,\K)$ denote the family of all such functions.
 But it must be interpreted as follows:   for any $A \models T_\forall$, we obtain a function
 $Y(A) \to \K_A$, namely $b \mapsto [f \inv(b)]$.   
 
  We  have the   presheaf-like property:
  
\<{statement} \lbl{presheaf}  If $[X_b]=[X'_b] \in K_b$ for any $b \in Y$,  then $[X]=[X'] \in \K$.  \>{statement} 
 
See \cite{HK}. 

Since we will consider various localizations and quotients of $\K(T)$, it will be useful
to discuss $\K(T)$-algebras in general.

Consider   functors  $A \mapsto R_A$ from models of $T_\forall$, together with
natural transformations $\K(T) \to R$.   Given $A \models T_\forall$
and a $T_A$-definable set $X$, 
  the image of $[X]_A$ in $R_A$ is denoted $[X]_A^R$.
We assume that  \ref{presheaf} holds for $R$, i.e. if $f: X \to Y$ is
$T_A$-definable, and if 
 if $[X_b]^R_b=[X'_b]^R_b$ for any $b \in Y$,  then $[X]_A^R=[X']_A^R$.  
 Call such functors Grothendieck algebras over $T$.

\ssec{Localizing by a definable family}  \lbl{localize}

Let $\mathcal{N}$ be an Ind-definable family of definable sets
  Assume $\mathcal{N}$ is closed under products.  Then
$\{[X]: X \in \mathcal{N} \}$ is a  a multiplicative subset of the Grothendieck ring.  
 
Let $\K$ be a Grothendieck algebra for $T$.
To define the localization of $\K$ by $\mathcal{N}$, consider the family of  
all Grothendieck algebras $R$,  such that if $A \models T_\forall$
and $X \in \mathcal{N}(A)$ then $[X]_A^R$ is invertible in $R_A$.  We let $\K[\mathcal{N} \inv]$ be the universal object of ${\mathcal C}$.  

It is natural to assume that each set in $\mathcal{N}$ has a definable element.
In this case, in any finite structure $\f$, the number of points of $X \in \mathcal{N}$ is positive; so any zeta function defined on $\K$ factors through the localization 
$\K[\mathcal{N} \inv]$.  In the application the elements of $\mathcal{N}$ will have
a  distinguished element $1$, and indeed will be essentially classes of  definable groups.

In practice we will only use the following consequence of the existence of inverses:

\<{statement} \lbl{211}
{ \em Let $(A_y: y \in Y) \subseteq \mathcal{N}$, let $(X_y),(X'_y)$ be two families of definable sets, and assume:
  $$[A_y] [ X_y] = [X_y] ^2 , \ \ \ [A_y] [ X'_y] = [X'_y] ^2$$
  $$[X_y] ^2= [X_y] [X'_y] = [X'_y]^2$$ for $y \in Y$.  Then }
$$ {\sum_{y \in Y} [X_y] = \sum_{y \in Y} [X'_y]. }$$
 \>{statement}

Here is a proof that the relation holds if
 division by $[A_y]$ is possible, i.e.   elements $e_y = [X_y] / [A_y], e_y' = [X'_y]  [A_y]$ exist with $e_y [A_y] = [X_y], e_y' [A_y] = [X_y']$.      
  Then 
$e_y= e_y^2=  e_ye_y' = (e'_y)^2=e'_y$.  So 
$$\sum_{y \in Y} [X_y] = \sum_{y \in Y} e_y [A_y] =  \sum_{y \in Y} e_y' [A_y] = \sum_{y \in Y} [X'_y]$$
 
\<{remark} When $A$ is a group,  
the relation $[A][X] = [X]^2$ is typical of principal homogeneous spaces $X$. 
For two torsors $X,X'$, the relation $[X][X'] = [X]^2  = [X']^2$ holds if $X,X'$
generate the same subgroup  of the Galois cohomology group, so that over
any field extension, one represents the zero class iff the other does.  \>{remark}

Note that
localization (even by families) commutes with quotients.  Also note that
$\K[\mathcal{N} \inv] = \K[\mathcal{N} \inv][G\mathcal{N} \inv]$ canonically.

\ssec{Representable quotients.}

Recall  the notion of ``piecewise definable''  or Ind-definable  from Appendix 1.   
In terms of saturated models $M$, an Ind-definable set is a   union $\union_{i \in I} X_i(M)$
of definable sets $X_i(M)$; where $I$ is an index set, small compared to $M$.

Let $E$ be an Ind-definable equivalence relation on an Ind- definable set $V$

Define a {\em weakly representative set} for $(V ,E)$ to be a definable set $Y$
such that for some definable $X$ and $f: X \to V$ and surjective $g: X \to Y$, every element of $V$ is $E$-equivalent to some element $f(x)$, and $g(x) = g(x')$ iff $f(x) E f(x')$.  

We require that $g$ is surjective only in the geometric sense, i.e. in models of $T$. 
For $b \in Y$ we let $V_b$ be the $E$-equivalence class of $f(a)$, for any
$a$ with $g(a)=b$.
 
If $Y,Y'$  are two weakly  representative sets for $(V, E)$, then $Y,Y'$ are definably isomorphic; moreover the isomorphism $y \mapsto y'$ is such that $V_y = V_{y'}$.    
 We write $Y=V/E$. 

If $V$ is a definable group   and $E$ is the equivalence relation corresponding
to a normal subgroup, then $Y$ carries a definable group structure.

\<{lem} \lbl{wrep}   $(V,E)$ is weakly representable iff  
  \<{enumerate}
\item  There exists a definable $X$ and $f: X \to V$   such that any element of
$V$ is $E$-equivalent to some element of $f(X)$. 
\item  For any such $f,X$, the equivalence relation $f^{-1} E$ is a definable relation on $X$.
\>{enumerate}
\>{lem}
\prf Clear. \eprf

We say that $(V,E)$ admits a {\em set of unique representatives} if above one can
choose $X=Y, g =Id_X$.  Only in this case can we be sure of a bijection $V(\f)/E \to Y(\f)$.

More generally, let $S$ be a subset of $\K$, closed under multiplication.  We say that
$(V,E)$ is $S$-representable  
 if it   is weakly representable by $(X,Y,f,g)$ as above, and for some definable set  $Z$ with $[Z] \in S$,  whenever $(U,Y,f_1,g_1)$ is another weak representation of $(V,E)$, $X_b= g\inv(b)$, $U_b=g_1 \inv(b)$,
 we have: 
 $[U_b ] [Z] = [U_b] [X_b] \in \K(T_b)$ for any $b \in Y$.   \footnote{The   weaker statement  $[U \times Z] = [U \times_Y X]$ appears to  suffice for our purposes.}
 
The case we have in mind is with $Z$ a definable group, and $X_b$
a $Z$-torsor.  Assume $X_b$ becomes trivial over any point of $U_b$.
Then we have an isomorphism $U_b \times Z \to U_b \times X_b$, over $U_b$.

We assign to $(V,E) $ the class $[X] / [Z]$ in the localization $\K[S \inv]$, and denote it  by $[V:E]$.   If $E$ is the orbit equivalence relation of the action of a group $H$ on $V$, we also write $[V:H]$.

  If $(V,E) $ is  $S$-representable,  then in particular 
$(V,E) $ is weakly represented, so $V/E$ is defined as well as $[V:E]$.  But the image of $[V/E]$ in $\K[S \inv]$ is not necessarily  equal to $[V:E]$.

\<{lem} \lbl{rep1}   $[V:E]$ is well-defined.  I.e. if  $(X',Y',Z',f',g',h')$ is another $S$-representation,  then $[X]/[Z] = [X']/[Z']$ in $\K[S \inv]$.   

If $V=V_1 \dot{\union} V_2, E_i = E | V_i$, and $(V_i,E_i)$ is $S$-representable, then so is $(V,E)$,
and $[V:E] = [V_1:E_1]+ [V_2:E_2]$.  \>{lem}

\prf  
$ [X' \times  Z ] = [X' \times_Y X  ] = [X \times_Y X'] =  [X  \times Z']$
so dividing by $[Z \times Z']$ we find:  $[X']/ [Z']    = [X]/[Z]  $.

 If $(V,E) $ is $S$-representable via $(X,Y,f,g,Z)$, then for
 any $Z' \in S$ it is also $S$-representable 
 via $(X \times Z', f, g, Z \times Z')$. Hence by taking common denominators we may  assume   $(V_i,E_i)$ is $S$-representable via the same denominator $Z$. 
 In this case 
the statement on addition is immediate.

 \eprf

If $S \subset \K$ is not closed under multiplication, we let $<S>$ be the set of products of elements
of $S$, and define $S$-representable to mean $<S>$-representable.

\<{lem} \lbl{srep} Let $Z$ be a definable group. 
Assume $(V,E) $ is weakly representable, and
let $Y =  V/E$.  Assume:     for any $b \in Y$,
there exists an $\f(b)$-definable $Z$-torsor $R_b$ such that there exists a definable function
$f_b: R_b \to   V_b$, and for any $\f' $ with
$\f(b) \leq \f'$ and $V_b(\f')  \neq \emptyset$, we have $R_b(\f') \neq \emptyset$.  
 
Then $(V,E)$ is   $[Z]$-representable.  
 
\>{lem}

\prf   Say $(V,E) $ is weakly representable via  $(X,f,g,Y)$.  Let $X_b = g\inv(b)$. By assumption,
for any $b \in Y$  there exists an $\f(b)$-definable $Z$-torsor 
$R_b$, $f_b: R_b \to V_b$, and an $\f(b)$-definable map $h_b: X_b \to R_b$.
By compactness and glueing, we can take $(R_b,f_b,h_b)$ to be uniformly definable.

Let  $(U,Y,f_1,g_1)$ be another weak representation of $(V,E) $, $X_b= g\inv(b)$, $U_b=g_1 \inv(b)$.    Then  $[U_b ] [Z] = [U_b] [X_b] \in \K(T_b)$.  Indeed
for any $c \in U_b$, using the point $h_b(c)$ we find a bijection $j_b: Z \to X_b$;
these can be glued to give a bijection $U_b \times Z \to U_b \times X_b$,
over the identity on $U_b$. 
\eprf  
 
 Let  $V$ is Ind-definable,   $E$ an Ind-definable equivalence relation on $V$, and $X$ a definable set.    By a definable function $X \to V/E$ we mean 
is definable relation $F \subset X \times V$ whose projection  $p$ to $X$ is surjective
and  1-1 modulo $E$ (i.e. if $(x,v),(x,v') \in F$ then $(v,v') \in E$.)   In this
case if $x \in X$, we let $f(x)$ be the $E$-class of $v$ for any $v$ with $(x,v) \in F$.

We say that $E$ is definable-in-definable families if
for any definable $U \subset V$, the restriction of $E$ to $U$ is definable. 
In this case the quotient $V/E$ is Ind-definable.

 \<{remark}  \rm  Let $f: X \to Y$ be a definable map between
definable sets,.  Let $E$ be the equivalence relation $f(x)=f(y)$.  If
$[X:E]$ is defined, write $f_*[X]=[X:E]$.  
For a given structure $A$ we may be interested in $f(X(A))$,
but within the Grothendieck ring of quantifier-free formulas, or of formulas up to $T$-equivalence when $A$ is not a model of $T$, we have no direct way to describe it.    
The class $f_*[X]$, when defined, offers a substitute.    
 
\>{remark}

\ssec{Essentially representable sets}  \lbl{T*}

Let $T^*$ be a universal theory containing the set $T_\forall$ of universal consequences
 of $T$.  An  Ind-definable set (of $T$) is called {\em formally empty} if it is the union of  definable sets $U$ such that $T^* \models U  = \emptyset$.   
 A structure $\f'$ is said to be {\em negligible} if 
$W(\f') \neq \emptyset$ for some formally empty $W$, or equivalently if
$\f' \not \models T^*$.  
 
 Let $\K^* = \K / I$, where $I = \{[X] \in \K:  T^* \models X = \emptyset \}$.

 Let $X$ be an Ind-definable set; write $X = \union_i X_i$ with $X_i$ definable.
We say that $X$ is {\em $T^*$-limited} 
if for some finite $I_0 \subseteq I$, letting $X_0 = \union_{i \in I_0} X_i$,
we have for all $j$, $T^* \models X_j \subseteq X_0$.    We   write $T^* \models X = X_0$ for short.   In this case let   $[X]$ be the image
of $X_0$ in $\K^*$.  The definition clearly does not depend on the choice of representation $ \union_i X_i$ or
on the finite set $I_0$.


 We say that $V/E$ is  {\em  $(T^*,S)$-representable} if there exists an Ind-definable
$V' \subseteq V$ and a   formally empty Ind-definable $V''$ with $V = V' \union V''$,
 and such that $V'/E$ is $S$-representable.  The image of   $[V':E]$ in $\K^*[S \inv]$
 is then well-defined, and denoted $[V:E]$.

A further variant of \lemref{srep} will be useful.  

 \<{lem} \lbl{srep3} Let $E$ an Ind-definable equivalence relation 
 on the Ind-definable set $V$.       Let  $Z$ be a definable group acting freely on a $T^*$-limited set $X$,  and  $f: X \to V/E$  an Ind-definable function whose  fibers    are   $Z$-orbits.   
  Assume:   for any  $\f' \geq \f$ and any $\f'$-definable $Z$-orbit $U$,
 $\union_{c \in X(\f')} f(c) \subseteq \union_{c \in V(\f')} cE$,
 where $cE$ is the $E$-class of $c$; with equality if  $\f' \models T^*$.
   Then $V/E$ is   $(T^*,Z)$-representable, and $[V:E] = [X]/[Z]$.

 \>{lem}
 
\prf   We may express $X$ as a direct limit of definable sets $X_i$.  We have $X_i \subseteq ZX_i \subseteq X_{i'}$ for some $i' \geq i$.  Replacing $X_i$ by $ZX_i$, we may assume the $X_i$ are $Z$-sets.  Since $X$ is $T^*$-limited, for some $X_0$ we have $X_i \m X_0$ formally empty, for all $i \geq 0$.    
 Let $f_0 = f | X_0$. Then for any  $\f' \geq \f$ and any $\f'$-definable $Z$-orbit $U \subseteq X_0$,
 $\union_{c \in X(\f')} f_0(c) \subseteq \union_{c \in V(\f')} cE$. Moreover if $\f' \models T^*$, $c \in V(\f')$
 and  $y \in cE$, then $y \in f(d)$ for some $d \in X(\f')$.  Since $X_i \m X_0$ is formally empty and
 $\f' \models T^*$, we have $d \notin (X_i \m X_0)$ so $d \in X_0(\f')$.  This shows that the hypotheses
 hold for $X_0,f_0$.  Since $[X]=[X_0]$ by definition, we are reduced to this case.  We may thus assume that
 $X$ is a definable set.

 For any $c \in X$,   there exists $b \in \f(c)$ 
such that    $b \in f(c)$.   By compactness, there exists a definable function 
 $g: X \to V  $ with $g(c) \in f(c)$.  
 
We have $g(x)Eg(x')$ iff $f(x)=f(x')$ iff $x,x' \in X'$ are $Z$-conjugate.  Let
$V'=g(X)$; then $V'$ is a definable subset of $V$, $V'/E$ is weakly representable, and $V'/E \cong X/Z$.  
  
For   $b \in V$, if $\f(b) \models T^*$ then    $b \in f(c)$ for some $c \in X(\f(b))$, so that
$b E g(c)$.   By  compactness there exists $V'' \subseteq V$ such that $v E g(v)$ for $v \in V''$, and $V \m V''$ is formally empty.

 The proof is now completed as in  the 2nd paragraph of 
 the proof of  \lemref{srep}.

 \eprf 
 
\<{rem} \lbl{srep3r}  \rm
If, in \lemref{srep3}, $Z$ does not act freely, but the stabilizer of each point in $X$ is one of 
finitely many groups $H_i \leq Z$, $Z_i  = Z/H_i$, then $V/E$ is essentially $\Pi_{i=1}^k Z_i$-representable,   and  $[V:E] = \sum_{i=1}^k [X_i]/[Z_i]$,
where $X_i$ is the union of the $f$-classes that are $Z_i$-orbits.   \>{rem}

\ssec{Absolute elements and invariant functions}  
A set $S$ of elements of $\K$ is called {\em absolute} if whenever $s,s' \in S$
are defined and distinct over $\f'$, and $\f' \leq \f''$, then   the images of $s,s'$
in $\K_{\f''}$ remain distinct. 
A function $\phi: E \to \K$ is said to take  absolute values on $E$ if $\{\phi(e): e \in E \}$ is absolute.   For instance, $\{0,1,[\k],[\k^2],\ldots\}$ forms an absolute set in $\K(ACF)$.

Let $G$ be a definable group, with a definable action on a definable set $D$, fibered
over $Ob G$.  Let $\phi: D \to \K$ be a definable function. 

\<{defn}\lbl{invariant-gen}  $\phi$ is {\em  $G$-invariant}
if for $c \in D, g \in G$ we have $\phi(c)=\phi(gc) \in \K_{c,g}$.  
 We say  $\phi$ is {\em  strongly $G$-invariant} if for any such $c,g$ 
 we have $\phi(c)=\phi(gc) \in \K_{c,gc}$.   \>{defn}
 
 The same definition could be made for a groupoid $G$.   In this
 case for each object $a \in Ob G$ we have a set $D_a$, and for each pair
 $a,b \in Ob G$ we have a definable function $Mor_G(a,b) \times D_a \to D_b$,
 such that the obvious associativity relations hold.   Given a family $(\phi_a)$
 of functions $D_a \to \K$, we have a definition of (strong) invariance as above.
 We will only use the case of
 two objects, with two corresponding definable division algebras $D, \dot{D}$; we will
 have $G=Aut(D)$as a division algebra, $\dot{G} = Aut(\dot{D})$, and
 also $M=Iso(D,D')$.  A pair $(\phi,\dot{\phi})$ of functions on $D,\dot{D}$ 
 is then said to be (strongly) {\em matching} if it is invariant under the groupoid.

\<{lem} \lbl{absolute} Let  $\phi: D \to \K$ be a definable function taking  absolute
 values.  If $\phi$ is invariant, then it is strongly invariant.  \>{lem}
 
 \prf  Clear. \eprf

While $G$-invariance depends on the group or groupoid action, the notion of strong $G$-invariance depends only on the equivalence relation $conj_G$ of $G$-conjugacy on $D$.
For an equivalence relation $E$ on $D$, 
say $\phi: D \to \K$ is $E$-invariant if $\phi(c)=\phi(c') \in \K_{c,c'}$ whenever $(c,c') \in E$.  Then strong $G$-invariance is the same as $conj_G$-invariance.
A still stronger notion is {\em descent to $D/E$}, i.e. existence of a definable $\bar{\phi}: D/E \to \K$ with $\phi(d) = \bar{\phi}(d/E) \in \K_d$.  We have however:

\<{lem}\lbl{invariance-descent}  Let $E$ be an equivalence relation on $D$.  Assume $\phi: D \to \K$ is
$E$-invariant.   Let ${\mathcal F}$ be the set of equivalence classes of $E$, 
and let $\widetilde{\K} = \K[{\mathcal F} \inv]$ be the localization, \secref{localization}.  Let
$\widetilde{\phi}: D \to \K'$ be the induced map.  Then $\widetilde{\phi}$ descends
to $D/E$.    \>{lem}

\prf  For an equivalence class $y$ of $E$, define $\bar{\phi} (y) = [y] \inv \sum_{d \in y}\phi(d)$.  If $d \in y$, then by $E$- invariance we have $\phi(d') = \phi(d) \in \K_{d,d'}$
for any $d' \in y$, so $\sum_{d' \in y} \phi(d') = [y] \phi(d) \in \K_d$.  It follows that
$\phi(d) = \bar{\phi}(y) \in \widetilde{\K}_d$ for any $d \in y$. \eprf

\ssec{Proof by cases.}   \lbl{lf}

 Let $T$ be a theory of fields, with base field $\f$; let $\lf$ be a finite Galois
 extension of this base field; so $\lf=L(\f)$ for some commutative definable
 algebra $L$.    Note that for any field extension $\f'$ of $\f$, $\f'$ contains
 a copy of $\lf$ over $\f$ iff $L(\f')$ is not a field.  At all events $L(\f)$
 has points - it is a $\dim(L)$-dimensional extension of $\f$ - and should not be 
 convused with the isomorphic copy $\lf$.
 
 We explain how an identity in the semi-group can be proved
 by cases, according to whether $L$ splits or not  in extension fields of $\f$.
 Similar considerations apply to definable finite sets in any theory.
   
 Let $I_L$ be a normal basis for $\lf$ over $\f$.  Then $I_L$ can be
 viewed as a finite definable set, so it has a class $[I_L] \in \K(T)$.
   If $I'_\lf$ is another  normal basis,
 there exists a bijection $f: I_L \to I'_L$ left invariant by the Galois action,
 and hence definable; so $[I_L] = [I'_L]$.  Thus the class $[I_L]$
 depends only on $L$.  Moreover $\e_L := \frac{1}{n} [I_L]$
is an idempotent in $\K[\frac{1}{n}]$.  

   Let   $\K$ be any
$\K_+(T)$-semi-algebra, such that:  \<{enumerate}
\item   $n$ has a multiplicative inverse. 
\item  $\e_L =0$.
\end{enumerate} 

\<{lem} \lbl{gr1}  If $V(\f')=0$  for any $\f'$ with $L(\f')$ a field, then $[V]=0 \in \K$. \>{lem}

\prf  Assume the condition holds.  Then whenever $b \in V$, $L(\f(b))$ is not a field;
hence $L(\f(b)) \cong \f(b)^n$; 
so $I_L \cong n$ over $\f(b)$. (where $n$ denotes a set of $n$ definable points.)
Hence $ V \times I_L \cong V \times n$.  Dividing by $n$ we obtain
 $[V] = [V] \e_L = 0$.  \eprf

 Typically $\K$ wil be a  
$\K(ACF_\f)$-algebra, obtained as a localization of the Grothendieck ring $\K(ACF_\f)$
or $\K_e(ACF_\f)$, and factoring out $\e_L$.     

\<{lem}\lbl{gr2}  Let  $\K$ be any
$\K(ACF_\f)[n \inv]$- algebra.  Assume $[V_1]=[V_2]$ holds in $\K_{\lf}$ and also in $\K/\e_L$.   Then
$[V_1]=[V_2] $ in $\K$.  \>{lem}
\prf An element of $\K_{\lf}$ is represented by an   element $X$ of $K$ and a 
definable function $g: X \to I_L$.  We saw that $[X] = [X] [I_L]/n = [X] \e$,
i.e. $[X] \in \K (1-\e)$. It is easy to see that $\K_{\lf}$ is isomrphic to
the ring $ \K \e$ with unit $\e$, or equivalently to $ \K / (1-\e)$.   Now if a class
vanishes modulo $\e$ and modulo $(1-\e)$, then it vanishes.  \eprf

\<{remark} \lbl{proof1}  \rm Let $TF$ be the theory of perfect fields.  Let $\K(TF)$ be the Grothendieck ring of
all formulas,  including quantified ones, and including   imaginary sorts.  Let $\K_c(PF)$ be the quotient
of $\K(TF)$ obtained by imposing a  
Cavalieri principle:  if $f: X \to Y$ is a definable function, and each fiber is provably isomorphic to $Z$,
then $[X] = [Z]$.     \thmref{A} admits a considerably simpler proof if values are taken in this ring; to begin with,
all issues regarding representability of quotients become superfluous.   The extra effort required
in using the quantifier-free Grothendieck ring is hopefully paid off in   geometrically
more precise answers. \>{remark}

\>{section}

\<{section}{The Grothendieck ring of exponential sums} 
 
Let $T$ be a theory of fields.  In particular, $T$ includes constants for a subfield $F$.
It is possible to allow additional relations, 
but in the present paper we will not use them so one can take the language to be the language
of $F$-algebras.  The field sort is denoted $\k$.  We will assume that the models of $T$ 
are perfect fields.  
    
In some parts of the paper 
 we will assume the existence of division algebras over $\k(C)$, where $C$ is a curve
over $\k$; in particular $\k$ is not algebraically closed.  Nonetheless constructions that do not require this assumptions are better carried out geometrically.   Thus for instance we will
define the Fourier transform of a test function, and the ``sum over rational points" functional, over the theory ACF of algebraically closed fields.   We define below a natural homomorphism
from the Grothendieck ring of exponential sums over ACF, to the 
Grothendieck ring of exponential sums over $T$; any equations holding at the ACF level will
thus continue to hold.

We define the   Grothendieck ring of exponential sums using generators and relations.

The generators are elements $[X,h]$ where $X$ is a definable set, and $h: X \to k$ a
definable function.

We write $\psi(c)=[\{c\},Id]$; we think of $\psi$ as an additive character, 
and   think of $[X,h]$ as  representing $\sum_{x \in X} \psi(h(x))$.  We will impose
the following relations:

\eq{1}{  [X,h][Y,g] = [X \times Y, h(x)+g(y)]; \, [0,0]=1 }
If $X,Y$ are disjoint
\eq{2} {    [X,h] + [Y,g] = [X \union Y, h \union g]    }
if $g: X \to Y$ is a definable bijection, 
\eq{3} {     [Y,h] = [X,h \circ g]   }
\eq{4} {                  [\k,Id] = 0                    }

Let $K_{exp}(T)$ be the ring presented by the generators $[X,h]$ and relations
(\ref{1})- (\ref{4}).   Define $\Ll = [k,0]$.  $K_{exp}(T)$ is naturally filtered by dimension: 
 ${F}_d K_{exp}(T) = \{[X,h]: \dim(X) \leq d \}$.   
Let  $K_e(T) = \K_{exp}(T)[\Ll \inv]$.

\<{lem} \lbl{l4}    
Let $(u,x) \mapsto u+x$ be a definable action of $(k,+)$  
 on a definable set $X$, and let $h(t+x)=t+h(x)$.  Then $[X,h] = 0$ in $K_e(T)$.
\>{lem}
 
\prf  For any   $c \in k$, and $X,h$ as above, 
we have 
$$[X,h] \psi(c) = [X,h] [\{c\},Id] = [X \times \{c\}, h(x)+c] =[X,h]$$
The last equality uses \eqref{3}, for the bijection $X \times \{c\} \to X$ given by the action of $c$. 
Thus $[X,h] (\psi(c)-1) = 0$.   Summing over all $c$ and using $\sum_{c \in \k} \psi(c)= [\k,Id] = 0$ we obtain $[X,h] (0 - [\k]) = - \Ll [X,h]$. Since $\Ll$ is invertible in $K_e(T)$
the result follows.  \eprf
 
 It will sometimes be useful to consider the semiring with the same generators as
 $K_{exp}(T)$, relations \eqref{1}-\eqref{4}, and in addition the relations $[X,h]=0$ 
 for $[X,h]$ as in \ref{l4}.  We will refer to these relations as \ref{4}'.

\<{remark}  Let  $K_{exp}'(T)$ be the ring presented by the generators $[X,h]$ and relations
(\ref{1})- (\ref{3}).   Then (4) holds in $K_{exp}'(T) [ (\psi(c)-c) \inv ]$, for any
definable element $c$ of $\k$.   In particular if $\psi(1)-1$ is inverted, relation (4) need
not be explicitly imposed.  \>{remark}

\prf  This follows from the computation $[X,h] ( \psi(c)-1)=0$ of the previous lemma. \eprf

\<{lem} \lbl{rem1}  
\<{enumerate} 

\item  If $k \models ACF_F$ then the conclusion of \lemref{l4} follows from Equations \ref{1}-\ref{4}, even without inverting $\psi(1)-1$ or $\Ll$.

\item  Every element of $K_{exp}(T)$ is represented in the form $[X,h]$ for some $X,h$.
\>{enumerate}
\>{lem}

\prf

(1)  In this case, there exists (in the constructible category) a quotient $Y$ of $X/ (k,+)$, and 
 by Hilbert's theorem 90,  $X \to Y$ admits a constructible section $f$.  So $X $ is definably isomorphic as a $(k,+)$-set to $Y \times k$.Thus $[X,h] \cong [Y, h \circ f] \times [k,x]$,
 and \lemref{l4}
 follows form Equations \ref{1} and the assumption $[\k,Id]=0$.
 
 (2)  The element $-1$ is represented by $[G_m,Id]$.  Hence 
 $[X,h] - [X',h'] = [X,h] + [G_m,Id ] [X',h']$.
 \eprf

If $h: Y \to \k$ and $f: Y \to X$  are definable functions, and $n \in \Nn$, for $a \in X$ we obtain 
an element $ L ^{-n} [Y_a,h |Y_a] $ of $K_e(T_a)$. By a {\em definable function} $\theta: X \to K_{e}(T) $
we mean a function $a \mapsto  \theta(a) \in K_{exp}(T_a)[\Ll \inv]$ of this form.  
Note in particular that this is not literally a function into $K_e(T)$.  
  Given such a definable function $\theta$
on $X$, and $g: X \to k$, 
we write 
\eq{summ}{ \sum_{x \in X} \theta(x) \psi(g(x)) := [Y, h(y) + g(f(y))] }

For $g=0$ we obtain a definition of $\sum_{x \in X} \theta(x)$.  

Let $A=Fn(X,K_e(T))$ be the set of definable functions $X \to K_e(T)$.  Given $f,g \in A$ we can define $(f,g) = \sum_{x \in X} f(x)g(x)$.  

\<{remark} \lbl{functionals} \rm For $a$ generating a substructure of a model of $T$,
let $F_a=Fn(X,K_e(T_a))$ be the set of $T_a$-definable functions $X \to K_e(T_a)$.  If $g \in Fn(X,K_e(T))$ we obtain, for each $a$, a homomorphism $\chi_g: F_a \to K_e(T_a)$, namely $\chi_g( f) = \sum_{x \in X} g(x)f(x)$.  Then if $f=f_t$ varies uniformly
in some definable family, $\chi_g(f_t)$ is a definable function of $t$.   It satisfies:

(*)  $\chi_g (\sum_{y \in Y} \theta(y) f_y) = \sum_{y \in Y}\theta(y) \chi_g( f_y) $
where $(f_y: y \in Y)$ is a definable family of definable functions into $K_e(T)$, and $\theta \in Fn(Y,K_e(T))$.  

(**) If $a =h(b)$ for a definable function $h$,  we have a natural homomorphism
$h^*: K_e(T_a) \to K_e(T_b)$, and  by composition also $h^*: F_a \ F_b$ Then
$\chi_g \circ h^* = h^* \circ \chi_g$.

Conversely, the $\chi_g$ are the only systerm of homomorphisms $F_a \to K_e(T_a)$,  given uniformly in 
$a$, satisfying the above properties.  For given such a system $(\chi)$, define $g(a) = \chi(1_{\{a\}})$
where $1_{\{a\}}$ is the characteristic function of the element $a \in X$.   Then from (*)
it follows that $\chi (f) = \chi ( \sum_{a \in X} f(a) 1_{\{a\}}) = \sum_{a \in X} f(a) g(a) = \chi_g(f)$.

Note that (*) is an analogue of $K_e(T)$-linearity, with finite additivity replaced by ``motivically finite" additivity.  In this sense the pairing $(f,g)$  
may be viewed as an isomorphism between $A$ and its ``motivic dual".

\>{remark}
 
\ssec{Polynomial maps on semigroups}

We  prove a general lemma that will be used to  extend the norm map, defined below, from the semiring to the ring $\Ke(T)$.

Let $A$ be a commutative semi-group,  $B$ be an Abelian group, and  $f: A \to B$ be a function.  We say that $f$ is 
    a {\em polynomial map of degree $0$} if $f$ is constant.  We say that $f$ is
          {\em polynomial   of degree $d+1$} if $\phi: A^2 \to B $ is polynomial of degree $d$,
          where $\phi(x,z) = f(x+z) - f(x) - f(z)$.  
          
 \<{lem} \lbl{polar}   Let $f: A \to B$ be a polynomial map of any degree $d$.  Let $a,b,c \in A$
 and assume $a+c=b+c$.   Then $f(a)=f(b)$.  \>{lem}
 
\prf  For $d=0$ this is clear.  Assume it is true for polynomial maps of degree $d$,
and $f$ has degree $d+1$.   
Let $\phi: A^2 \to B$, $\phi(x,z) = f(x+z)-f(x)-f(z)$.  Then $\phi$ is polynomial of degree $d$.
We have $(a,c)+(c,c)=(b,c)+(c,c)$ in the semigroup $A^2$.  Hence 
so $\phi(a,c) = \phi(b,c)$.  In other words $f(a+c)-f(a)-f(c) = f(b+c)-f(b)-f(c)$.
Since $a+c=b+c$, subtracting equal terms from this expression in the group $B$,
we obtain $f(b)=f(a)$.  \eprf

\ssec{Irreducible definable sets}

 A definable set $X$ is {\em irreducible} if $X \neq \emptyset$ and $X$ contains no proper,
nonempty definable sets.  (Model theoretically one says that $X$ {\em isolates a complete type.})   In $ACF_F$, of course, all irreducible sets are finite.

\ssec{Hilbert 90}

We will say  that {\em Hilbert 90 holds for $T$} if
any $T_F$-interpretable finite dimensional $\k$-vector space admits a basis
consisting of $F$-definable elements.  This is true in $ACF_F$ for any field $F$,
field, by the vanishing of 
  the first Galois cohomology of $GL_m(\k)$, see \cite{serre-loc} ch. 10, Prop. 3.  
On the other hand if  Hilbert 90 holds for $T$, then any definable, finite dimensional $\k$-algebra,
as well as any definable finite dimensional module over such an algebra, are already $T_F$-definable,
where $F$ is the field of definable points.  .

\def\uR{\underline{R}} \def\uA{\underline{A}} \def\uB{\underline{B}}
\<{lem}\lbl{90}  Let $T$ be a theory of   fields, and assume Hilbert 90 holds for $T$. Let $\uR$ be a definable  finite dimensional algebra (associative, with $1$). 
Then any definable $\uR^*$-torsor
has a definable point.  \>{lem}

\prf

If $\uA$ is a definable $\uR$-module,  let $\uA^*=\{a \in \uA: ra=0 \rightarrow r=0 \}$.  We say
$\uA$ is {\em free on one generator} if this is the case for $\uA(M)$, for some $M \models T$. 
If  $\uA$ is  free on one generator, then  $\uA^*$ is
an $\uR^*$-torsor.   Conversely, if $\uB$ is an $\uR^*$-torsor, let $\uA$ be the quotient of 
  $\uR \times \uB$ by the action of $\uR^*$, $(x,y) \mapsto (xr,r \inv y)$.  Then it is easy to define 
  an $\uR$-module structure on $\uA$, making $\uA$ into a form of $\uR$,with $\uA^*=\uB$.
We thus have to show that if $\uA$ is a definable $\uR$-module and $\uA$ is free on one generator $c$, then this generator can be chosen definable.

By the remark just above the lemma, $\uR$ and $\uA$ are $ACF_F$-definable.  The  non-generators of $\uA$ form
a proper Zariski closed subset of $\uA$.    When $F$ is infinite, $\uA(F)$ is Zariski dense in $F$, and it suffices
to choose any  generator in $\uA(F)$.   When $F$ is finite, we return to the connected algebraic group $R^*$, and use Lang's theorem instead.   \eprf
  
Note also that if Hilbert 90 holds for $T$, then   any definable $(\k,+)$-torsor $H$
has a definable point, since $H$ can be viewed as an affine line within a definable 2-dimensional $\k$-space.

\ssec{Norm map}

In this paragraph 
let $X$ be a {\em finite} definable set, and $Def_X$ the category of definable sets over $X$, 
i.e. definable sets $Y$ together with definable maps $Y \to X$; a morphism is a definable map $Y \to Y'$
commuting with the maps to $X$.  

We have a functor $N: Def_X \to Def$,with $NY= \prod_{x \in X} Y_x$ the set of sections of $Y  \to X$.  

In case $T=ACF_F$ and $X$ is irreducible,   this is just 
the Weil restriction of scalars of $Y_a$ from $F(a)$ to $F$ (where $a \in X$.)

We will assume now that  Hilbert 90 holds for $T$.  
Consider a triple  $(Y,f,h)$, with $(Y,f) \in Def_X$ and $h: Y \to \k$ a definable function.  
For $z \in NY$ define $Nh(z) = \sum_{x \in X} h(z(x))$; this is a finite sum taken in $(\k,+)$.  Let $N(Y,f,h) = [NY, {Nh}]$. 
 Let $S=Fn(X,K_{exp}(T))$.
  We wish to show that $N$ induces
  a map $N: S \to K_{exp}(T)$.  In other words, we need:
  
  \<{lem}  Assume Hilbert 90 holds for $T$.  Let  $(Y,f,h)$ represent  a definable function $\Phi: X \to K_{exp}(T)$.  Then $N(Y,f,h)$
depends only on $\Phi$.    \>{lem}

\prf

First we let $S^+$ be the semiring of definable sections $X \to K_{3}^+(T)$,
where $K_{3}^+$ is the semiring with generators and relations of  \eqref{3} over $X$,
with operations defined by Equations \ref{1},\ref{2}.    It is easy to see from functoriality that \eqref{3} is respected by $N$, i.e. the norm of definably isomorphic
pairs over $X$ are definably isomorphic.  We  
obtain a map $N: S^+ \to K_{3}(T)$.

We now verify that relation \eqref{4} is respected by $N$.  In fact we will
do more, and verify that the relations \eqref{4}' in   \lemref{l4} are respected too:   
 assume  $\k$ acts definably on each fiber $(X_a, h |X_a)$ of $Y \to X$,   via a 
map $\rho: \k \times Y \to Y$, with $f \circ \rho = f$ and $h(\rho(t,y))= t + h(y)$.
We will show that  if $[Y,h]+[W,g]=[W',g']$ in $S^+$, then $[NW,Ng] = [NW',Ng']$.
To begin with we show that $[NY,Nh]=0$.  
Let $N\k$ be the group of   maps $X \to \k$, and for $\bar{t} \in N\k$
let $\phi \bar{t}  = \sum_{x \in X} \bar{t}(x)$.  
We obtain by functoriality an action 
$N \rho:  N\k \times NY \to NY$, with 
\eq{7}{ 
Nh ( N \rho (\bar{t},\bar{y})) = \phi \bar{t} + Nh(\bar{y})}
   The $\k$-vector space structure on $\k$ is inherited
by $N\k$, and $\phi: N\k \to \k$ is a surjective $\k$-linear transformation.
The kernel $\ker \phi$ is a $\k$-vector space;  so
there exists a   basis of $\ker \phi$ consisting of definable elements.  Now using Hilbert 90, the torsor $\phi \inv (1)$ has 
  a definable element $\bar{t_1}$.  Define an action of $\k$ on $N\k$
  by $(\a,\bar{y}) \mapsto N \rho (\a \bar{t_1}, \bar{y})$.  Then
  the conditions of \lemref{l4}
 are met, so $[NY,Nh] = 0$ in $K_{exp}(T)$.

Next we show, for $Y$ as in the previous paragraph, not only that $[NY,Nh]=0$
but also that  if $[Y,h]+[W,g]=[W',g']$ in $S^+$, then $[NW,Ng] = [NW',Ng']$.  We can take $W'$  to be the disjoint union of $Y$ and $W$, and extend the action of $\k$ on $Y$ to an action on $W'$,
trivial on $W$.  Let $N_{i}W' $ be the set of sections $s: X \to W'$, such that $|s \inv (Y)| = i$.
Then $NW'$ is the disjoint union of $NW$ and of $N_iW'$ for $i \geq 1$.  So it suffices to show that
$N_iW' = 0 \in K_{exp}(T)$ for $i \geq 1$.  

Let $[X]^i$ be the set of $i$-element subsets of $X$.   For $w \in [X]^i$, consider the 
$i$-dimensional $\k$-space $\k^w$, and define $\phi_w: \k^w \to \k$ by $\phi_w(x)=\sum_{t \in w} x(t) $.
Let $B$ be the fiber product of all spaces $\k^w$ over $\k$, via the maps $\phi_w$.  I.e.
$B = \{ (a_w) \in  \prod_{w \in [X]^i} \k^w:  (\exists \a \in \k)(\forall w \in [X]^i) \phi_w(a_w) = \a \}$.
This is a $\k$-space of dimension $({{|X|}\choose{i}})(i-1) +1$.  We have a linear map
$\phi: B \to \k$, $\phi ((a_w)) = \phi_w(a_w)$ (for any $w \in [X]^i$.)
 We have an action of
$B$ on $N_iW'$, as follows.  We have a map $\psi: N_iW' \to  [X]^i$, $\psi(s) = s \inv(Y)$.  
$B$ will preserve the fibers of $\psi$; on $\psi  \inv (w)$, $B$ will act via the $w$'th coordinate, 
i.e.
 $(a_w) + s =s'$ with $s'(t) = a_{s \inv(Y)}(t)+s(t)$ if $s(t) \in Y$,  $s'(t)=s(t)$ if $t \in W$.  
 So Equation \eqref{7} holds.  From here the proof is the same as in the   paragraph following equation \eqref{7}.  
 
It follows that   $N$ induces a function on the semiring $S_1^+$ with generators \ref{1}-\ref{4}'.    
It is easy to see that $N: S_1^+ \to K_{exp}(T)$ is polynomial of degree $|X|$.   By \lemref{polar}
it induces a function on the image of $S_1^+$ in  the associated ring $S$.  This function is 
clearly multiplicative, $N(ab)=N(a)N(b)$.  
By  \remref{rem1} this image equals
$S$, so we obtain $N: S \to K_{exp}(T)$.    \eprf
 
Inverting $\Ll \inv$, we find a norm map from the ring of definable maps $X \to  K_e(T)$, 
into $K_e(T)$.   We will also denote the norm of $c$ by $\prod_{x \in X} c(x)$.

Let $Sym_n(X) = Y/Sym(n)$,
where $Y$ is the set of distinct $n$-tuples of $X$. 
\footnote{$Sym_n(C)$ is a definable set in a possibly imaginary sort of $T$.}
 We also treat elements of $Sym_n(X)$
as $n$-element subsets of $X$.  

\<{defn} \lbl{products} Let $c: X \to \Ke(T)$ be a definable function.  Define
$$\prod_{x \in X} (1+c(x)t)  = \sum_{n=0}^\infty b_n t^n$$ where 
$b_n = \sum_{s \in Sym_n(X)} \prod_{t \in s} c(t) $  \>{defn}

By a definable
function $X \to \Ke(T)[[t]]$ we mean one of the form $c(x) =\sum c_n(x) t^n$ with each $c_n$ a definable function into $\Ke(T)$.   Define $\prod_{x \in X} (1+c(x) t)$ where $c: X \to \Ke(T)[[t]]$ is a definable function in the natural way.  We have:

\<{lem}  \lbl{products-trans}
If $f: Y \to X$ is definable and $b: Y \to \Ke(T)[[t]]$, define $a(x)$ by:
  $1+ta(x) = \prod_{y \in Y_x} (1+t b(y))$.  Then $\prod_{x \in X} (1+ta(x)) = \prod_{y \in Y} (1+tb(y))$. \>{lem}

\prf
  If we allow rational coefficients, this can   be done using the analogous statement for sums,
  which is evident, and the
isomorphism $\exp:  t \Ke(T)_\Qq[[t]] \to 1 + t \Ke(T)_\Qq[[t]]$.  Otherwise, a routine  computation by coefficients.  \eprf

\ssec{Compatibility(1)} 

Assume $T$ admits quantifier-elimination (as may be achieved by 
Morleyzation.)  Let $T'$ be a theory extending $T_\forall$, possibly but not necessarily
in a richer language.  We assume any model of $T'$ is definably closed as a substructure
of a model of $T$.     The generators $[X,h]$  of $\Ke(T)$ can
be taken with $X,h$ quantifier-free definable.  As such they are also elements
of $\Ke(T')$, and we 
define  $\mu$ on generators by $[X,h]_T \mapsto [X, h]_{T'}$.
Then equations \ref{1}-\ref{4}
 are respected by $\mu$, and we obtain a ring  homomorphism $\mu=\mu_{T/T'}: \Ke(T) \to \Ke(T')$.

\<{example} \lbl{cm1} $T=ACF_F$ where  $F$ is a perfect field,  $T'=Th(F)$,
both in the language of $F$-algebras.  \rm This will interest us especially when
$F$ is finite or pseudo-finite.

When $F$ is finite, $\Ke(T')$ is the group ring $\Zz[(F,+)]$.  Hence if we choose
a homomorphism $\psi: (F,+) \to (\Cc,\times)$ we obtain by composition
a nontrivial homomorphism $\mu_\psi: \Ke(T) \to \Cc$.

When $F$ is pseudo-finite, $\K(T')$ is related to virtual Chow motives, cf. \cite{DL}.

\>{example}

The words ``definable", etc. continue to refer to $T$ unless otherwise indicated.

 If $X$ is finite, we let $|X|$ be the number of points of $X$ in a model of $T$.  
 
 \<{lem} \lbl{cm2}  Let $X$ be an irreducible definable set.  
  Let $M \models T, a \in X(M)$.  Then
 $\Ke(T_a)$ is naturally isomorphic to the ring $S$ of definable functions
 $X \to \Ke(T)$.  \>{lem}
 
 \prf  Given a definable function $X \to \Ke(T)$, we obtain by evaluation at $a$
 an element of 
  $\Ke(T_a)$.  This gives a homomorphism $ev: S \to \Ke(T_a)$.
 Any definable set of $T_a$  is $T_a$-definably isomorphic to one of the form
  $Y_a = f \inv(a)$, where
 $Y$ is a definable set and $f: Y \to X$ a definable function.  
 This shows that the evaluation map $S \to \Ke(T_a)$ is surjective.  
  Moreover any 
 definable function $f$ on $Y_a$ is the restriction of a definable function on $Y$.  
 By
 irreducibility of $X$, if $f$ restricts to a bijection $Y_a \to Y'_a$
 then it must be a bijection $Y \to Y'$ over $X$.    The required isomorphism 
 follows already at the level of semirings, and hence extends to an isomorphism of rings.
 \eprf
 
 We write $\Ke(T_X)$ for either of the rings in \lemref{cm2}.  If $X$ is finite, we have the norm map
 $N_{T_X/T}: \Ke(T_X) \to \Ke(T)$.  Let $T'=Th(F)$, as above.
 
\<{lem}\lbl{frobenius}    Let $X$ be a finite irreducible definable set.  The composition $\mu_{T,Th(F)} \circ  N_{T_X/T}$ is a ring homomorphism.
\>{lem} 

\prf   Say $|X|=n$.  Write $\mu = \mu_{T,Th(F)} $ and $N = N_{T_X/T}$.  Since $N$ is multiplicative and
$\mu$ is a ring homomorphism, it is clear that $\mu \circ N$ is multiplicative.  
Let $a,b: X \to \Ke(T)$
be definable functions, and $c=a+b$.  Then 
$$N(c)=N(a)+N(b) + \sum_{m=1}^{n-1} \sum_{s \in Sym_m(X)} \prod_{x \in s} a(x) \prod_{x \in X \setminus s} b(x)$$
Now $\mu(  \sum_{s \in Sym_m(X)} \prod_{x \in s} a(x) \prod_{x \in X \setminus s} b(x)) = 0$
simply because $Sym_m(X) (F) = \emptyset$ for $0<m<|X|$.  Hence $\mu N (c) = \mu N(a) + \mu N(b)$.  
\eprf

The ring homomorphism $\mu \circ N$ induces a ring homomorphism on the power  series rings,
$\Ke(T_X)[[t]] \to \Ke(T')[[t]]$, with $t \mapsto t^n$.   This homomorphism coincides with
$a \mapsto \mu( \Pi_{x \in X} a(x))$.  
We note the corollary:

\<{lem} \lbl{irr-prod} Let $X$ be a finite irreducible definable set, $n =|X|$.  
Let $a_k: X \to \Ke(T)$ be a definable function, with $a_0=1$.  
Let $b_k(x) =  \prod_{x \in X} a_k(x)$.  
Then 
$$  \mu ( \prod_{x \in X} \sum_{k=0}^\infty a_k(x) t^k) = \mu ( \sum_{k=0}^\infty b_k(x) t^{nk})$$
\>{lem}

\prf  Follows from \lemref{frobenius}.  Alternatively it can be proved directly, 
using the fact that there are no $F$-definable nontrivial partitions of $X$.
\eprf 
  
  We mention in passing an additional, straightforward  compatibility of motivic volumes with ultraproducts,
in particular of formulas over pseudo-finite fields with corresponding motivic formulas over finite fields.  In the next paragraph we compare the   latter with classical  adelic integration.

\ssec{Compatibility(2)}

In this section we let $T=ACF_F$, with $F$ a finite field.  We wish to compare motivic adelic volumes
to classical ones.  For this purpose we do not need exponential sums, so
let $\K$ be any algebra over the Grothendieck ring of $T$ (typically,
obtained by inverting the class of the affine line.)

While the classical treatment of adeles has a  factor for each closed schematic point of the curve $C$,
ours has a factor for each point of $C$ in $F^{alg}$.   Motivic volumes are nevertheless
compatible with adelic volumes; for this the ``Frobenius'' \lemref{frobenius} is essential.

let $X$ be a scheme of finite type over $F$, and $\cU$ a scheme
of finite type over $X$.  For $x \in X(M)$ (where $M \models T$) 
let $a(x) = [\cU_x] \in \Ke(F(x))$.  This gives a definable function $X \to \Ke(F(x))$.  
On the other hand, let $X_{closed}$ be the set of closed schematic points
of $X$.  For $v \in X_{closed}$, let $F_v$ be the residue field, a finite extension of $F$, with $q_v$ points.
Let $U_v$ be the fiber of $\cU$ above $v$, and $a(v) = | U_v(F_v)|$.  
 Recall    \defref{products}. 

\<{lem} \lbl{comp3}  $$ \prod_{v \in X_{closed}} (1+a(v) q_v^{-s}) = \mu_F \prod_{x \in X} (1+a(x)t)_{t=p^{-s}}$$
Hence if
$ \prod_{v \in X_{closed}} (1+a(v) q_v^{-s})$ converges absolutely to some $r \in \Cc$;
  then 
$ \mu_F \prod_{x \in X} (1+a(x)t)$ converges absolutely to $r$ at $t=p \inv$.  
    \>{lem}

\prf   Any $v \in X_{closed}$ corresponds to a finite definable subset $X_v$ of $X$,
a Galois orbit.    We have $\mu_F \prod_{x \in X_v} (1+a(x)t) = 1 + \mu_F \prod_{x \in X_v} a(x) t^{deg(v)} $, by \lemref{irr-prod}.

\claim{1}We have an equality of formal series:
$$\mu_F \prod_{x \in X} (1+a(x)t) = \Pi_{v \in X_{closed}} \mu_F \prod_{x \in X_v} (1+a(x)t)$$ 
 
\prf  If $X$ is {\em finite}, then $X= \union _{v} X_v$, and the claim is a special case of 
\lemref{products-trans}, even without applying $\mu_F$.
  In general,  consider the coefficient $b_N$ of $t^N$ in  $\prod_{x \in X} (1+a(x)t)$.  
By \defref{products} we have $\mu_F b_N =   \sum_{s \in Sym_n(X)(F)} \mu_F \prod_{u \in s} a(u)$.
This is a finite sum, so for a sufficiently large finite definable $X' \subset X$,
the coefficient of $t^N$ in the finite product $\mu_F \prod_{x \in X'} (1+a(x)t)$
is the same.   Take $X'$ so large that $X''$ contains no finite definable set of size $\leq N$.  Let $X''=X \m X'$.  On the right hand side, the product decomposes into a product
over $X'_{closed}$ and a product over $X''_{closed}$; the latter has no nonconstant
terms of degree $\leq N$; so on the right hand side too the $t^N$ term for
$X$ and for $X'$ is the same.  The claim follows. \eprf

On the other hand $\mu_F \prod_{x \in X_v} (1+a(x)t) = 1 + \mu_F \prod_{x \in X_v} a(x) t^{deg(v)}  = 1 + a(v) t^{deg(v)}$, so $(\mu_F \prod_{x \in X_v} (1+a(x)t))_{t=p^{-s}} =
1+a(v)q_v^{-s}$.  The lemma follows.    \eprf

The above lemma (slightly generalized) will permit the comparison of   classical Tamagawa measures
(with convergence factors) to motivic ones.    

A similar comparison is valid for $\K_e$.  
 Let $\bF$ be a finite field $GF(q)$, together with a choice of character $\psi: (F,+) \to \Cc^*$.

 Define $\mu_{\bF}([X,h]) = \sum_{x \in X(F)} \psi(h(x))$.   Then equations \ref{1}-\ref{4}'
 are respected by $\mu_{\bF}$; so a homomorphism $\mu_{\bF}:  K_{exp}(T) \to \Cc$.  We have $\mu_{\bF}(\Ll)=q$,
and hence $\mu_{\bF}$ extends to a ring homomorphism $K_{exp}(T)(\Ll \inv) \to \Cc$.  
We also write $\mu_{\bF}$ for the natural extension to $K_{exp}(T)[[t]] \to \Cc[[t]]$.
We have:

	\<{lem} \lbl{comp2}  Let $X$ be a finite irreducible definable set.  Let $a: X \to \Ke(T)$ be a definable function. Let $x_0$ be a point of $X(F')$, where $F'$ is a finite field extending
	$F$, and assume $F'$ is generated by $x_0$ over $F$.  
	Let $\psi'(x') = \psi ( tr_{F'/F} (x'))$. 
	  Then $\mu_{\bF} ( \prod_{x \in X} a(x) ) = \mu_{F'}  a(x_0)$.
	\>{lem} 
The proof is similar to the case of counting.

\>{section}

 \<{section}{Local test functions and integration on linear spaces}  \lbl{localtestfns}
 
 Let $k \models T$.  Let $K$ be a discrete valued field, with residue field $k$.   In any valued field, $\Oo$ will always denote the valuation ring, $\Mm$ the maximal ideal.
 A {\em parameter} is an element $t$ of minimal positive valuation $1$.   Let 
 $$K(N;M) = \{x: \val(x) \geq -N \} / \{x: \val(x) \geq M \} =t^{-N} \Oo / t^M \Oo$$

We have an isomorphism $\alpha_{N,M}: k^{N+M} \to K(N;M)$, 
 $$x=(x_{-N},\ldots,x_M) \mapsto \sum x_i t^{i} $$
By  a {\em test function of level (N,M)} we mean a  function $\phi: K(N;M) \to K_e(T)$, such that
  $\phi^t :  = \phi \circ \alpha_{N,M}$ is a 
    definable function $\phi^t: k^{N+M} \to K_e(T)$.    We view $\phi$ as a function 
    on $K$ whose support is contained
 in $t^{-N} \Oo$, and invariant under $t^M \Oo$.   Let $\CS(K;N,M)$ be the set
 of test functions of level $(N,M)$, and let  $\CS(K) = \union_{N,M \in \Nn} \CS(K;N,M)$.

  We define
$\int^t:  \CS(K;N,M) \to K_e(T)$ by
$$\int^t \phi =  \Ll^{-M}  \sum_{x \in  k^{N+M}} \phi^t(x)  $$  

  If $t'$ also has valuation $1$, then for each integer $r$,  
  $t'= \sum_{i=1}^{r} a_it^i (\mod t^{r+1})$, 
 for some (uniquely defined) $a_1,\ldots,a_{r} \in k$.   We call $t,t'$ {\em definably
 equivalent} if for each $r$, each $a_i$ is $T$-definable.

\<{lem} \lbl{local1}  If $t,t'$ are definably equivalent parameters, then  $\int^t, \int^{t'}$ agree on $\CS(K;N,M)$.  \>{lem} 

\prf   $(\phi^t)^{-1} \phi^{t'} $ is definable and hence induces the identity
on $K_e(T)$.   \eprf

Hence $\int^t$ depends only on the definable equivalence class of $t$.  We will later
describe a canonical choice of a definable equivalence class of parameters.  
 We thus cease to write $t$ in the superscript, letting $\int = \int^t$.

\<{lem} \lbl{local2}
\begin{enumerate}
  \item   Let $i: K(N,M) \to K(N+1,M)$ be the natural inclusion.  It induces
$i_*: \CS(K;N,M) \to \CS(K;N+1,M)$ (extension by $0$.)   We have $\int \circ i_*  = \int$. 
  \item Let $\pi: K(N,M+1) \to K(N,M)$ be the natural projection.  It induces
  $\pi^*:  \CS(K;N,M) \to \CS(K;N,M+1)$ (composition).  
  The image of $\pi^*$ is the set of elements $\phi$ of $\CS(K;N,M+1)$ invariant
  under $\ker(\pi)$.  
  We have $\int \circ \pi^*   = \int$. 
\end{enumerate} 
 \>{lem}  

\prf Clear.  Regarding the image of $\pi^*$,  note that $\pi$ has a definable section. 
In general  it is   clear   that if
$E$ is an equivalence relation on $Y$ with a definable section $Y/E \to Y$, then 
an $E$-invariant definable map on $Y$ arises from a definable map on $Y/E$. \eprf

We view both $i_*$ and $\pi^*$ as inclusion maps, and let $\CS(K) = \union_{N,M \in \Nn} \CS(K;N,M)$.   We view the elements of $\CS(K)$ as (``smooth, bounded'') functions on $K$; 
  the integral  as a function on   $\CS(K)$.  
 
The same goes for families.   If $Y$ is a definable subset of $k^n$,  $\alpha_{N,M}$ yields a notion of definable function
on $Y \times K(N;M)$; and given such a definable function $\phi(y,x)$ we can integrate
with respect to the $K(N;M)$ variables to obtain a definable function $Y  \to K_e(T)$.

Let $K^m(N,M) =K(N,M)^m$.  Define
 $\CS(K^m;N,M)$,  $\CS(K^m) = \union_{M,N} \CS(K^m;M,N)$, and 
$\int: \CS(K^m) \to K_e(T)$ in a similar way.  If $U$ is a vector space defined over $K$   with a distinguished basis, 
we identify $U$ with $K^m$, and in this way define $\CS(U;N,M)=\CS(K^n;N,M)$ etc.

\<{rem} \rm In fact a definable element of $GL_m(K)$ induces an isomorphism
$\CS(U) \to \CS(U)$ (not in general preserving levels), and using this system of isomorphisms it is possible
to define $\CS(U)$ for any $U$ having a definable basis (but not a distinguished one); but we will not need this. \>{rem}

\ssec{Finitely many valued fields}

Assume given finitely many discrete valued fields $(K_i: i \in S)$, each
 with a parameter $t_i$ and residue field $k$.    We write $K_S$ for the ring $\Pi_{i \in S} K_i$.  Let
 $$K_S(N,M) = \Pi_{i \in S} {K_i}(N,M) $$  
 We  define 
 $\CS(K_S^m)$ and  $\int: \CS(K_S^m) \to K_e(T)$ in the same way as for one field above.

 In practice the $K_i$ will extend $(F,v_i)$ where $F$ is a fixed field, and $v_i$ are valuations
 of $F$.  We will also have a vector space $U$ over $F$ of dimension $m$, with an $F$-basis, and write
 $\CS(U;K_S) = \CS(\Pi_{i \in S} K_i^m)$ using the   basis for the identification.

\<{rem} \rm Classically, the map $\CS(K_1) \tensor \ldots \tensor  \CS(K_n) \to \CS(K_1\times \ldots \times K_n)$ is surjective,
but this will not be the case for us.   Nor is the image of this homomorphism preserved under the ``sum of rational points" maps of \S \ref{5.4}. \>{rem}

\ssec{Fourier transform}  \lbl{fourier-local}  Let $k,K,t,N,M$ be as above.  We also   fix a
nonzero linear map $r:K \to k$, vanishing on $t^{-M} \Oo$ for some $M$.     The dual of $\Oo$
with respect to $r$ is $\Oo ^{\perp} = \{x: (\forall x \in \Oo) r(xy) =0 \}$; it contains $t^{M} \Oo$,
and is an $\Oo$-module, so it must have the form $t^{\nu} \Oo$ for some $\nu \in \Zz$.
We assume $\nu$ is even.

We  define the local Fourier transform of $\phi \in \CS(K;N,M) $
by
\eq{fourier}{  \Ff(\phi) (x)  = \Ll^{-\nu/2} \int_y \phi(y) \psi(r(xy))   }

It is clear that $\Ff(\phi)$ is invariant under $\{x: \val(x) \geq  N+M \}$.
Using \lemref{l4}
    one sees that $\Ff(\phi)$ is bounded.  Hence
$\Ff: \CS(K) \to \CS(K)$.

The inversion formula is easily proved.
 
On $K^n$ we define 
  define $\Ff: \CS(K^n) \to \CS(K^n)$ by the same equation \ref{fourier}, with $xy$ interpreted as the standard dot product.
  
 More generally, given finitely many valued fields $K_1,\ldots,K_n$, each with a parameter
 $t_i$ and a linear map $r_i: K_i \to k$, we define for $\phi \in \CS(K_1 \times \ldots \times K_n)$ 
 $$\Ff(\phi)(x_1,\ldots,x_n) = \int_y \phi(y_1,\ldots,y_n) \psi (\sum r_i(x_iy_i)) $$
so that the Fourier transform can be computed one variable at a time.

\ssec{Comparison with \cite{HK}}

Let $K$ be a discrete valued field, with distinguished subfield $F \leq \res(K)$
and a distinguished parameter $t$.  Assume $char(F)=0$.  Then for test functions
$\phi$ both the above theory and the  integration theory of \cite{HK} apply.  
  We explain the connection.  
Let $i(\phi)$ denote the integral in the sense of \cite{HK}.  Let $T$ denote the theory
of the residue field of $K$ (including constants for $F$), and let $T'$ be the theory
of $K$ (including  the $T$-structure on the residue field and a constants for  $t$.)
 Let $j$ denote
the ``rational points'' functor of \cite{HK}, chapter 10,  towards the theory $T'$.   In general
$ji$ takes values in $K^+(T) \tensor K^+(\G) / I_{sp}$.  
Let $K^+_{fin}(\G) $ be the subsemiring of $K^+(\G)$ represented by finite definable sets.
Let $L = [G_a] \in K^+(T)$.   
 We have a ``weighted counting''
homomorphism $K^+_{fin}(\G) \to \Nn[L \inv ] \leq K(T)[L \inv]$, where $L =[G_a]$; namely, the point $n \in \Zz \leq \G$ is assigned the value $[G_m(k)]L^{-n}$; sum over finite sets. 
This induces
$h: K^+(T) \tensor K^+_{fin}(\G)  \to K(T)[L \inv]$   
respecting the restriction of the congruence $I_{sp}$.  (In fact $h$ is an isomorphism
between the image of $ K^+(T) \tensor K^+_{fin}(\G)$ in   $K^+(T) \tensor K^+(\G) / I_{sp}$,
and $K^+(T)$.)

\<{lem}  (1)     Suppose     $\phi$ takes values in $K^+(T)$.  Then $ji(\phi)$
 lies in $K^+(T) \tensor K^+_{fin}(\G) / I_{sp}$, and $hji(\phi) = \int \phi$.  \>{lem}
 
 \prf    The maps $h,j,i$  and $\int$ are  completely additive
over the residue field, i.e. if $W$ is definable over $T$
and $\phi = \sum_{w \in W} \phi_w$ then 
$hji(\phi) = \sum_{w \in W} hji(\phi_w)$, and similarly for  $\int \phi$.  Hence
the statement on $ji(\phi)$ and the equality   $hji(\phi) = \int \phi$ reduce to the case of characteristic functions of a
point in $K(n,m)$; this is an immediate computation.  \eprf

This immediately gives a change of variable formula, and the ability to integrate over varieties.  
Unfortunately it only applies in characteristic $0$ so in general we are obliged to 
do this from scratch.

\>{section}

\<{section}{Global theory}
   
We continue with the theory of fields $T$, containing constants for a field $F$.

Let  $C$ be a smooth, projective, absolutely irreducible  curve over $F$.  Let $K=k(C)$ be the function field.   We will define global test functions and construct a Fourier transform operator,  and a ``sum over rational points'' functional, on them.

 We also fix a nonzero 1-form $\omega$ on $C$, defined over $F$.   We assume  $\omega$ can be chosen in such a way that 
   every zero or pole of $\omega$ has even multiplicity.  
      For $g=0$ one can choose $\omega$
  with one double pole, and no zeros; for $g=1$ one can choose $\omega$ regular, with no
  zeros.  When $g>1$ such a form may not exist over $F$; the next paragraph contains a proof   that it can always be found over a finite   extension.  Various remedies are possible when one wants to remain over $F$, including the introduction of a formal square root of the affine line;
since we will not require this case we will not enter into this discussion.  
   
  We have $\deg(\omega)=2g-2$, where $g$ is the genus of $C$.  
Assume $g>1$, and let $ Jac(C)_{(2g-2)}$ be the 
   divisor classes of degree $2g-2$, a torsor of the Jacobian $Jac(C)$.   
  Consider the map $f: C^{g} \to   Jac(C)_{(2g-2)}$ given by
 $(c_1,...,c_{g-1},c_g) \mapsto  2(c_1+...+c_{g-2}-c_{g-1}) + 4 c_g$. 
   By subtracting $f(c_0)$ for some $c_0 \in C^g$
 one obtains a map into the Jacobian,
and if the image has dimension $h<g$ then the Jacobian is easily seen to have dimension
$h$, a contradiction.   Hence the image $f(C^{g})$ has dimension $g$.  
Now $\dim(Jac(C)_{(2g-2)}) = \dim(Jac(C))=g$.  
   Since $C^{g}$ is complete, the image is closed, so $f(C^g)= (Jac(C)_{(2g-2)})$.  
    So we can find $(c_1,...,c_{g})$ such that 
$ 2(c_1+...+c_{g-2}-c_{g-1}) + 4 c_g$ represents the canonical class. Then there exists a form 
$\omega$ such that $div(\omega) = 2(c_1+...+c_{g-2}-c_{g-1}) + 4 c_g$.

For simplicity we assume that  $T$ contains $ACF_F$.  This is not a serious restriction
since the language may be larger than that of $F$-algebras, and may in particular
include a predicate for a subfield $K$.  

If $v$ is a valuation on $F(C)$, the residue field of $v$ is a finite
 extension $F_v$ of $F$.   The assumption that $T \models ACF_F$ is used to conclude that $F_v$ is a subfield  of $k$; this simplifies the notation.   
    It is not however contained in $F$, so a direct translation of the classical theory would
    lead to integrals with values in  of $K_e(T_{F_v})$ rather than $K_e(T_F)$.  It is not clear
    how to multiply elements of $K_e(T_{F_v})$ for distinct $v$.  
  However adelic integration with 
  our ``redundant'' definition of the adeles  involves taking products  over several conjugate
  representatives of the same valuation; 
 this means that the integral factors through the norm map,  hence does belong to $K_e(T_F)$,
 and further products over distinct valuations make sense.  
   Similarly if the language contains a predicate
for a subfield $K$, the integrals of quantities defined over $K$ will themselves be over $K$.  
 
  Our ``sum over rational points'' is actually a sum over $k(C)$, not $F(C)$, 
  including notably $F^{alg}(C)$-rational points.     This is necessary
  to allow uniformity in definable families,  e.g. \lemref{g2}.     Nevertheless we show 
  compatibility with the classical sum, via $\mu_F$.  And the Poisson summation formula
  holds motivically (in the sense of motivic integration), before $\mu_F$ is applied.  
 
\ssec{$k(C)$ as an Ind- definable field}.

Let $L^g_{F}$ be the language including the following symbols, {\em all viewed as relation symbols.}
 
1)  A sort $K$ intended to denote $k(C)$.  On $K$, the language  of $F$-algebras,  made relational, i.e.  with a relation symbol for the zero-set of each polynomial over $F$;
A unary predicate symbol for $k \subseteq K$.

2)  Relation symbols $V_n \subseteq C(k) \times K$.  Intended meaning:
$V_n(\a,f)$ if $ord_\a(f)=n$, i.e. the order of vanishing of $f$ at $\a$ equals $n$.

4)  ${\bf t} \subseteq V_1$ intended to pick out a parameter $t_c$ for $K_{v_c}$, uniformly in $c$. 

Let $F \leq k$.  We impose the natural $L^g_F$ structure on $k(C)$ (with auxiliary sort $k^{alg}$, and $K$ interpreted as $k(C)$.)  

5)  $W = \{(\a,f,\g) \in C(k) \times k(C) \times k: \ \g =  res_\a (f \omega)  \}$.

 According to \lemref{piecewise-g},  with this structure, $k(C)$ is piecewise definable over $k$.

\ssec{Global test functions}  \lbl{global-test-functions}

 For each $u \in C(k)$ let  $k(C)_u$ be the valued field $k(C)$, with valuation   coresponding to the point $u$, i.e. $\val(f) >0$ iff $f(u)=0$.  We could take the completion but for
 our immediate purposes it is not important since we really use only the vector spaces 
 $t^{-m} \Oo_u / t^m \Oo_u$, where
   $\Oo_u$ denotes the valuation ring and $t$ is a parameter.  
   Let $\Aa$ be the   restricted product of fields $K_u$ relative to the rings $\Oo_u$.  
(Classically one takes only algebraic $u$, and only one copy for each conjugacy class;
this would suffice to tell apart our test functions.)

A {\em   global  test function on $\Aa^n$} is given by a finite $S \subset C(k)$ and
an element $\phi$ of $\CS((\Pi_{u \in S} k(C)_u)^n)$, with the understanding 
that if $S'$ is disjoint from $S$, then $(S,\phi)$ is identified with $(S \union S', \phi \tensor 1_{\Oo_{S'}^n})$.   Here $\Oo_{S'} = \Pi_{u \in S'} \Oo_u$.    
 
 A {\em definable subset} of $\Aa^n$ is 
defined similarly, so that the characteristic function of a definable subset is a test function.   

The form $\omega$ on $C$ provides us with a linear map $k(C)_u \to k$, namely
$f \mapsto res_u(f \omega)$.  It is this map that we use in \ref{fourier-local}, to obtain
a local Fourier transform $\Ff_u: \CS(k(C)_u) \to \CS(k(C)_u)$.    If $S$ is a finite
subset of $C$, we define $r_S(f) = \sum_{u \in S} res_u(f \omega)$.

\ssec{Fourier transform}

The Fourier transform   of a global test function  on $\Aa$ is defined by
choosing a representative $(S,\phi)$ for the global test function such that
$\omega$ is regular and nonzero outside $S$, and letting 
$$\Ff((S,\phi)) = (S,   \Ff(\phi))$$
 the latter $\Ff$ being the semi-local Fourier transform defined
at the end of \ref{fourier-local}.  It is easy to check that this is well-defined.  

We can also define $\Ff((S,\phi) = \phi'$ directly, $\phi' =  \L^{1-g} \int \phi(r_S(xy)) \phi(x) dx$ .

\ssec{Summation over rational points}  \lbl{5.4}

Let $u \in C(k)$.  View $k(C)_u$ as an piecewise-definable valued field in $T$, with distinguished
parameter $t=t_u$.   For $f \in k(C)$, we write $f_u$ for $f$ viewed as an element of $k(C)_u$.

Fix a definable global test function $\phi$, represented as $(S,\phi_S)$ for some $S$.  
We will define $\phi(f)$ for  $f \in k(C)$. 

Let $k(C)_S =  \{f \in k(C):  (\forall u \notin S)(v_u(f) \geq 0 ) \}$.   For $f \notin k(C)_S$
we define $\phi(f)=0$.  This is forced by the definition of global test functions,
since if $v_u(f)<0$, then $\phi$ is also represented by $(S \union \{u\}, \phi_S \tensor 1_{\Oo_u})$ and $1_{\Oo_u}$  vanishes at $f$.  
 
For $f \in k(C)_S$,   let  $f_S = (f_u: u \in S) \in \Pi_{u \in S} k(C)_u$.  Define
$\phi(f) = \phi_S(f_S)$.  

It is clear that $\phi(f)$  does not depend on the choice of $S$. 

If $f$ is definable, then $\phi(f) \in K_e(T)$.  In general if $f,\phi$ are $F'$-definable then
$\phi(f) \in K_e(T_{F'})$.

Let $Y$ be a limited subset of $k(C)$.  Then $y \mapsto \phi(y)$ is clearly a definable
function $Y \to K_e(T)$.
We thus  have an element $\sum_{y \in Y} \phi(y) \in K_e(T)$.  (cf. Equation \eqref{summ}).

Let $m$ be an integer such that $\phi_S$ is supported on $\prod_{u \in s} t_u^{-m} \Oo_u$.
Let
 $$Y_0 = \{f \in k(C)_S: (\forall u \in S)  v_u(f) \geq -m \}$$
  Then $Y_0$ is a  limited subset of   $k(C)$.  We have $\phi(f)=0$ for any $f \in k(C)_S \setminus Y_0$, hence for any
$f \in k(C) \setminus Y_0$.  

Define
$\CR(\phi) = \sum_{y \in k(C)} \phi(y) = \sum_{y \in Y_0} \phi(y)$.

If $\phi=\phi(y,x)$ depends on other variables $x=(x_1,\ldots,x_n)$, we define $\sum_{y \in 
k(C)} \phi(y,x)$ in the same way, and denote it $\CR_y \phi $.

\<{lem}  \lbl{g1} Let $\phi$ be a    test function in $n+1$ variables $x_1,\ldots,x_n,y$.  Then 
$\CR_y \phi = \sum_{y \in k(C)}^{} \phi(x,y) $ is a   test function in $n$ variables.  \>{lem}

\prf   Definability is clear.  
 By assumption, $\phi_S(x,y)$ is invariant under $C_S^{n+1}$ for some congruence
 subgroup $C_S$.   From this it
 is clear that $\CR_y\phi $ is $C_S^n$-invariant.  Since $\phi_S$ is supported on $D^{n+1}$ for some bounded set $D$, it is clear that $\phi'$ is supported on $D^n$.   
 
\eprf

Call a global test function {\em simple} if it is represented by $(S,\phi_S)$, where
$\phi_S$ is the characteristic function of a single coset of $t_v^m \Oo_v$, for some $m$.

\<{lem} \lbl{g2} (1)  Let $W$ be a definable set of $T$, and let 
 $\phi_a$ be a global test function, defined uniformly in $a \in W$.  Let 
$\chi: W \to K_e(T)$ be a definable function.  Then
$\sum_{a \in W} \chi(a) \phi_a$ is a global test function.

(2) $\sum_{x \in k(C)} \sum_{a \in W} \phi_a(x) = \sum_{a \in W} \sum_{x \in k(C)} \phi_a(x)$.

(3)  $\int_x \sum_{a \in W} \chi(a) \phi_a(x) = \sum_{a \in W} \chi(a)  \int_x \phi_a(x)$

(4)  $\Ff \sum_{a \in W} \phi_a = \sum_{a \in W} \Ff \phi_a$

(5)  Any global test function can be expressed as in (1). 
\>{lem}

\prf  (1)-(4) are clear.  As for (5), the test function is  represented by $(S,\phi_S)$ for some $S$;
and  for some $m$,
 $\phi_S$ factors through a function $\phi'$ on $\prod_{u \in S} t_u^{-m} \Oo_u / t_u^m \Oo_u$, pulled back to $\prod_{u \in S} t_u^{-m} \Oo_u$ and extended
 by $0$ to $\prod_{u \in S} K_u$. 
 Identify 
  $t_u^{-m} \Oo_u / t_u^m \Oo_u$ with  $ \k^{2m} $ via the  basis $t_{m-1},\ldots, t_{-m}$.  
   Let $W= (\k^{2m})^S$.  For 
$a \in W$ let $\phi_a$ be the simple test function concentrating on $a$,
and $\chi(a) = \phi'(a)$.
Then clearly $\phi = \sum_{a \in W} \chi(a) \phi_a$.  

\eprf

\<{lem} \lbl{comp4}  Assume $F$ is a finite field, and fix a nontrivial character $\psi$ 
of $F$.  Let $\mu_\psi$ denote $\mu_{T/Th(F)}$ composed with $\psi: \Zz[(F,+)] \to \Cc$.
Let $\phi$ be a definable global test function in $n$ variables.  Then 
$\mu_\psi (\sum_{x \in k(C)^n} \phi (x))= \sum_{x \in F(C)^n} \mu_\psi \phi (x)$.  
\>{lem}

\prf By opening up the definitions.  The Ind-definable sum $\sum_{x \in k(C)^n}$ reduces
to a certain definable sum $\sum_{y \in Y_a}$ (with limited $Y_a$ depending on $\phi_a$.)
Now the union $Y=\union_a Y_a$ is still limited, and we may write $\sum_{x \in k(C)^n} = \sum_{y \in Y}$ for any of the test functions in question.  
 In general $\mu_\psi$ commutes with definable sums
$\sum_{y \in Y}$.  \eprf

\<{lem}\lbl{g3}  Let $a \in F(C)$, and $\phi'(x)=\phi(x+a)$.   Then $\CR \Ff (\phi') = \CR \Ff(\phi)$.
\>{lem}
\prf We may assume $a \omega$ is holomorphic at $u$  for $u \notin S$.  Then
$\Ff(\phi') = \psi(r_s(ab)) \Ff(\phi)$.  In particular for any $b \in k(C)_S$
we have
$\Ff( \phi')(b) = \psi(r_s(ab)) \Ff(\phi)(b)$.  
  But the sum of residues $r_s(ab)= \sum_{u \in F(C)} \res_u(ab \omega)  = 0$.  \eprf

\ssec{Poisson summation formula}

\<{lem}   \lbl{poisson1}  Let $\phi(x)$ be a definable global test function,
and $\psi(y) = \Ff \phi(x)$.  Then $\CR \psi = \CR \phi$.  \>{lem}
 
\prf In view of \lemref{g2} we may assume $\phi$ is simple.  We compute using representatives in $ \prod_{u \in S} k(C)_u$, where $S$ is a large enough
finite set.  Say $\phi$
is represented by $(S,\phi_S)$ with $\phi_S$ the characteristic function of
$W=a+ \prod_{u \in S} t_u^{m_u} \Oo_u$, $a \in \prod_{u \in S} k(C)_u$.
  Let $D$ be the divisor on $C$ supported at $S$,
with multiplicity $m_u$ at $u$.   Then $L(D) := \{f \in k(C): (\forall u \in C)(v_u(f) \geq m_u)\}$
is a finite-dimensional subspace of $k(C)$, defined over $F$.  Let $D' = div(\omega)-D$
be the dual divisor.  We write $k(C)_S$ for the image in $\prod_{u \in S} k(C)_u$
of $\{a \in k(C): (\forall u \notin S)( v_u(a) \geq 0) \}$.  

\case{1}  There exists $a' \in k(C)_S \meet W$.  

Since $ k(C) \meet W$ is a torsor
for $L(D)$, by Hilbert 90, $F(C) \meet (a+ \prod_{u \in S} t_u^{m_u} \Oo_u) \neq \emptyset$.
So we may take $a' \in F(C)$.  By \lemref{g3} we may translate by $a'$; so we may 
assume $a'=0$.  In this case by direct computation we see that
$\Ff(\phi) = \Ll^{deg(D)+1-g} \phi'$ where $\phi'$ is the simple global test function concentrating
on $D'$.  The equality asserted in the lemma is precisely Riemann-Roch.

 \case{2}  $k(C)_S \meet W = \emptyset$.
In this case we have $\CR \phi = 0$ and we have to show that $\CR \Ff \phi = 0$.

Let $A = \{x \in \prod_{u \in S} k(C)_u  (\forall u \in S) (v_u(x) \geq m_u) \}$.  
Using Riemann-Roch, $A / (k(C)_S \meet A)$ is finite dimensional.  We use the
form $r_S(xy)$ on $ \prod_{u \in S} k(C)_u$.  By definition of the Fourier transform 
and \lemref{l4}
, $\Ff \phi$ is supported on $A^{\perp}$.   For $y \in A^{\perp}$,
$r_S(xy)$ takes a constant value $\rho(y)$ for all $x \in W$.   This $\rho$ is a linear
map on $A^{\perp}$ and in particular on $B=A^{\perp} \meet k(C)_S$.  
We have $\CR \Ff \phi =vol(W)  \sum_{y \in B}   \psi \rho(y)$.  If $\rho$ 
 vanishes on $B$   then $W \subseteq B^{\perp} = A + (k(C)_S)^{\perp}$ (as one obtains 
 by factoring out  to reduce to a finite dimensional situation.)  But $k(C)_S$ is self-dual
 for $r$, since the sum of residues equals zero.  Thus $W \subseteq A + k(C)_S$,
 contradicting the case definition.   Thus $\rho$ is not constant on $B$; but then
 by \lemref{l4}
 we have  $\sum_{y \in B}   \psi \rho(y) = 0$, so $\CR \Ff \phi =0$.  
 
\eprf

\<{lem}   \lbl{poisson2}   Let $\phi(x,z)$ be a definable global test function, of several variables,
and $\psi(y,z) = \Ff_x \phi(x,z)$.  Then $\CR_y \psi = \CR_x \phi$.  
\>{lem}

\prf  We take a representative $\phi_S$ of $\phi$, with $S$ sufficiently large, as usual.
Since $\phi_S$ is smooth and of bounded support, we can view it as a function on 
a finite-dimensional $\k$-space; the same goes for $\CR_y \psi$ and $\CR_x \phi$.
To show that such functions into $K_e(T)$ are equal it suffices to show that
at any value $b$ of $z$ in any model of $T$, we have $\CR_y \psi(y,b) = \CR_x \phi(y,b)$.
This is just \lemref{poisson1} applied in $T_b$.  \eprf

\ssec{Characterization of $\CR$}

We remark that another proof of the Poisson summation formula is possible, using
a self-dual {\em characterization} of $\CR$ among definable distributions.  

The rational points functional clearly enjoys the following properties, where the equalities
(1,2) take  place in $\Ke(T_a)$:

\<{enumerate}
\item  $\CR (\psi(r(ax)) \phi(x)) = \CR(\phi)$ for $a \in K$.
\item  $\CR(   \phi(ax) = \CR(\phi)$ for $a \in K$.
\item  $\CR ( 1_\Oo) = \Ll = [G_a]$
\>{enumerate}

This is in fact a characterization of $\CR$ among definable distributions, provided
that we invert $\psi(c)-1$ for every $0 \neq c \in \k$ uniformly, and that $r$ is chosen
so that $K$ is self-dual for $r(xy)$.    We sketch the proof.

 Property (1) implies that $\CR$ concentrates on $K$-points.
This uses the fact that $K^\perp =K$ and the invertibility of $\psi(c)-1$ for $c \neq 0$.
 Property (3), along with (2),  implies that $\CR(\phi) = 1 $ for $\phi$ concentrated near $0$, and with $\phi(0)=1$.    Using (2) again we obtain this for any rational point.
 
Properties (1),(2) are exchanged by the Fourier transform, while (3) is left invariant.
Thus $\CR   \circ \Ff$ enjoys the same properties, giving (under the stronger assumptions) another proof that $\CR \circ \Ff = \CR$.

\>{section}

\<{section}{Theory of valued fields over a curve}  \lbl{qeSec}

\ssec{Valued fields with a field of representatives}  

We begin with the local ingredient of the logical theory we will use.  We take a three-sorted
language of valued fields, with a sort $VF$ for the valued field, a sort $
\G$ for the value group, and a sort $\res(VF)$ for the residue field.
We take the usual language of valued fields, including, after Delon, 
a binary function symbol
$\res(\frac{x}{y})$; defined to be $0$ when $\val(x) < \val(y)$, and otherwise
the residue of $\frac{x}{y}$.    In addition, 
the value group has a distinguished element $1>0$;  and there is  
a function symbol $i: \res(VF) \to VF$ for a section 
of $\res$.  Thus $k: = i(\res(VF))$ is a   distinguished
subfield $k \subset \Oo$, and   the residue map is bijective on $k$.   

The theory $T_{loc}$ asserts that $K$ is an  algebraically closed   field, $\val$ is a valuation,
with valuation ring $\Oo$ and maximal ideal $\Mm$, and $\Mm \oplus k = \Oo$.
 
Quantifier elimination for   pairs of valued fields much more complex than ours is known; see 
 \cite{leloup1}, who attributes the case of $T_{loc}$ to Delon.   However
 the quantifier elimination in \cite{leloup1} takes to be basic formulas such as
 $$(\exists x_1,x_2  \in k) (\val(y-x_1y_1-x_2y_2) > \val(y_3))$$
 asserting that $y$ is closer to the vector space $y_1k+y_2k$ than $y_3$ to $0$.
The language we use allows only to take the coefficients of actual members
of this vector space, not of nearby points.  It seems simplest to give a direct proof of QE.

\<{lem}\lbl{lin.disjoint}  Let $A$ be a subfield of $M \models T_{loc}$,  closed
under $i \circ \res$.   Then $A,k(M)$ are linearly disjoint over $k_A := A \meet k(M)$.  \>{lem}

\prf  We show by induction on $n$ that if $a_1,\ldots,a_n \in A$ are linearly dependent
over $k(M)$, they are  linearly dependent over $k_A$.  If some $a_i=0$ this is clear;
this covers the case $n=1$.    Assume
the statement holds below $n$, and let $a_1,\ldots,a_n \in A$ be linearly dependent
over $k(M)$.   Reordering, we may assume $\val(a_i) \geq \val(a_n)$ for $i \leq n$.  
Dividing by $a_n$ we may assume $a_n=1$ and $\val(a_i) \geq 0$ for each $i$.
Let $b_i = i \res(a_i)$; then $b_i \in k_A$.  Performing the column operation $a_j \mapsto
(a_j -b_j a_n)$, we may replace $a_j$ by $a_j-b_j$ for $j<n$, so we may assume 
$\res(a_j)=0$ for $j<n$.  Now for some $c_i \in k(M)$ we have
$\sum_{i=1}^n c_i a_i = 0$, or $c_n = - \sum_{i=1}^{n-1} c_ia_i$.  
  So $\val(c_n) >0$.  But $c_n \in k$; so $c_n=0$.  Thus $\sum_{i=1}^{n-1} c_i a_i = 0$,
  and by induction $c_1,\ldots,c_{n-1}$ are linearly dependent over $k_A$.  
\eprf

\<{lem}\lbl{qe-loc}  $T_{loc}$ admits quantifier elimination.  $k,\G$ are stably embedded and 
strongly orthogonal.  Their induced structure is the field and ordered group structure, respectively. \>{lem}

\prf 
Let $\Uu,\Uu'$ be saturated.   Write $k=k(\Uu), k'=k(\Uu')$ 

\claim{}  
 Let $f_k: k \to k'$ and $f_\G: \G(\Uu) \to \G(\Uu')$
be given isomorphisms.   Also let $f:A  \to A' $ be an isomorphism betwee small subrings of $\Uu,\Uu'$, such that:  \<{enumerate}
\item[i)]   $f$ is 
 a partial isomorphism of valued fields, carrying $k$ to $k'$.  
\item[ii)]
   (Compatibility) $f(x)= f_k(x)$ for $x \in k_A$, and    $\val \circ f =  f_\G \circ \val $ on $A$.  

\item[iii)]
 If $0 \neq a,b \in A$, $\val(a)=\val(b)$ then  $i \circ \res (a/b) \in A$ (hence the same holds for $A'$).
\>{enumerate}

   Then   $f_k \union f_\G \union f$ extends
to an isomorphism $\Uu \to \Uu'$.  

\prf   Note in (iii) that   $i \circ \res(a/b)$ is the unique element $c$ of $k_A$
such that $\val(a-bc) > \val(a)$; so $f(c)= i \circ \res(f(a)/f(b))$.
 
  It suffices to show that $f$ extends to a partial isomorphism with (i-iii)
  on a subring containing  $A$ and a given   element of $c$ of $\Uu$.   
 If this subring  is not small, we can always replace it by a small subring still
satisfying  (i)-(iii), by   L\"owenheim-Skolem.

0)  $f$ extends to the field generated by $A$:  it clearly extends to a valued field isomorphism,
commuting with $f_\G$.  Any ratio $b/b'$ of elements of the field of fractions has the form $a/a'$
for some $a,a' \in A$, so (since $\res({x/y})$ was taken to be a function of two variables)
$\res(b/b') = \res(a/a')$ and commutativity with $f_k$ is clear too.  Hence we may assume 
$A$ is a field.

1)  
 By \lemref{lin.disjoint},  $k,A$ are linearly disjoint over $k \meet A$,
and $k', A'$ are linearly disjoint over $k' \meet A'$.   Hence there exists a (unique) isomorphism $f_{kA} : kA \to k'A'$.
Recall  (\cite{leloup1}, 1.1.2) that 
$kA$ is separated over $k$ in the sense of \cite{baur}:
any finite dimensional $k$-subspace of $kA$ has a basis $c_1,\ldots,c_n$,
such that $\val(\sum_i a_i c_i )= \min \val(a_ic_i)$ for any $a_1,\ldots,a_n \in k$.  
It follows that $f_\G \val(x) \leq \val ( f_{kA}(x))$ for $x \in kA \setminus (0)$.   By symmetry the other inequality holds.  
 Hence $f_{kA}$ 
is a valued field isomorphism compatible with $f_\G$. 
Note that $\res(kA) = \res(k)$.  

2)  Extend $f_{kA}$ to the field generated by $kA$ using (0), and then to 
 a valued field isomorphism $f_2: \acl(kA) \to \acl(kA')$.  Note that (i),(iii)
holds trivially.
 
3)  Let $c \in \Uu$.  If for some  $a \in \acl(kA)$ one has $\gamma= \val(c-a) \notin \val(A)=\val(\acl(kA))$, let $c' \in \Uu'$
be any element with $\val(c'-f(a))=f_\G(\gamma)$.  Then $f_2$ extends
to $f_3: \acl(kA)(c) \to \acl(kA')(c')$ uniquely with $c \mapsto c'$.   Indeed any element
of $\acl(kA)(c)$ is a product of elements of the form $c-d$, $d \in A$, and for such elements
the data forces $\val(c-d) = \min (\val(d),\gamma)$; and similarly on the $A'$-side.

4)  Otherwise, $kA(c)$ is an immediate extension of $\acl(kA)$.     Note that $\G(\acl(kA)) = \G(A)$
and this is a small subset of $\G$.  Let $b$ be the set of $\gamma \in \G(A)$ such that
the ball $B(\gamma)=B_\gamma(c)$ contains a point of $\acl(kA)$; in this case 
$B(\gamma)$ is defined over $\acl(kA)$, so $f_3(B(\gamma))$ is a ball over 
$\acl(i(k) A')$.
 Using saturation, 
let $c'$ be any point of $\meet_{\gamma \in Y} f_3(B(\gamma))$.  Then as in (3), $f_2$
extends uniquely to an isomorphism on $\acl(kA)(c)$ with $c \mapsto c'$.   
Then $f'$ satisfies the conditions for $f$.   
\eprf
 To prove  quantifier elimination, we need to extend a partial  isomorphism on small subtructures.
 Let $f:A  \to A' $ be an isomorphism between small subrings of $\Uu,\Uu'$. 
 By (0) of the Claim we may take $A,A'$ to be fields.  
Find $f_\G$ and $f_k$ compatible with $f$.  Then the Claim implies that $f$ extends
to an isomorphism  $\Uu \to \Uu'$.

Taking $\Uu =\Uu'$, the Claim gives the 
stable embeddedness (cf. \cite{CH}, Appendix.)  
Since any automorphism of $\G$ and any automorphism of $k$ lift,
the statements on strong orthogonality and the induced structure on $\G$ and $k$ follow.

 \eprf

If $A \leq M \models T$, $A(c)$ denotes the smallest substructure of $M$ containing $A,c$.
By {\em transcendence degree} of $B/A$ we mean the transcendence degree of $B \meet \VF$
over $A \meet \VF$.  

\<{remark}   In general, for $c \in \VF$, $A(c)$ can have infinite 
transcendence degree over $A$.  \rm (Let $\val(t)>0$ and consider $c=\sum a_i t^i$, with
$a_i \in k$.  Then each $a_i \in \Qq(t,c)$.) \>{remark}

Recall that $\RV = \VF / (1+ \Mm)$,   $\rv:  \VF \to \RV$ and
$\valr: \RV \to \G$ are the natural maps, 
and  $\RES$ is the subset of $\RV$
consisting of points whose image in $\G$ is definable.  $\RES$ is a strict Ind-definable
set, i.e. a union of definable sets.   

\<{lem} \lbl{rv1}  Let $\g$ be a definable point of $\G$, $V_\gamma =  \valr \inv (\gamma)$.
Then there exists a definable map $g: V_\gamma \to \VF$, with $\rv \circ g = id$.
    \>{lem}

\prf  For some $m$, and some definable $c \in \VF$,  we have $m \g = \val(c)$.
  Recall that $k$ is embedded in $\Oo$.  
    Let $A(\gamma) = \{y \in \VF: \val(x) = \g \}$,
 and $B(\gamma,c) = \{y \in A(\gamma): y^m \in c \k^* \}$.  
Then $\rv$ is injective on $B(\gamma,c)$; if $\rv(y)=\rv(y')$ then $y/y' \in \k^*$, and
$\rv(y/y') = \res(y/y') = 1$, so $y=y'$.  Thus the restriction of $\rv $  to $B(\gamma,c)$ 
defines a bijection $r: B(\g ,c) \to V_\g$, whose inverse 
is a section $V_\gamma \to B(\gamma,c) \subseteq A(\gamma)$.    
\eprf

Thus  not 
only $\k$ but also  $\RES$ admit a section into $\VF$:

\<{lem}  \lbl{rv2} $\RV$ is stably embedded, with the same induced structure as from
$ACVF$.  \>{lem}

\prf  This can be seen by extending 
a given automorphism $f_{rv}$  of $\RV(\Uu)$, as 
 in \lemref{qe-loc}.   
 In Step 1, note that when $\val(\sum_i a_i c_i) = \min \val (a_ic_i)$,
 it follows that $\rv(\sum_i a_ic_i) = \sum_{i \in I_{min}} \rv(a_ic_i)$; where
 $I_{min}$ is the set of indices $i$ such that $\val(a_ic_i) $ is minimal, and where
 addition is defined naturally on elements of $\rv$ by $ex+ey = e(x+y)$, where
 $x,y,x+y \in k ^*$ and $e \in \rv $.  In (2) we choose an extension to the algebraic
 closure compatible with the given isomorphism on $\rv$; this uses the stable embeddedness
 and known induced structure of $\RV$ in $ACVF$.  Step (3) is the same, noting that
 the data determines $\rv(c-d)$ too.  Step (4) is identical; note that $\RV$ does not grow
 in immediate extensions.      \eprf
 
\ssec{Structure of definable sets} 

We take a further look at the structure of definable sets and raise some questions; the material here will not be needed for the main theorem, where only smooth functions will be used.

\<{defn} \lbl{normal-def} Let $T=T_{loc,A}$. 
 Let $X \subseteq \VF^n \times \G^l$   be a   definable set.  A 
 {\em normal form} for $X$ is an expression 
 $$X = g(X^*)$$ where $\tX  \subseteq \VF^{n+m} \times \G^l$ is an 
 $ACVF_A$- definable set, and $g$ an $ACVF_A$- definable function on $\tX$, 
 and $X^* = \tX \meet (\VF^n \times k^m) \times \G^l $, such that:

(N)  $g: X^* \to X$ is bijective.   \>{defn}

We will say that $X$ is normal {\em via} $\tX$, or via $g$.  
  
\<{lem}  \lbl{normal1}  (1) 
If $X \subseteq \VF^n \times \G^l$ is the disjoint union of two definable sets $X_1,X_2$, each having a normal form, then $X$ has a normal form.
\def\tD{{\tilde{D}}}

(2)  Say $D \subseteq \VF^n$ is weakly normal via  $\tD \subseteq \VF^{n+m+l}$ if
$\tD$ is ACVF-definable,   and the projection $\pi': \VF^{n+m+l} \to \VF^n  $ induces
a bijection $\tD \meet (\VF^n \times \k^m \times \VF^l) \to D$.   Then a weakly normal 
set $D$ has a normal form, with $g$ a coordinate projection.  

(3)  If $D$ has a normal form, then it has one with $g$ a coordinate projection. \>{lem}

\prf  (1) By adding $0$'s we may assume $X_i$ is normal via $\pi: \tX_i \to X_i$ with $\tX_i 
\VF^{n+m} \times \G^l$.  Let $\tX_i' = \tX_i \times \{i\}$ (with $i \in \k$) and $\tX= \tX_1 \union \tX_2 \subseteq \VF^{n+(m+1)} \times \G^l$.  Then $X$ is normal via $\tX \to X$.
 
 (2)   Let $D'$ be a quantifier free formula equivalent in $ACVF_A$
to the projection of $D^*$ to $\VF^{n+m}$.  Then $D$ is normal via $D'$; and if 
$g$ is a function on $D'$, then $g$ is $ACVF_A$-definable iff $g \circ \pi''$ is  $ACVF_A$-definable. 

(3)  Say $X=g(X^*)$, as in \defref{normal-def}.  Let $X'' = \{(g(x),x): x \in \tX \}$.
Then $X$ has weakly normal form via $X''$.  By (2) it has normal form using a projection.

 \eprf

\<{lem} \lbl{normal2}  (1)  Every definable subset of $\VF^n \times \G^j$ has a normal form.

(2)if $X$ is a definable set and $f: X \to \VF$ a definable function,
there exists $\tX$ with $X$ normal via $n: \tX \to X$, and an $ACVF_A$-definable
function $F$ on $\tX$ such that $F$ agrees with $f \circ n$ on $\tX \meet X^*$.

\>{lem}

\prf  We
add to the language a new function symbol $D(x,y)$, defined to by $\frac{x}{y}$
when $y \neq 0$ and $0$ when $y=0$.  Also  replace $\res(\frac{x}{y})$ by a unary
function $R$, defined to be $i \circ \res$ on $\Oo$, and elsewhere $0$.  
Thus $\res(\frac{x}{y})$  is now a composed term, $R(D(x,y))$.  

To prove (1),  let $\phi$ be a formula in variables $\VF^n \times \G^j$.  
The $\G$ coordinates play no role, and we will ignore them to simplify notation.

If $\phi$ does not involve $R$ at all, it is already an ACVF formula.  Otherwise, 
$R $ occurs in some term in $\phi$, and the innermost occurence must have the form
$R \circ t$, with $t$ a rational function (more precisely, a term using $+,\cdot,D$.)

Let $\phi'(x,y)$ be the formula obtained from $\phi$ by replacing this instance of $R \circ t$
by $y$, and adding a conjunct:  
$$(t \notin \Oo \implies y=0) \& (t \in \Oo \implies \val(y -t )> 0)$$
  Then $\phi'(x,y) \& (y \in k)$ projects
bijectively to $\phi$.  By induction, there exists an ACVF formula $\phi''(x,y,y')$
such that  the projection $(\phi''(x,y,y') \& y' \in \k^l) \to \phi'$ is bijective.
 
 Then
 the projection $(\phi'' \& (y,y') \in \k^{l+1}) \to \phi$ is also bijective, and shows that $\phi$
 is normal. 
 
  For (2), apply (1) to the graph of $f$.  We obtain a partition of $f$ into sets $Y_i$,
  and ACVF-definable $\tilde{Y}_i  \subseteq \VF^n \times \VF \times \VF^m$,
  such that $\tilde{Y}_i \meet  (\VF^n \times \VF \times \k^m)$ projects injectively
  onto $Y_i$.   Let $U_i = \{(x,z) \in \VF^n \times \VF^m: (\exists ! y)(x,y,z) \in \tilde{Y}_i \}$.
  Then $U_i$ is $ACVF_A$-definable, and 
   $\tilde{Y}_i \meet  (\VF^n \times \VF \times \k^m)$ is contained in the pullback 
   $\tilde{U}_i$ of  $U_i$.  Replacing    $\tilde{Y}_i $ by $\tilde{Y}_i  \meet \tilde{U}_i$,
   we may assume the projection $\tilde{Y}_i \to \VF^n \times \VF^m$ is injective.
 Let $\tX_i$ be the image of this projection.  Then the projection 
 $n_i: \tX_i \meet (\VF^n \times \k^m) \to X$ is injective; let $X_i$ be the image.  
 The composition $f \circ n_i$  is   $ACVF_A$-definable, since it is the 
 section of the injective map $\tilde{Y_i} \to \tX_i$.     \eprf

\<{cor} \lbl{normal3} Any  definable set $X$ admits a definable map $\xi: X \to \res(\VF)^*$, whose fibers are  
 definable by valued field formulas; $\xi \inv(c)$ is $ACVF_{F(i(c))}$-definable.  \>{cor}
 
\prf  Let $g: \tX \to X$ be a normal form for $X$, with $\tX \subseteq X \times \VF^m$.
Let $\pi: \tX \to \VF^m$ be the projection, and let $h: X \to X^*$ be the inverse of the injective map 
$g | X^*$.  Let $\xi(x) = \res \pi h (x)$.  Then $\xi \inv(c) = g ( \pi \inv( i(c)) )$.  \eprf

\<{defn} \lbl{vf-dimension} A definable $X$ has {\em $\VF$-dimension $\leq n$ } if there exists a definable
$f: X \to \VF^n$ whose fibers are internal to $\RV$. \>{defn}

Unlike the case of pure Henselian fields of residue characteristic zero, the Zariski closure of $X \subseteq \VF^m$ can have larger $\VF$-dimension than $X$; for instance $\k$ is Zariski dense, of $\VF$-dimension zero.  

We will use the valuation topology on $\VF$, the discrete topology on $\RV$, and the product topology on products.

\<{lem}  \lbl{dimension}  Let $X$ be a definable subset of $\VF^n$.  Then the boundary of $X$ has dimension $<n$.
\>{lem}

\prf  We show that outside of a set of dimension $<n$, every point of $X$ is an interior point.  An $\RV^*$-union of sets of dimension $<n$ still has dimension $<n$, so we may   fiber over $\RV^*$.  By \corref{normal3}, we can take
$X$ to be defined by a valued field formula; but then the statement follows from the known one for ACVF.  The same  applies to the complement of $X$, so almost every point is interior either in $X$ or in the complement. \eprf

\ssec{Integration}

In this subsection, we assume residue characteristic zero.  We discuss an integration
theory for more general sets than test functions.  These results will not be required
for the proof of \thmref{A}.

 Recall the   measure-preserving   bijections of \cite{HK}.  These are   definable bijections between subsets $X \subseteq \VF^n \times \RV^m$.  We repeat the definition here, 
with the difference that we explicitly allow any bijection on the ``discrete'' $\RV$ sort.  

Let $\RV^*$ denote $(\RV \union \G)^m$.  (For some   unspecified $m$.)

 \<{defn}  \lbl{fc0}  Fix $n$.  An {\em elementary admissible transformations} of $\T_A$ is a
 $\T_A$-
 definable function of   one of the following types:
 \begin{enumerate}
  \item Maps 
  $$(x_1,\ldots,x_n,y_1,\ldots,y_l) \mapsto (x_1,\ldots,x_{i-1},x_i+a,x_{i+1},\ldots,x_n,y_1,\ldots,y_l)$$
with $a=a(x_1,\ldots,x_{i-1}, 
y_1,\ldots,y_l): \VF^{i-1} \times \RV^* \to \VF$
an $A$- definable function of the coordinates $y,x_1,\ldots,x_{i-1}$. 
 
  \item Maps $(x_1,\ldots,x_n,y_1,\ldots,y_l) \mapsto  (x_1,\ldots,x_n,y_1,\ldots,y_l,b(x,y))$
where $b: \VF^n \times \RV^{*} \to \RV^*$ is any definable function.  
\end{enumerate}    \end{defn}
 
A function  generated by elementary  admissible transformations
 over $A$ will be called an {\em admissible  transformation} over $A$.    
   \def\RVvolg{{\RV}_{\Gamma {\rm -vol}}}
 \def\RVni{{  \RV[n,\idot]}}

\def\CRV#1{{C\RV[#1]}}
Recall also the category ${{\RVvolg[n,\cdot]}}$ of definable subsets $X$ of $\RV^*$,
along with definable functions $f: X \to RV^n$, and the lift $\L$ to 
the category of definable subsets of $\VF^* \times \RV^*$ of dimension
$\leq n$.  Morphisms are definable maps.   The lift of $(X,f)$ is just
$X \times_{f,\rv} (VF^*)^n$.  The  $\cdot$
indicates that no bound is placed on the dimension ``discrete" $\RV^*$-component
of $\VF^* \times \RV^*$, or correspondingly on the fibers of $f$.   We will omit this  symbol, as well as the letters referring to $\G$-indexed volume since only such 
volumes will be considered, and write $CRV[\leq n]$ for ${{\RVvolg[n,\cdot]}}$.
$\CRV{\leq n} $ is the direct sum of $\CRV{k}$ over $k \leq n$.
(cf. Def.   5.21; see also Theorem 8.29.)
Here we will only consider volumes indexed by $\Gamma$. 
 


\<{lem} \lbl{normal4} Every definable $X \subseteq \VF^n \times \RV^* $ 
is in admissible bijection with the lift $\lambda U$ of some object $U $ of $\CRV{\leq n}$.  
 \>{lem}

\prf  Let $\xi: X \to \res(\VF)^l$ be 
as in \corref{normal3}.  By \cite{HK}, for each $c \in \xi(X)$,
there exists an $ACVF_{A(i(c))}$-definable object $Y(c)$ of $\CRV{\leq n}$
 $f_c: X(c) \to \L Y(c)$.  Note that $Y(c)$ is actually $A_{\RV}(c)$-definable,
by stable embeddedness.  
Putting these together  (as in \cite{HK} Lemma 2.3 ), 
we obtain an object $Y$ of $\CRV{\leq n}$, and an admissible bijection 
$X \to U$.  \eprf

Let $E_\mu$ be the equivalence relation on definable sets generated by the following steps:

1)  If there exists an admissible $f: X \to Y$, then $(X,Y) \in E_\mu$. 

2)  Let $X,Y$ differ only by a set of $\VF$-dimension $<n$.  Then $(X,Y) \in E_\mu$.

\<{cor} \lbl{normal5}  Every 
 definable $X \subseteq \VF^n \times \RV^* $ 
is   $E_\mu$-equivalent to a lift $\L U$ of some object $U $ of $\CRV{n}$.   \>{cor}

\prf The lift of objects of $\RV[k]$ for $k<n$ has $\VF$-dimension $<n$, and can be discarded. \eprf

Let $K_+(T_A)$ be the semi-group of $E_\mu$-classes of $T_A$-definable subsets 
of $\VF^n \times \RV^*$.

\<{cor} \lbl{n6}  There exists an isomorphism $ K_+(T_A) \to K_+(\CRV{n})/I$, for
some congruence $I$.  \>{cor}

\<{question}  Determine generators for $I$.   Do the generators of   $I_{sp}$ of 
\cite{HK} suffice?
\>{question}

Let $Res$ denote the residue field.

\<{lem}  If $X,Y$ are definable subsets of $Res ^n$, and $([X],[Y]) \in I$, $X,Y$
defined over a finitely generated domain $R \subseteq A$, then for some $0 \neq d \in R$, for all
homomorphisms $h: R[d \inv] \to F_q$ into a finite field, the $F_q$-varieties $X_h=X \tensor _R F_q, Y_hY \tensor _R F_q$
have the same number of points in $F_q$.  \>{lem}

\prf 
Otherwise, embed $R$ into the ultraproduct of the $F_q$;
interpret $X(F_q)$ as the volume of $\L X$; the fact that $([\L X] ,[\L Y]) \in E_\mu$ implies
that (for some $d$)  for   all such $h$, $\L X_h, \L Y_h$ have the same measure.   Here $\L$ and $E_\mu$
are interpreted in $F_ q((t))^a$ with the natural section $F_q \to F_q((t))^a$.  \eprf

Assume $\val(A) = \Zz$.
By \lemref{dimension}, any definable subset of $\VF^n$ has a well-defined Kontsevich motivic integral,
with values in the completion $A[[L \inv]]$ where $A$ is the Grothenieck ring of varieties over $A_{\res}$,
and $L$ is the class of the affine line.  
This integral is linear on constructible functions.  Moreover the power series obtained represent
rational functions, by the proof of Denef-Loeser.

\ssec{Valued fields over a curve}  
\lbl{valued fields over a curve}
Let $\f$ be a field, and $C$ a  smooth curve over $\f$.  

We describe here a first-order theory $\T=ACVF_{C;\f}$ convenient for adelic work.    It has the following sorts.

$k$ - an algebraically closed field with a distinguished field of constants $\f$.  $k$ is
endowed with the language of $\f$-algebras.  

$C(k)$ (when $C$ has a distinguished point, or genus $\geq 2$, we can take $C \leq \Pp^n $ so a special sort is not  necessary, but we take one nonetheless.)  

$\G$ - an ordered Abelian group, with distinguished element $1>0$.  
 
$VF$.  This sort comes with a map $VF \to C(k)$; the fibers are denoted $VF_x$.
Each $VF_x$ comes with valuation ring $\Oo_x$, a surjective homomorphism 
$\res_x: \Oo_x \to k$, and a ring embedding $i_x: k \to \Oo_x$, such that $\res_x \circ i_x = Id_k$.  Also, a map $v_x: VF_x \setminus \{0\} \to \G$, denoting a valuation with valuation ring $\Oo_x$.   For any variety $V$ over $\f$, we obtain using $i_x$
a variety over $VF_x$; let $V(VF) = \union _{x \in C(\k)} V(VF_x)$, a set fibered over $C(\k)$.  

We identify $k$ with its image $i_x(k)$.  

As a final element
of structure, we have a function $c: C(k) \to C(VF)$, such that $c(x) \in C(VF_x)$; and 
for any   $f \in k(C)$,  $\val f(c(x)) = ord_x (f) \cdot 1$.   

Technically, the above depends on a specific chart for $C$ as an abstract algebraic variety over $\f$.   We take $C$ to be a complete curve.
If $C$ is given as a union of open subvarieties $C_i$ embedded in $n_i$-dimensional affine space, $i=1,\ldots,m$, then $c$ is 
by $n_i$-tuples  $c_i^j$ of functions $c_i^j: C(k) \to VF$; the theory will state the natural
compatibilities, and up to obvious bi-interpretation will not depend on the chart. 

We note that $\G$ serves as the value group of each of the valued fields $VF_x$.
This will be important later, for instance when considering divisors on $C$ of degree $0$.
The identification of the various value groups

\<{lem} \lbl{vfc1}  $\T$ is complete.  $k$ and $\G$ are stably embedded.
The induced structure on $k$ is the $\f$-algebra structure.
The structure on $\G$ is the ordered group structure, with distinguished element $1$.

Moreover $\T$    has quantifier-elimination.  
 \>{lem}

\prf  Let $\Uu,\Uu'$ be saturated models of $\T$ of the same cardinality.  We wish to show that $\Uu,\Uu'$ are isomorphic.  Since $k \cong k'$ and $\G \cong \G'$ we may assume they
are equal.  In particular $C(k)=C(k')$.  Now for any $x \in C(k)$ we have $VF_x$ and $VF'_x$.
The structure $(VF_x,k,\G,c(x)) $ is a saturated model of $T_{loc}$, and
so is $(VF'_x,k,\G,c'(x))$.  Moreover by quantifier-elimination, $tp_{T_{loc}}(c(x) / k,\G) = tp_{T_{loc}}(c'(x)/ k,\G)$, since the valued fields $k(c(x)),k(c'(x))$ are both isomorphic
to $k(C)$ with the valuation corresponding to the point $x$.  Hence by stable embeddedness
of $k,\G$ there exists an isomorphism $f_x: VF_x \to VF'_x$, with $f_x(c(x)) = c'(x)$.   Putting together these isomorphisms we obtain an isomorphism $f: \Uu \to \Uu'$.  

The same proof shows that any automorphism of $k, \G$ extends to an automorphism of $\T$.
Hence they are stably embedded and their induced structure is as stated.   In the same way we can extend
partial isomorphisms, hence quantifier elimination.  
 \eprf

 See Appendix 1 for the notion of ``piecewise definable''.

\<{lem} \lbl{piecewise-g}  $k(C)$ is piecewise definable over $k$.  
\>{lem}

\prf  For $\Pp^1$ this is completely elementary.  The elements of $k(\Pp^1)$ 
can be identified with pairs $(f,g)$ of polynomials, relatively prime, and with $g$ monic
or $g=0, f=1$.  For any given bound on the degrees, this is clearly a constructible set.  
Given a pair $f,g$ of polynomials
of degree $\leq n$, the valuation at $0$  of $f/g$ is bounded between $-n$ and $n$,
and each of the possible values is constructible.  Moreover $PGL_2(\k)$ acts constructibly
on the polynomials of degree $\leq n$ and on the set of valuations, and for $\phi \in PGL_2(\k)$,  the valuation 
of $f/g $ at $c=\phi^{-1} (0)$ equals the valuation of $f^\phi / g^\phi$ at $0$.

For other curves,  the proof is more easily carried out using 
 quantifier elimination for  $T^{loc}$.   We   show that when $k \models ACF_\f$, the structure
  $k(C)$ in the language $L^g_\f$  is piecewise definable over $k$.  Since the induced structure on $k$
  from $ACVF_{C;\f}$ is the $\f$-algebra structure, it suffices to piecewise interpret $k(C)$ in a model 
  of $ACVF_{C;\f}$, provided that the interpreted copy of $k(C)$  is contained in $\dcl_M(\k)$.

The field $k(C)$ as an $\f$-algebra was treated in \lemref{limited-1} (2); moreover the proof
there shows that $k(C)$ as a differential algebra, i.e. with the additional ternary relation $df = h dg$, is also piecewise definable over $k$.

 We have to show that on each limited subset $Y$ of $k(C)$, and each $n$,  the relation
$V_n$
restricted to $Y \times C(k)$ is definable.   We may take $C$ to be embedded
in $\Pp^m$.  For some $d$, each element   $f$ of $Y$ is a quotient of two homogeneous polynomials of 
degree $d=d(Y)$.  Now $V_n(f,\alpha)$ holds iff $ord_\alpha (f) = n$
iff $\val (f(c(\alpha))) = n \cdot 1$.  Hence $V_n \meet (Y \meet C(k))$ is definable
in   $M$, over $\f$.  By quantifier-elimination in $ACVF_{C;\f}$ and the 
fact that the induced structure on $k$
  from $ACVF_{C;\f}$ is the $\f$-algebra structure, $V_n$ is definable.

We now define ${\bf t}$.   Let $t_0 \in F(C)$ be non-constant.
For any $\a \in C(k^a)$ such that $t_0: C \to \Pp^1$ is unramified at $\a$, we let $t_\a=t_0 - t_0(\a)$.
There are finitely many values of $\a$ where $t_0$ is ramifed; there we make some choice
of parameter $t_\a$, in a Galois invariant way.  Let ${\bf t} = \{(\a,t_\a): \a \in C(k^a)\}$.

Given again a limited subset $Y$, we need to compute
$res_a(f \omega)$, uniformly in $f \in Y$ and  $\a \in C(k)$.  Using the definable parameter $t_\a$ and the differential structure, we find the unique $g \in C(k)$ with $f \omega = g dt_\a$.
Now $g$ lies in some limited definable set $Y'$.  For some $n$, for each $y \in Y'$,
we have $v_\a(g -(  c_n(y) t_\a ^{-n} + \ldots + c_{1}(y) t_\a^{-1} )\geq 0$ 
for some (unique) $c_n(y),\ldots,c_{1}(y) \in \k$.  Then $res_a(f \omega) = c_1(g)$.  
 \eprf

\ssec{Adelic definable sets.}  \lbl{adelic}
 
 Let $C$ be a curve over a field $\f$, $F=\f(C)$.  
 We continue working with the theory $ACVF_{C;\f}$;  definability relates to this theory.   For a subset of $\k^n$, it is equivalent to $ACF_\f$-definability, while for a subset of $VF_v^n$ (with $v$ a point of $C$) it is equivalent to $T^{loc}_{\f(v)}$-definability.  
These theories
admit quantifier-elimination; when we use a quantified formula, we mean the quantifier-free equivalent.  The witness is not assumed to exist rationally
over $F$.

\<{example} \rm \lbl{x-1}  We will define algebras $D_v$ with subrings $R_v$ for each $v$, and consider $R_v^*$-conjugacy classes.
On $D_v(F)$-points, this will not be the same as $R_v^*(F)$-conjugacy.  Moreover,
if $\dot{D_v},\dot{R_v}$ is another such pair, we define a point $x$ of $D_v(F)$ 
to be conjugate to a point $\dot{x}$ of $\dot{D_v}(F)$ if there exists an isomorphism
$h: D_v \to \dot{D_v}$ with $h(R_v)=\dot{R}_v$ and $h(x)=\dot{x}$.  We are interested
only in the case where $D_v(F),\dot{D}_v(F)$ are non-isomorphic division algebras;
so such an isomorphism never exists rationally; the notion of integral  conjugacy across forms will nevertheless be useful.  \>{example}     %

Let $X$ be a definable set, $V$ a variety over $F$.   By a {\em definable function} $f: X \to \Pi_{v \in C} V({VF}_v )$ we mean
a definable function $f$ on $X \times C$, such that $f(x,v) \in V({VF}_v)$.   We   view
$f$ as a function into the product.  We will only consider the case that $X \subseteq \dcl(\k)$.  

Let $\phi$ be a formula in the language of $ACVF_{C;\f}$ {\em enriched with 
additional unary function symbols $\xi_1,\ldots,\xi_n: C \to   \VF$.}  Then $\phi$ defines a subset 
$Z= \phi(\Pi_{v \in C} {VF}_v^m) := \{ z \in \Pi_{v \in C} {VF}_v ^n: \phi(z) \}$.  
Note that $Z$  {\em definable-in-definable-families}, i.e.  $f \inv(Z)$ is a definable
subset of $X$ whenever   $f: X \to \Pi_{v \in C} {VF}_v ^m$ is definable.
Similarly if $Z$ is defined by an infinite disjunction of formulas $\phi_k$, then 
$Z$ is     Ind-definable-in-definable-families.  This   extends  
  to subsets of $\Pi_{v \in C} V({VF}_v)$ where $V$ is any variety over $F$.
  

Assume $V$ is an affine variety, or at any rate that a subset $V(\Oo_v)$ of integers
is given in some way for each $v$, and $V(\Oo_v)$ is definable uniformly in $v$.  
Then $ \Pi_v V(\Oo_v)$ is  definable in the above sense (by the formula 
$(\forall v \in C)(\xi(v) \in \Oo_v)$), and so if $f: X \to \Pi_v V({VF}_v)$ is definable
then $f \inv ( \Pi_v V(\Oo_v)))$ is a definable subset of $X$.

\<{example} \lbl{orbits} \rm Let $D_v,R_v$ be as in \exref{x-1}, and write $R^*=\Pi_v R_v$.   Let $a \in D(\Aa)$ be definable.  The $R^*$-orbit $O$ of $a$ is then definable; for $x \in D$, 
$$x \in O \iff (\forall v)(\exists u \in R_v^*)(u a u \inv = x)$$  In particular, 
$O(K):= O \meet D(K)$ is Ind-definable.

If $R_v$ is a bounded subset of $D_v$ for all $v$, then 
  $O(K)$ is a limited subset of $K$; hence by \lemref{vfc1} it is definable over $\k$.

 In practice we will have $R_v$ bounded only  for all $v\neq v_0$; in this case $O(K)$ is not limited, but is nearly so, and we will still be able to assign to $O(K)$ a class in 
 an appropriate quotient of the Grothendieck ring.   
 
Of course one can have a  definable $R^*$-orbit with no definable element, and the
same considerations hold.  
\>{example}

 \<{example} \rm
Let $K_v^h$ be the Henselization within $VF_v$ of $K=k(C)$.  
This is an Ind-definable set
with parameter $v$:     In the notation of \secref{valued fields over a curve}, an element of $K_v^h$ has the form $h(a,c(v))$ for some definable function $h$, and an $n$-tuple $a$ from $k$.  

Note that 
%
any definable function $f: X \to \Pi_v VF_v$ has image contained
in $\Pi_v K_v^h$.  

 The adelic points of $\Pi_v V( VF_v)$ are defined as the union over all $m$ of the
 set $ \union_{|w|=m} \Pi_{v \in w} V(VF_v) \times \Pi_{v \notin w} V(\Oo_v)$.   
We let $V(\Aa)$ be the set of adelic points of $\Pi_{v \in C(\k)} V(K_v^h)$.  
For any subset $Y$ of $C(k)$ we also let $V(\Aa_Y) $ be the set of adelic points of
  $\Pi_{v \in Y} V(K_v^h)$.  
We saw that any definable function on $X$ into the adelic points maps into $V(\Aa)$. 
 
Let $f: X \to \Pi_{v \in C} V({VF}_v)$ be a definable function.  If $f(x) \in V(\Aa)$  for any $x$ we say that $f: X \to V(\Aa)$.   By compactness, this implies the existence
of a fixed $N$, such that for each $x$, for some $w(x)\subseteq V$ of size $N$,
$f(x) \in \Pi_{v \in w(x)} V(K_v^h) \times \Pi_{x \notin w} V(\Oo_v)$; moreover by stable embeddedness of $\k$, if $X=W(\k)$ is a constructible set over $\k$, the map $x \mapsto w(x)$ can be taken to be constructible.

We will be interested in $G(\Aa) / H(\Aa)$ for certain congruence subgroups $H$;
since $K$ is dense in $K_v^h$, for our purposes $K_v^h$ and $K$ could be used
interchangeably, and we discuss $K_v^h$ only   to clarify the link with the classical presentation. 
\>{example}

 When $H$ is a definable subgroup of a group $G$, and $X \subseteq G$,
we write $X/H$ for the image of $X$ in $G/H$.

\<{example}  \rm \lbl{x-cycles}
    Let $G$ be a group scheme   over $F$, and $H_v = G(\Oo_v)$; or more generally 
assume $H_v$ is a   uniformly $ACVF_v$- definable subgroup of $G$. Write $H=\Pi_v H_v$.  Given a finite subset $s$ of $C$, let $(G/H)_s = \Pi_{v \in s} G(K_v^h)/H_v $; 
in our examples $H_v$ will be open in $G$, and we will have 
$G(K_v^h)/H_v  = G(K)/H_v$.   
  At all events, $G(K)$ or $G(K_v)$ are Ind-definable sets; write
$G(K) = \union_n G[n]$, with $G[n]$ definable.  Let $(G/H)_s[n] = \Pi_{v \in s} G[n]/H_v $.  Then  $(G/H)_s = \union_n (G/H)_s[n]$ is Ind-definable.   

When $s \subset s'$ we have a map $(G/H)_s \to (G/H)_{s'}$, mapping
$a \mapsto a'$ where $a'(v) = H_v/H_v$ for $v \in s' \m s$.  
Consider  the direct limit of $(G/H)_s$  over all finite subsets $s$ of $C$.  The 
directed set here is the set $P_\omega C$ of finite subsets of $C$; this is itself Ind-definable, limit of the definable sets $P_{\leq n} C $ of all $\leq n$-element subsets of $C$.  Since  
$(P_{\omega}C) (M)$ depends on the model $M$,  the direct 
limit $\lim_{s \in P_{\omega}(C)} (G/H)_s$ is not, as presented, an Ind-definable set.  However, for fixed $n$,  the disjoint union $\sqcup_{s \in P_{\leq n}(C)} (G/H)_s[n]$
is a definable set.  There is a natural isomorphism 
$\lim_{s \in P_{\omega}C} (G/H)_s = \lim_n \sqcup_{s \in P_{\leq n}C} (G/H)_s[n]$.  
The latter is Ind-definable, and we denote it $(G/H)(\Aa)$.   

 Note
that when $G$ is Abelian, $(G/H)(\Aa)$ is an Ind-definable group, but is not in general
an Ind-(definable group).  
A basic case:  $G=G_m, H=G_m(\Oo)$; in this
case $(G/H)(\Aa)$ is the group of cycles on $C$.   \>{example}

\<{example}  \lbl{x-tori-1} \rm  Let $V=T$ be a multiplicative torus.  We have a uniformly definable homomorphism 
$T(VF_v) \to X_*(T) \tensor \G$; taking sums we obtain a homomorphism $\val_T: T(\Aa) \to X_*(T) \tensor \G$; the kernel is denoted $T(\Aa)^0$.  Then $T(\Aa)^0$ is Ind-definable.  \>{example}

\<{example} \rm Let $C$ be a curve, $S$ a nonempty finite definable subset.
  Let $I=  \Pi_{v \in C} G_m(K_v^h)$, $T_S = \Pi_{v \in S} G_m(K_v^h) \times \Pi_{v \notin S} G_m(\Oo_v)$.
  We also have a diagonal embedding of $K^*$ into $I$.  Then $I/K T_S$ is   weakly representable if and only if for any $a,b \in S$, the image of $a-b$ in the Jacobian of $C$ is a torsion point.  In particular, if $|S|=1$ then $I/K T_S$ is  weakly representable.
    
  Indeed for checking weak representability we may change the base so that $C$ has a rational point $0$, in fact we can take $0 \in S$.  In this case the Jacobian $J$ can be identified with $I^0 / K^* \Pi_{v \in C} G_m(\Oo_v)$
  where $I^0$ is the group of ideles of degree $0$. 
We have a natural   homomorphism $J \to I / K T_S$.  It is surjective, since $I_0K_0 = I$.  The  kernel generated by the points $a-b$ with $a,b \in S$.  If one of these points is not torsion, then condition (2) of
\lemref{wrep} will not hold.   If all points $a-b$ are torsion then the kernel is finite, and (1,2) of \lemref{wrep} are clear.
  \>{example}

\<{example}    \lbl{x-tori-2}  \rm 
  The double-coset equivalence relation $\TR x T(K) = \TR y T(K)$ is
 Ind-definable.   In this case for definable $t,t': C \to T(\Aa)$ we have
$t E t'$ iff $(\exists k \in T(K)) (t \inv t k \in \TR)$  which is Ind-definable.  The formula inside
the quantifier is definable: $(\forall v \in C)(y \in \TR_v)$. 

 Let $\CT = T(\Aa) / \TR$ as in \exref{x-cycles}.  It follows that the 
 embedding of $T(K)$ in $\CT$   is Ind-definable.
 
 When $T=G_m$, and $\TR = G_m(\Oo_v)$, $\TR \bs T(\Aa)^0 / T(K)$   is the Jacobian
 of $C$.   With more general $\TR$ one obtains    Rosenlicht generalized Jacobians of $C$.  If $T=R_{C'/C} G_m$ is obtained from a cover $C'$ of $C$ by reduction of scalars of $G_m$, this is a generalized Jacobian of $C'$.   \>{example}

 \>{section} 

 \<{section}{Division rings}

\ssec{Adelic structure of cyclic division rings.}  \lbl{norm-s}

$\f$ is a perfect field.  $\k$ is a model of $ACF_\f$.  We denote $F=\f(t)$, $K=\k(t)$.
This section is purely algebraic, and adeles (or repartitions), when they are mentioned,  are treated essentially classically. 

When only one valuation is involved, we denote the valuation ring by $\Oo$.
 When many  valuations $v$ are involved, we denote the Henselization of  $F$ as a valued field by $F_v$, the valuation ring by $\Oo_v$, and let $\Oo = \Pi_v \Oo_v$.

Let  $L$ be a commutative semi-simple algebra defined over $\f$, of dimension $n$,
with $Aut_\f(L)$ a cyclic group; 
  let $g$ be a generator.   Let 
$$ D_{g,t} =  L[s] / (sa=g(a)s, s^n=t)$$
We view $D_{g,t}$ as an $ACF_{\f(t)}$-definable
algebra.  For most of the discussion, we fix $g$ and denote $D=D_{g,t}$.   
Eventually we will compare $D$ to another form 
  $\dot{D}=D_{\dot{g},t}$ 
 with $\dot{g}$ another generator of $Aut(L)$.  

 We are mainly interested in the case that $L(\f)$ is a cyclic Galois field extension
of $\f$; in this case we let $\lf = L(\f)$, see \secref{lf}.  As we will see, $D(F)$ is then  a division algebra over $F=\f(t)$.  

Let $d_1,\ldots,d_n \in L(\f)$ be an $\f$-basis
for $L(\f)$.  Then $(d_is^j: 1 \leq i \leq n, j \in \Nn)$ is a  basis for $D$ over $K$.

   Let $v_0,v_\infty$ be the valuations of $\k(t)$ with $v_0(t)>0$, $v_\infty(t)<0$, and let
   $v_1$ be the valuation of $F=\f(t)$ over $\f$ with $v_1(t-1)>0$.

Let $N=  R^* \meet D^*(K)$.   Note that $\bL$ is a normal subgroup of $N$,
(indeed $N$ is the normalizer of $\bL$ in $D^*$.)

Over $\lf(t)$ there exists a representation of $D$ on $L$; namely $L$ acts diagonally,
while $s$ is the product of the permutation corresponding to $g$ with a diagonal matrix
$(1,\ldots,1,t)$.  This defines an isomorphism of $D$ with $M_n$ over $\lf$.

 For any valuation $v$ on $\f(t)$ with $v(t) = 0$, we define a subring $R_v$ of $D$:  
 
 \<{statement} \lbl{S}    For  $v(t)=0$, for 
 $x_{ij} \in \k(t)_v$, $i,j=0,\ldots,n-1$ we have:

 $ \sum x_{ij} d_j s^i \in R _v  \iff  \bigwedge_{i,j=0}^{n-1} x_{ij} \in \Oo_v  $
\>{statement}

When $v(t) \neq 0$, we let $R_v = D$.  On two occasions we will refer to the ring defined by the same formula
\ref{S} at $v=0$; but in this case we will denote it $S_0$.

The family of rings is uniformly definable in the theory $ACVF_{\f,C}$ of \secref{qeSec}.
  We refer to this choice of 
subring, for each valuation  $v$, as the adelic structure.   

 We write $\Oo_\Aa$ for $\Pi_v \Oo_v$, and let $R=\Pi_v R_v$.   Denote the $n \times n$ matrix algebra by $M_n$.   Let $\bM_v = \Mm_v R_v$.  Let $Z$ be the center of $D^*$.
 Write $\bL=L(\k)$.
 
\<{lem}   \lbl{adeldiv}  \begin{enumerate}
  \item  For $v(t)=0$,  $(D,L,R_v) \cong_{\lf(t)} (M_n,L,M_n(\Oo_v)$.    This characterizes $R_v$ uniquely up to  conjugacy by \lemref{r3} (3). 
 
   \item  $R(K) = \bL[s,s\inv]$.    If $x \in D(F)$ and $x \in Z R_v^*$ for each $v$
   with $v(t)=0$, then $x \in Z(F) R(F)$. 
   
    When $L(\f)$ has no 0-divisors, 
we   have $L[s]^*(\f) = L^*(\f)$, and $R(F)^* = \bL[s,s \inv ]^* (\f)= L^*(\f) s^{\Zz} = N$.

\item  a)       Let $v$ be such that $v(t)=0$.  
 Let $F'$ be a Henselian valued field extension of $(F,v)$ with residue field $\f'$.  
 Assume (*):   $D_{g,\bt}(\f) \cong M_n(\f)$.   Then $(D,R_v)  \cong (M_n,M_n(\Oo))$ by an $ACVF_{F'}$-definable isomorphism.   \\
b)   If $v(t-1)>0$,  or if $\f'$ has a unique (and Galois) field extension of order $n$, then (*) holds. \\

\item  Fix $v$, and denote by $\bt$ the image of $t$ in $R_v/\bM_v$.  
  If $v(t)=0$, then $R_v / \bM_v   \cong D_{g,\bt}  $. 
  Assume $v(t) > 0$, and $L(\f)$ is a field.  Then $D(F_v)$ is a valued division ring, and  $s^\Zz U = D(F_v)^*$ where
$U = \{x: v(x)=0 \}$.

\item  Assume  $v(t)=0$, $n$ prime.
Let $\tau \in F_v$ be a uniformizer, i.e. $v(F_v) = \Zz v(\tau)$.    
  If $D_{g,\bt}(\f)$ is a division ring,   
   then any nonzero element of $D(F_v)$ can be written
as $z+ \tau^m y$, where $z \in F_v, m \in \Zz$, $y \in R_v$ and the residue $y+\bM_v$ of $y$ is regular semisimple. 
\item     Assume $L(\f)$ has no zero divisors.
Then every left ideal of
 $\bL[s]$ is principal.

\end{enumerate}

\>{lem} 

\prf      
(1)  
  Over $\lf$ we have a definable
 basis $e_1,\ldots,e_n$ of $L$ consisting of idempotents, such that $g(e_i) = e_j$ (where $j=i+1 \mod n$.)  The change of basis from $(d_i)$ to $(e_i)$ is effected by a matrix in $GL_n(\lf)$ 
 and so it does not effect the adelic structure.  We have an $\lf$-definable isomorphism $D \to GL_n$, mapping $e_i$ to the standard
 matrix $e_{ii}$ (with $1$ at $(i,i)$ and zeroes elsewhere), and mapping $s$ to the product of 
 the cyclic permutation matrix, with the diagonal matrix $(1,1,\ldots,t)$.   This is an isomorphism
 between $D$ and $GL_n$;   when $v(t)=0$, it is straightforward to verify that the ring $R_v$ defined in  \eqref{S}
 maps to $GL_n(\Oo)$.

(2)
Let $x=  x_{ij} d_j s^i  \in D(k(t)) $, and assume $x  \in  R_v$ for all 
 $v \neq 0,\infty$.  Then $v(x_{ij}) \geq 0$ for all such $v$, so $x_{ij} \in \k[t,t \inv]$.  Hence
 $ \sum x_{ij} d_j s^i \in L[s,s\inv]$.  
 
Let $x \in D(F)$ and suppose for each $v$ with $v(t) =0$ we have
 $x = z_v r_v$ with $z_v \in Z$, $r_v \in  R_v^*$.   Let $a=\det(x) \in F$.
 Then $v(a)=v(  z_v^n) $  for each $v \neq 0,1$; so $v(a)$ is divisible by $n$
 in $\val F_v$, for each such $v$, and we can find $b \in F$ with $n v(b) = v(a)$
 for each such $v$.  Dividing by $b$ we may assume $v(a)=0$ for $v \neq 0,1$.
 So $v(z_v)=0$, so $z_v \in \Oo_v^*$ for each such $v$; but then $x \in R_v^*$. 
  
 
When $L(\f)$ has no 0-divisors, we  
  have $L[s]^*(\f) = L^*(\f)$, and $L[s,s \inv ]^* (\f)= L^*(\f) s^{\Zz}$:  it suffices to see
 that if $f,g \in L[s](\f)$ and $fg=s^n$ then $f,g \in L s^\Nn$; this is clear by viewing $f,g$
 as non-commuting polynomials in $s$; the product of a non-monomial with a polynomial is always non-monomial, by considering lowest and highest terms.

(3) 
(b)  Note that $D_{g,1} = M_n$ over any field, with the integral structure  \eqref{S} coinciding
with  $M_n(\Oo_v)$.  
 
a)  Since $D_{g,\bt} (\f) \cong M_n(\f)$, we have $\bt   = N_{L(\f)/\f} (\bc)$ for some $\bc \in L(\f)^*$.
(Proof: find an element $s'$ normalizing $L$ such that conjugation by $s'$ has the effect of $g$ on $L$, and such that $(s')^n =1$; such an element exists in the matrix ring.  It follows
that $\bc=s \inv s' $ centralizes $L$, so $s \inv s' \in L$; and we compute $N(\bc)=s^n = \bt$.)  
Let $h(X) = X^n + \ldots \pm \bt $ be the minimal polynomial of $\bc$ over $\f$.  Since
$\f$ is perfect, $h$ must have some nonzero monomials of degree strictly between
$0$ and $n$.  So $h'(\bc) \neq 0$.
Lift $h$
to $H (X) = X^n + \ldots \pm t \in \Oo[X]$.  Using  
Hensel's lemma find $c \in L(\f_v)$ with $H(c)=0$.  This is an invertible element of $R_v$, and
  $t = N_{L(F_v)/F_v} (c)$.  Dividing by $c$ we reduce to the case of $D_{g,1}$.  
  
If $v(t-1)=1$, i.e. $\bt=1$, we may take $\bc=1$.  If the residue field has a unique (and Galois)
extension of order $n$, then the norm map from this extension to $\f'$ is surjective.

(4)      Clearly $L(F_v) = L F_v$ is a field, namely an unramified extension
of $F_v$; we view it as a valued field.    Any nonzero element of $a \in D(F_v)$ may be written uniquely as $\sum_{i=0}^{n-1} u_i s^i$ with $u_i \in LF_v $; 
define $v(a) = \min_i \{v(u_i) +  \frac{i}{n} \val(t) \}$.  
Note that this minimum is attained at a unique $i$, since 
the summands are distinct elements of $\Qq \val(t)$, even modulo $\Zz \val(t)$.   
From this it follows that $v(ab)=v(a)+v(b)$.  In particular $D(F_v)$ has no zero-divisors;
since $D$ is finite-dimensional over the center, 
it follows that $D(F_v)$ is a division ring.   We have $s^\Zz U = D(F_v)^*$ where
$U = \{x: v(x)=0 \}$.  

(5)
  Note first that if $\bD$ is an $n$-dimensional division ring over a perfect  field $\f$
  (or any field $\f$ if   $n \neq char(\f)$),  then any nonzero element $\ba$ of $E$ is either central or regular semisimple.  Indeed over  $\f^{alg}$
there is an isomorphism $\a: \bD \to GL_n$, and the set $s$ of eigenvalues
of $\a(\ba)$ does not depend on the choice of $\a$.  Since $n$ is prime, these eigenvalues   are all equal or all distinct.    If all are equal, say to $\gamma$, then $\gamma \in \f$  and $ \ba - \gamma$ is  non-invertible, hence equal to $0$.   

Note also that $\bM_v(F_v) = \tau R_v (F_v)$.
For $x \in D(F_v)$, let $v(x)$ be the unique $m$ with $\tau^{-m} x \in R_v \m \tau R_v$.
Since $R_v (F_v ) / \tau R_v(F_v))$ is a division ring, we have $v(xy)=0$ whenever
$v(x)=v(y)=0$, and it follows that $v(xy)=v(x)+v(y)$ in general. 

  Let $a \in D(F_v)$.  If $a$ is central there is nothing to prove.  Otherwise
 $ab \neq ba$ for some $b \in D(F_v)$.  For any central $z$ we have $ab-ba= (a-z)b-b(a-z)$,   so $v(ab-ba) \geq \min v((a-z)b),v(b(a-z)) = v(a-z) + v(b)$,
 or $v(a-z) \leq v(ab-ba) - v(b)$.  Thus there exists a central $z$ with $v(a-z)$ maximal; subtracting $z$,
 we may assume it is zero, i.e. $v(a) \geq v(a-z)$ for any central $z$.  Multiplying
 by $\tau^{-v(a)}$ we may further assume that $v(a)=0$.  Now $v(a-z) \leq 0$
 for any central $z \in R_v$;  as central elements of $\bD=R_v/\bM_v$ lift to central elements of $R_v$,
   the residue $\ba$ of $a$ cannot be central.  By the first paragraph
 we have $\ba$ regular semisimple.

(6)  
    The proof of the Euclidean algorithm   for the commutative polynomial ring works equally well for the 
 twisted polynomial ring $\lf[s]$, and so does the proof that every left ideal
 in a Euclidean ring is principal.
  \eprf

 \def\Oof{\Oo'}  \def\Aaf{\Aa'}
 \def\Rf{{R'}}

Let $\Oof = \Pi_{v \neq \infty} \Oo_v$, let $\Aaf_F$ be the  adeles over $F$ without the factor at $v=\infty$, and let $\Rf = S_0 \times \Pi_{v \neq 0, \infty} R_v$.

Let $\Rf(F)$ be the set of elements of $D(F)$, whose images in $D_v$
fall into $R_v$ for every $v$, {\em and} whose image in $D_0$ falls into $S_0$.    
The rings $R_v$, 
$\Rf(F)$ have the property that left invertible elements are right invertible, and their invertible elements are denoted $(R_v)^*,(\Rf)^*$.

\<{lem}  \lbl{dec1.1}    Assume $L(\f)$ is a field.   Then the natural  map 
$  D^*(F)  \to  (\Rf)^* \bs D^*(\Aaf)$ is surjective.   In particular,
$D^*(F) \to R^* \bs D^*(\Aa)$ is surjective.
  \>{lem}
 
\prf    

To prove surjectivity, let $c \in D(\Aaf_F)$.  Multiplying by  a nonzero element of $\f[t]$, we may assume  $c \in \Rf \meet D(\Aaf_F)^*$.
   Let $I = {\Rf} c \meet D(F)$.  This is an ideal in ${\Rf}  \meet D(F)$.  By \lemref{adeldiv} (6) we have $I = {\Rf}(F) d$ for some $d \in D(F)$.  
So  ${\Rf} c \meet D(f) = ({\Rf}\meet D(F)) d$.
 Using  $c \in D(\Aaf_F)^*$ we see that $\det(c) \neq 0$ and 
 $v(\det(c))>0$ for at most finitely many $v$, so $I \neq (0)$, indeed $I \meet \f[t] \neq (0)$; thus $d \neq 0$.  We have
  ${\Rf} c d \inv \meet D(F)= {\Rf}  \meet D(F)$.  
 Thus $c d \inv \in {\Rf} $, and  $1 \in {\Rf} cd \inv$ so $cd \inv \in ({\Rf})^*$.
\eprf

We aim to show that the class in a localized Grothendieck ring of the 
set of rational points on a given integral conjugacy class, does not depend
on the form of the division ring.  
At this point we   prove a special case:  if one of these sets is nonempty, so is the other.

When $c$ is a regular semi-simple element, let $T_c$ denote the centralizer of $c$.  
 
\<{lem} \lbl{local-conj}    
 Let $D=D_{g,t}, {\,\dot{D}}= D_{\dot{g},t}$, where $\dot{g}$ is another generator of $Aut(L)$.  Fix a place $v$ and let $c \in  D(F_v)$ be a regular semisimple element. 
    
There exists $\dot{c} \in \dot{D}(F_v)$ such that $(D,R_v,c)$ and
  $(\dot{D},\dot{R}_v,\dot{c})$ are $\k(t)_v$-isomorphic.    (We then say that $c,\dot{c}$
  match integrally at $v$.)
 
 Moreover, there exists  a definable
  isomorphism $i:T_{c} \to T_{\dot{c}}$,   mapping
  $T_{c} \meet R_v $ to $T_{\dot{c}} \meet \dot{R}_v$. \>{lem}
  
\prf   
Note:

\noindent (\#) the class of $D_{\dot{g},\bt}$ is a multiple of the class of $D_{g,\bt}$ in the Brauer group of $\f$.

If $D_{g,\bt} \cong M_n$ then by  \lemref{adeldiv} (1,3) we have $(D,R_v) \cong (M_n,M_n(\Oo))$.  Moreover
by (\#), $D_{\dot{g},\bt} \cong M_n$  and the statement is clear.   

Assume therefore that $D_{g,\bt}$ is a division ring.  Then so is $D_{\dot{g},\bt}$.  Moreover
the same conjugacy classes are represented in  $D_{g,\bt},   D_{\dot{g},\bt}$.  The reason
is that a field extension $\f'$ of $\f$  embeds in $D_{g,\bt}$ iff $[\f':f] = n$ and
$D_{g,\bt}$ splits over $\f'$ (\cite{jacobson}, Theorem 4.8, p. 221); it is the same for $D_{\dot{g},\bt}$; but by (\#) it is clear that $D_{g,\bt}$ splits over $\f'$ iff $D_{\g,\bt}$ does.  

If $D(F_v)$ is a division ring, then by the same argument so is $\dot{D}(F_v)$, and they
represent the same classes.  In this case there exists $\dot{c} \in \dot{D}(F_v)$
with $(D,c), (\dot{D},\dot{c})$ isomorphic over $\k(t)$.  In view of the definition of
$R_v$ this suffices in case  $v(t) \neq 0$.   Assume now that $v(t)=0$.  

  By \lemref{adeldiv} (5) we may assume $c \in R_v$ and $c$ has regular semisimple residue
  $\bar{c}$.  By the above, there exists $'\bar{c} \in D_{\dot{g},\bt}$ conjugate to $\bar{c}$.  Lift $'\bar{c}$ to some $'c \in \dot{R}_v(F_v)$, with $'c$ lying in
  some unramified field extension of 
  $F_v$; this extension must be isomorphic to $F_v(c)$, and
  so there exists an element $\dot{c}$ in this extension, with residue $'\bar{c}$, and such that $\dot{c},c$ satisfy the
  same minimal polynomial over $F_v$.   Over $\k(t)_v$ we have $R_v \cong M_n(\Oo) \cong \dot{R}_v$,
  and moving the question to $M_n(\Oo)$, the two conjugate elements $c,\dot{c}$ with   regular semi-simple residues are  clearly $GL_n(\Oo)$-conjugate. (A $K$-basis of eigenvectors $v_i$ for $c$, with   $v_i \in \Oo^n \m \Mm \Oo^n$, is an $\Oo$-basis
  for $\Oo^n$ by Nakayama.)

   The ``moreover" follows from the main statement:   
 any two $\k(t)_v$-isomorphisms $(D,c) \to (\dot{D},\dot{c})$ differ by a conjugation of $D$ by an element of $T_{c}$; such a conjugation induces the identity on $T_{c}$; hence all
  such isomorphisms induce the same isomorphism $i:T_{c} \to T_{\dot{c}}$, and $i$ is definable.   Since 
  some $\k(t)_v$-isomorphism $(D,c) \to (\dot{D},\dot{c})$ respects the integral structure this must be true for $i$.                      
                                      \eprf

\<{remark}  \rm The analogue of \lemref{local-conj} is not true for $S_0^*$ in place
of $R_0^*$-conjugacy or matching.
Indeed the element $s \in D$ does not match any element of $\dot{D}$ in this sense.  
The result does become true if one uses the group $s^\Zz S_0^*$ (an extension of order $n$), {\em and} replaces $\k((t))$ by $\k((t))^{alg}$ in the definition of matching.  \>{remark}

\<{lem}\lbl{el-conj}   
  Let $D=D_{g,t}, {\,\dot{D}}= D_{\dot{g},t}$, where ${\dot{g}}$ is another generator of $Aut(L)$.
  Let $c \in  D(F)$.  Then  there exists $\dot{c} \in \dot{D}(F)$ such that for any $v$,    $(D,R_v,c)$ and
  $(\dot{D},\dot{R}_v,\dot{c})$ are $\k(t)_v$-isomorphic.   (We will say that $c,\dot{c}$ match adelically.)
  
  Moreover, there exists  a definable
  isomorphism $T_{c} \to T_{\dot{c}}$, preserving adelic structure, i.e. mapping
  $T_{c} \meet R_v $ to $T_{\dot{c}} \meet \dot{R}_v$ for each $v$.  \>{lem}
  
\prf  If $L(\f)$ is not a field, then $D, \dot{D}$ are   definably adelically isomorphic over $F$,
so the statement is clear. Assume $L(\f)$ is a field.

  Choose $\dot{b} \in \dot{D}(F)$ with $(D,c) \cong_{ACF_{\k(t)}} (\dot{D},\dot{b})$.  
By \lemref{local-conj}, for any $v$,   there exists
  $\dot{c}_v \in \dot{D}(F_v)$ 
with $(D,R_v,c) \cong_{ACVF_{k(t)_v}} (\dot{D},\dot{R}_v,\dot{c}_v)$.  
  In particular, whenever $v(t)=0$ and $c \in R_v$, we have $\dot{c} \in \dot{R}_v$, so $\dot{c}  \in \dot{D}(\Aa)$, and $\dot{c}, \dot{b}$ are $\dot{D}(\Aa)$-conjugate.  

By \lemref{dec1.1},  they are $\dot{D}^*(F) \dot{R}^*$-conjugate.   So we can find $b$ which is $\dot{D}^*(F)$ -conjugate
to $\dot{b}$ and for each $v$, $G(R_v)$-conjugate to $\dot{c}$.    

The  ``moreover" is as in the proof of the corresponding statement  in \lemref{local-conj}.   \eprf

 \>{section} 

\<{section}{Rational points in integral conjugacy classes}  \lbl{[O]}
 \def\ba{\mathbf{a}}

  In this section we define   integral conjugacy classes $O$, and   compute the class in  the appropriate
  Grothendieck ring of $O(K)=O \meet D(K)$.   The formula obtained will be independent from
  the form of the division ring.  We begin with a discussion of local conjugacy classes.
Let $\f,\k,F=\f(t),K=\k(t),L,L(\f),D=L[s]$  be as in the previous section.

 We will use a Grothendieck ring associated with $ACF_\f$.    Recall $\e_L$ from \secref{lf}.  Let   $\K_0$ be any
$\K_+(ACF_\f)$-semi-algebra, such that $\e_L =0$; this can be achieved by passing to the quotient $\K^*$ of $\K$ associated
to the theory $T^*$ whose models are fields $\f' \geq \f$ such that $\f'(t)$ has no zero
divisors (cf. \secref{T*}).   Let $Gr$ be the multiplicative monoid
generated by  the class $[\CT : T(K)]$ for any torus $T$, and also by
all classes of group varieties; this includes especially $n$ and     $[\bL^*]$.     
Let $\K = \K_0[Gr \inv]$ (cf. \secref{localize}).

If $\f'$ is an extension field of $\f$, we will write 
$\K_{\f'} = \K(ACF_{\f'}) \tensor_{\K(ACF_\f)} \K$.   We also write $\K_c$ for $\K_{\f(c)}$.   All classes $[X]$  in this section refer to $\K$ or some $\K_c$, as specified.

 \<{defn} \lbl{invariant}   Let $\phi$ be a local test function on $D_v$.
   $\phi$ is {\em   $R_v^*$-invariant}
if for any $\f' \models T^*$, for $c \in D(\f'(t)_v)$ and $d \in R_v^*(\f'(t)_v)$
  we have $\phi(c)=\phi(^d c) \in \K_{\f'}$
\>{defn}

 Note that
any $\phi$ whose domain is formally empty counts as invariant.

It would also be possible to  define $R_v^*$-invariance of 
local test functions, in terms of  the action of $R_v^*$ on the subquotients $K(N,M)^n=t^{-N} k[[t]]^n / t^M k[[t]]^n$; see the proof of \lemref{invariance} (c).
 
\<{lem} \lbl{NS}  Let $E$ be a finite dimensional algebra (with $1$) over  a perfect field $F$.
    If $c,c' \in E(F)$ and $c,c'$ are $E^*(F^{alg})$-conjugate, then $c,c'$ are
$E^*(F)$-conjugate.  The same is true for the loop ring $E[[x]]$, or in general for any pro-finite dimensional algebra.  \>{lem}

\prf      Let $X= \{a \in E^*: ac a \inv =c' \}$, viewed as a pro-algebraic variety.
We have to show that $X(F) \neq \emptyset$.  If $E=\liminv E_n$
then $X = \liminv X_n$ with $X_n$ defined in the same way; so the question
reduces to the finite-dimensional case.  
 
Note that $X$ is coset of  
$T_c = \{a \in E^*: aca \inv = c \}$.  Now $T_c$ is the group of units
of the ring $\bar{T}_c = \{a \in E: ac = ca \}$. By \lemref{90}, $X$ has a rational point. \eprf

\<{defn} \lbl{match}   Let   $c \in D(\f'(t)_v)$, $\dot{c} \in \dot{D}(\f'(t)_v)$.  We say that
$c,\dot{c}$ are in {\em matching conjugacy classes} if  $(D,R_v,c) , (\dot{D},\dot{R_v},\dot{c})$ are isomorphic over  $\k(t)_v$.

 Let $\phi,\dot{\phi}$ be $R_v^*$-invariant     functions on $D,\dot{D}$
respectively.  We say that $\phi,\dot{\phi}$ {\em match}
if for any $\f' \models T^*$, $\phi(c)=\phi(\dot{c})  \in \K_{\f'}$ whenever   
$c \in D(\f'(t)_v)$, $\dot{c} \in \dot{D}(\f'(t)_v)$ are in matching conjugacy classes.
  \>{defn}  
  
This definition works thanks to \lemref{match0}, \lemref{invariance} (b) below.  

\<{lem} \lbl{match0}  For $c,c' \in D(\f'), \f' \models T^*$, we have $(D,R_v,c) \cong (D,R_v,c')$   over $\k(t)$
iff $c,c'$ are $R_v^*(\f')$-conjugate. \>{lem}

\prf     For $v(t) \neq 0$ there is nothing to show, since $R_v = D$.  Assume $v(t)=0$.  

Note  that over $\k(t)$, $R_v$ is a matrix ring over a valuation ring, \lemref{adeldiv} (1).   The stabilizer of $R_v$ in $D_v^*$ is $Z R_v^*$, where $Z$ is the center of $D^*$. 
Now  if $g \in M_n(\k(t))$ and over $K^{alg}$ we have $g \in Z GL_n(\Oo)$,
 then $g \in Z GL_n(\Oo_K)$, as one easily sees by considering the valuation of
 matrix coefficients.  Hence $Z R_v^* (K^{alg}) \meet M_n(K) = Z(K) R_v^*(K)$.  One can also check directly that the $GL_n(K)$-stabilizer of $M_n(\Oo_K)$ is
 $Z M_n(\Oo_K)$.
 
So $c,c'$ are $Z R_v^* ( \k(t))$-conjugate, and hence also $R_v^* (\k(t))$-conjugate.
$R_v^*(\k(t))$ may be viewed as the projective limit of $R_v^* (\k(t))/(1+ \Mm_\alpha)$ over the various congruence ideals $\Mm_\alpha$ of $R_v$.
By \lemref{NS}, $c,c'$   are $R_v^*(\f'(t))$-conjugate.  \eprf

In order to match anything,  $\phi$  must be strongly $R_v^*$-invariant; but by 
\lemref{invariance} (b) this is automatic.

\<{lem} \lbl{invariance}  (a)  Let $a,a' \in D(\f(t)_v) $ and assume $a',a'$ are $R_v^*$-conjugate
over $\k(t)_v$.  Then they are conjugate by an element of $R_v^*((\f(t)_v)$.  

(b)    Let $\phi:D(K_v) \to \K$ be $R_v^*$-invariant  (\defref{invariant}).  
Then $\phi$ is strongly $R_v^*$-invariant (\defref{invariant-gen}).  In other
words, if   $\f' \models T^*$, for $c \in D(\f'(t)_v)$ and $c' \in R_v^*(\f'(t)_v)$
  we have $\phi(c)=\phi(c') \in \K_{\f'}$ provided that $c,c'$ are $R_v^*$-conjugate
  over $\k(t)_v$.
  
(c)  Let $\phi,\dot{\phi}$ be $R_v^*$- (resp. $\dot{R}_v^*$-) invariant local test functions on $D_v,\dot{D}_v$ respectively.
To test them for matching, it suffices to show that $\phi(c) = \dot{\phi}(\dot{c}) \in \K_{\f'}$ 
for matching $c,c' \in D(F')$, $\f' \models T^*$, $F'  = \f'(t)$.  Also, if $W$ is a proper subvariety of $D$, it suffices to consider $c,c' \in D(F') \m W$.
\>{lem}

\prf (a)   If $v(t) \neq 0$, this is immediate from \lemref{NS}.  Otherwise, multiplying by a scalar element we may assume $a,a' \in R_v(\f(t)_v)$.   Then the 
statement follows from the profinite part of  \lemref{NS}. 

(b)  Immediate from (a) and the definitions.

(c)  We prove the case $v(t) \neq 0$; the case $v(t)=0$ is similar but easier, and will
not be needed.   Take $v(t)>0$ and $\f'=\f$ to simplify notation; write $F=\f(t)$.   

Consider the ring  defined by the formula \ref{S} at $v(t)>0$; we will refer to it as
$S_0$. 
Multiplying an element of $D$ by a suficiently large power of $t$ will put it in $S_0$, so it suffices to consider $D_0^*$-conjugacy on elements of $S_0$. 

  $S_0^*$-conjugacy is a finer relation than $D_0^*$-conjugacy, so
$\phi$ is $S_0^*$-invariant, and similarly $\dot{\phi}$.  As in (b), $\phi$ is strongly
$S_0^*$-invariant.  On the other hand $\phi$ is invariant under additive $t^m S_0$ translations for large enough $m$, and in fact descends to a function $\bar{\phi}$ on the quotient $S[m] = S_0 / t^m S_0$.  It follows that $\bar{\phi}$ is $S_0 / t^m S_0$-translation invariant; 
using \lemref{NS} for the ring $ S_0/ t^m S_0$ as in (a,b), we see that 
$\bar{\phi}$ is  strongly $S_0^*$-invariant.  By \lemref{invariance-descent}, $\bar{\phi}$  descends to a function $\psi$ on the set $Y$ of conjugacy classes of
 $S_0/t^mS_0$.  (We use here that $S[m]^*$ and its quotient subgroups have invertible classes  in $\K=\K[Gr \inv]$.)   Let $\dot{Y}$ be the corresponding construction for $\dot{D}$.
 
Let $E$ be the equivalence relation of $D^*$-conjugacy, restricted to $S_0(F_v)$.  Then $E$ respects $t^mS_0$-translation, so descends to
$S_0(F_v) / t^m S_0(F_v)$.  Being coarser than  $S_0^*$-conjugacy,   it induces an equivalence relation $E$ on  the image $Y_{F_v}$ of $S_0(F_v)$ in
$Y$.     Matching induces a bijection 
$Y_{F_v}/E \to \dot{Y}_{F_v} / \dot{E}$.  Now $\psi$ and the corresponding function $\dot{\psi}$
are $E$-invariant (resp. $\dot{E}$-invariant), and the question is
whether the functions they induce on $Y_{F_v} / E , \dot{Y}_{F_v} / \dot{E}$ correspond.
But $S_0(F_v) / t^m S_0(F_v) = S_0(F) / t^m S_0(F)$, so
every element of $Y_{F_v}, \dot{Y}_{F_v}$ can be represented by an $F$-rational point
(outside a given subvariety $W$); hence it suffices to check invariance at these points.     
  \eprf
  
\ssec{Global conjugacy classes}

 By the {\em integral conjugacy class} of a regular semi-simple element $c \in D(K)$
 we mean the set $O$ of all elements $d \in D(K)$ such that  for any place $v$, $(D_v,R_v,c)$ and
 $(D_v,R_v,d)$ are isomorphic over $K$.

Note that for almost all $v$, we have $c \in R_v$, and so  $d \in R_v$, both  with regular semisimple reduction.

  Let  $\bL = L(\k)$.  Let $N_{L/\k}$ be the norm map on $L$.  
    Let $L_1$ be the kernel of the norm map $\bL \to \k$:
$$L_1 = \{a \in \bL^*: \Pi_{i=1}^n g^i(a)=1 \}$$
By Hilbert 90  $x \mapsto g(x)/x$ is a surjective map $\bL^* \to L_1$; 
there   exists a constructible section   of this function; in particular  $[L^*] =  [L_1] [G_m]$.    In this expression, $L^*$ and $L_1$ are 
viewed as $ACF_\f$-definable sets; to avoid confusion, we will use the notation
 $\bL=L(\k)$ when interested in $\k$-points.   Similarly $\bL_1 = L_1(\k)$.

Let $N=  R^* \meet D^*(K)$, and let $\bN$ be the image of $N$ in $PD=D^*/Z$.
  By \lemref{adeldiv} (2), we have $N = R^* \meet D^*(K) = L(\k)[s,s\inv]^* = L(\k)^*s^\Zz$. 
  We view $\bN=N /   t^\Zz$ as a subgroup
of $PD$, containing $\bL:=L(\k)^* / \k^*$ as a subgroup of   index $n$.   Note
that we have   definable representatives $1,s,\ldots,s^{n-1}$ for the cosets of $L(\k)^*/\k^*$ in $\bN$.  Thus in  the Grothendieck ring $\K(ACF_\f)$ we have
$[\bL] = n [L^*/G_m] = n [L_1]$.

  While $N$ is not a limited set, it is easy to describe a limited set $N'$ of elements
  of $N$, such that $N,N'$ have the same image in $PD$.     So $\bN$ can be seen as a constructible set over $\f$.

Let $CN = \union_{b \in N} T_b$, where $T_b$ is the centralizer of $b$ in $D^*$.
The definition of $CN$ is geometric, i.e. if $c \in D(F')$ and $c$ commutes with  some $n \in N(K^{alg})$, then $c$ commutes with some $n' \in N(F')$.   Indeed if $n = as^i$
with $a \in \bL$, then $c$ commutes with any element $xs^i$ where $a'$ lies in
the linear space  $\{x \in L: cxs^i = xs^i c \}$, intersected with the  Zariski open set $\det(x) \neq 0$.

\<{lem}\lbl{o-o}  Let $O$ be an integral conjugacy class.  If $O \meet CN \neq \emptyset$
then $O \subset CN$.  \>{lem}

\prf  Let $c \in O \meet CN$, and let $c' \in O$, with say $c,c' \in D(F')$.  Then $c' = r c r \inv$ for some 
$r \in D^*(F')$, and for each $v$ we have $r T_c = r_v T_c$ for some $r_v \in R_v$.
Since $c \in CN$, there exists $n \in N$ with $n \in T_c$.  So $r n r \inv  = r_v n r_v \inv$.
It follows that $rn r \inv \in R_v$ for each $v$ (since $N \subseteq R_v$.)  So
$r n r \inv \in \meet_v R_v = N$.  Thus $r c r \inv $ commutes with an element of $N$,
so $c' \in CN$. \eprf

  Let $D=D_{g,t}, \dot{D}=D_{\dot{g},t}$; define $\dot{N} = \bL \dot{s}^\Zz$, and
  $\dot{CN} = \union _{\dot{b} \in \dot{N} \m {\bL}}  T_{\dot{b}}(K)$. 

\<{lem} \lbl{orb} There exists a  definable bijection $CN \to \dot{CN}$, mapping
integral conjugacy classes to matching integral conjugacy classes.  \>{lem}

\prf    Over $\lf$ we have a bijection $h: D \to \dot{D}$ preserving adelic structure, i.e. preserving $R_v$ for each $v$.  In particular as noted above we have $h(N)=\dot{N}$,
and $h$ maps integral conjugacy classes to matching ones.  
 Moreover $h({\bL}) = {\bL}$
since e.g. ${\bL}$ is the subspace generated by $N \meet (1+N)$.  
$h | \bL$ must be a Galois automorphism; by pre-composing with conjugation by a power of $s$ we may assume that $h|\bL=Id_{\bL}$. 
 It remains to consider $CN \m \bL$.

 If $a \in  CN \m {\bL}$ then 
there exists some element $bs^i \in T_a \meet {\bL}s^i  $, $i \neq 0$.
We may take $i=1$, since if $b_1 \in {\bL}$ is any element with $g^i(b_1)/b_1 = g(b)/b$
then $b_1s \in T_{bs^i}$ (and using Hilbert 90).   Hence 
$ CN \m {\bL} = \union_{b \in {\bL}^*}( T_{bs}(K) )$, and similarly on the $\dot{D}$ side. 

Given $b \in {\bL}^*$, the elements $bs \in D^*, b\dot{s} \in \dot{D}^*$ have the
characteristic polynomial $X^n-N_{{\bL}/\k} (b) t $.  If $b's$ has the same characteristic
polynomial as $bs$, then $b'=cb$  for some $c \in L_1$.  In this case by Hilbert 90
we have  $a \inv g(a) = c$ for some $c \in \bL$; so conjugation by $a$ takes $bs$ to $b's$
while preserving integral conjugacy.  This shows that the integral conjugacy class of $bs$ is precisely $L_1 bs$.  It follows that $h(L_1bs) = L_1 b \dot{s}$.  

In particular $h(L_1s) = L_1 \dot{s}$; we may assume $h(s) = \dot{s}$.   So $h(bs)=h(b) \dot{s} = b \dot{s}$, for $b \in L^*$.

The isomorphism $h| T_{bs}$ takes $bs$ to $b \dot{s}$, and as such, 
is determined up to conjugation by the centralizer of $bs$.  Hence the restriction
  $h| T_{bs}$ is determined uniquely, given $b$; we denote it $\alpha_b: T_{bs} \to T_{b\dot{s}}$.  So $\alpha_b$ is $\f(b)$-definable.  

Note that if $k \in \k^*$ then $T_{kbs} = T_{bs}, T_{b\dot{s}} = T_{kb\dot{s}}$
and $\alpha_{kb} = \alpha_{b}$.  

Note that $\alpha_b$ takes any element of $T_{bs}$ to an element of the matching
integral conjugacy class (since this is true of $h$.)

Also, two centralizers $T_{b_1s},T_{b_2s}$ are disjoint or equal; they are equal
when $b_1s,b_2s$ commute, i.e. when $b_1/b_2 = g(b_1) / g(b_2)$; equivalently when
  $b_1/b_2 \in \k^*$.
  
Thus we may define $\alpha:   \union_{b \in \bL^*}  T_{bs} \to \union_{b \in \bL^*} T_{b\dot{s}}$
by $\alpha(x) = \alpha_b (x)$ when $x \in T_{bs}$.  It is clearly a bijection,
and for any two matching integral conjugacy classes $O, \dot{O}$, it restricts to a bijection   $\union_{b \in \bL^*}(O \meet T_{bs}(K) ) \to
\union_{b \in \bL^*}(\dot{O} \meet T_{b\dot{s}}(K) )$.

\eprf

Let $O$ be an integral   conjugacy class  of regular semi-simple elements.  $O$ is a  $\T_f=ACVF_{C;\f}$-Ind-definable subset of $D(K)$ (\S \ref{valued fields over a curve}.)  
  By \lemref{vfc1}, $O$ is also Ind-definable for $ACF_\f$.
$O$ is $T^*$-limited, since   $R_v$ is   bounded   for any $v$ with $v(t) = 0$, and $R_v / Z$ is $T^*$-limited  for $v(t) \geq 0$,   where
$Z$ is the center of $D^*$.


 We define the canonical torus with adelic structure  $T$ associated to
$O(K)$ as follows.  For any $c \in O(K)$ let $T_c= \{a \in D^*: ac=ca \}$.
If $c,c' \in O \meet D(K)$ there exists $d \in D(K^{alg})$  with 
$dc d \inv = c'$; such a $d$ gives an isomorphism $ad_d: T_c \to T_{c'}$, preserving
adelic structure,
$ad_d(x) = d x d \inv$; but $ad_d$ does not depend on the choice of $d$, so we can
write $f_{c,c'} =ad_d$.  Clearly $f_{c',c''} f_{c,c'} = f_{c,c''}$.  We can factor our this system to obtain a torus $T$, with adelic structure; given any $c \in O$ we have an 
isomorphism $f_c:  T \to T_c$; and $T$ is definable over the field of definition of $O$.

Let $\CT=   (Z(\Aa)  \TR ) \bs T(\Aa)$.   This is an Ind-definable group. 
The diagonal embedding $T(K) \to T(\Aa)$ induces a homomorphism
$T(K) \to \CT$.   The image of $T(K)$ in $\CT$ is Ind-definable, and hence
so is the corresponding coset equivalence relation $E$ on $\CT$.

 \<{prop} \lbl{E2}    Assume $O(K) \meet CN = \emptyset$, and let $c \in O(K)$.   
  Then  $\CT/ T(K)$ is 
 $(T^*,n L_1)$-representable, and we have
$$[O(K)]  =     n [L_1]   [\CT: T(K)] $$ in $\K_c$. 
  \>{prop} 

\prf   Fix $c \in O(K)$.  We identify $T$ with $T_c$.  We will verify the conditions of \lemref{srep3}, with
$V,X,Z$ of that lemma corresponding here to 
 $\CT,O(K),\bN$, respectively, and $E$ being the $T(K)$-coset equivalence
relation.

We identify $\bL^* / k^*$ with $L_1$ , and $\bN =N / t^{\Zz}$ with the semi-direct product of $L_1$
with $s^\Zz / s^{n\Zz}$.    $N$ acts on $O(K)$ by conjugation inducing an action
of $\bN$; the latter action is free  because of the assumption:  $O \meet CN  = \emptyset$.

Define $f: O(K) \to \CT/T(K)$ as follows.  Let $d \in O(K)$.  Then $d=a c a \inv = b c b \inv$ for some $a \in D^*(K), b \in R^*$.  We have $b \inv a \in T(\Aa)$.  If also
$d = {a'} c {a'} \inv = {b'} c {b'} \inv$ with ${a'} \in D^*(K), {b'} \in R^*$, then $a \inv {a'} \in T_c(K)$
and $b \inv {b'} \in (T_c \meet R) = \TR$, so $b \inv a, {b'} \inv {a'}$ have the same class 
in $\TR \bs T(\Aa) / T(K)$ .
and in particular in $(\TR Z) \bs T(\Aa) / T(K)$.  
Let $f(d)$ denote this class.  The graph of $f$ pulls back
to an Ind- definable subset of $O(K) \times \CT$. So $f$ is Ind-definable.  

If $f(d) = f({d'})$, then for some $a,{a'} \in D^*(K), b,{b'} \in R^*$ we have 
$b \inv a = {b'} \inv {a'}$ away from $0$,  and $d=a c a \inv = b c b \inv$, 
${d'} = {a'} c {a'} \inv = {b'} c {b'} \inv$.  So ${b'} b \inv  = {a'} a \inv \in D^*(K) \meet Z(\Aa) R^* $.  By \lemref{adeldiv} (2), we have ${a'} a \inv \in   = Z N$. 
it follows that $Z Na= Z Na'$ so   $d,{d'}$ are $\bN$-conjugate.
So the fibers of $f$ are   $\bN$-orbits.   

 We now show:   for any $\f' \geq \f$  
and any $\f'$-definable $\bN$-orbit $U$,
 $\union_{d \in O(K)(\f')} f(d)  \subseteq  \union_{d \in \CT(\f')} d T(K)$.   

Let $d \in O(K)(\f')$.  We have to show that $f(d)$ contains
a point of $\CT(\f')$.  Indeed by   \lemref{NS} there exists $a \in D^*(\f')$ with $a c a \inv = d$.  By definition of $O(K)$ there exists $b \in R^*$ with $b c b \inv = d$;
$b \TR$ is uniquely determined.  Thus $\TR b \inv a  \in \TR \bs  T(\Aa)  = \CT$
is determined, and  $Z_0 \TR b \inv a \in \CT(\f')$.  

Now assume $\f' \models T^*$; so $L(\f')$ is a field.
Let $Z_0 \TR e \in \CT(\f')$, $e \in T(\Aa)$.   We have $\TR e \in \CT= \TR \bs T(\Aa) \leq R^* \bs D^*(\Aa)$.   
 By \lemref{dec1.1} there exists   
$a \in D^*(\f'(t))$ with  $R^* a = R^* e$; so $b \inv a = e$ for some $b \in R^*$.
Since $e \in T(\Aa)$ we have $ac a \inv = bc b \inv$.  
  Let $d = a c a \inv$.  Then $d \in O(K)$ and  $f(d) = \TR e T(K)$.
\eprf

\<{rem}  \rm It follows from the proposition, in particular, that there exists
a definable $W$ with $\CT \m W$ formally empty,  and such that 
$E$ is definable on $W$.  It is not the case that $E$ is definable-in-definable-families.
It is possible, but unnecessary for our purposes, to modify $\TR$ by using a definable
subgroup $H_0$ at $0$ such that $H_0$ is bounded modulo the center, but
$H_0(\f'((t)))$ contains $D(\f'((t)))$ for any $\f' \models T^*$.  One can
take $H_0 =  = s^\Zz Z S_0^*$.  
This has no effect on classes in $\K$ since   $\e_L=0 \in \K$, but yields an equivalence relation that
is  definable-in-definable-families.  \>{rem}

    If $D,\dot{D}$ are two forms of $M_n$ with the same adelic structure ( \S \ref{norm-s})  it is possible to match their integral conjugacy classes; see \lemref{match}
    for the matching of $R_v^*$-classes  of $D$ to those of $\dot{D}$, for any
    $v$.   We say that $c,\dot{c}$ match if their $R_v^*$-conjugacy classes match
    for any $v$.  
 Write $\dot{O}$ for the $\dot{R}^*$- class corresponding to $O$.  
   
 \<{cor} \lbl{equal}   Assume $O(K)$ has an $F$-rational point.  Then $[O(K)]=[\dot{O}(\K)]$.\>{cor}
 \prf  By \lemref{el-conj}, $\dot{O}(K)$ also has an $F$-rational point.    So \propref{E2} is valid
 in $\K$, and by either this or  \lemref{orb} we obtain $[O(K)]=[\dot{O}(\K)]$. \eprf

\<{remark} 
 \lbl{H1} \rm  Let $\tJ = \CT^0 / T(K)$; it is an algebraic group over $\f$.
The relation between $[\CT^0 : T(K)]$ and $[\tJ]$ is very close. The interesting
case is that the torus $T$ does not split, and we have $T=R_{C'/C} G_m$
for a certain curve $C'$ over $C=\Pp^1$; i.e. $T$ is obtained by restriction of
scalars from $\f(C')$ to $\f(C)$.   The adelic constructions commute with restriction
of scalars, and $\tJ$ can be identified with a certain Rosenlicht generalized Jacobian
of $C'$.  For any field $\f'$ such that $C'(\f') \neq \emptyset$, it can be shown
that there exists a rational section $\tJ \to \CT^0$, and therefore 
$[\CT^0: T(K)] = [\tJ]$, so we simply have the class of an algebraic group.
In general while $\CT^0, \tJ$ and the exact sequence $T(K) \to \CT^0 \to \tJ$
are all defined over $\f$, no section exists, so the quotient group cannot quite
be identified with the algebraic group; cf. \cite{lichtenbaum}.  
\>{remark}

\<{rem} \rm  on a group theoretic level, we are using a special case of the bijection
$$ (R \meet K) \bs (RS \meet KT) / (S \meet T) \cong (R \meet S) \bs (RK \meet ST) / (K \meet T)$$
valid for any subgroups $R,S,K,T$ of a group $G$.
The bijection maps the $(R \meet K),(S \meet T)$-double coset of $rs=kt$ to the 
$(R \meet S), (K \meet T)$-double coset of $r \inv k = s t \inv$.
  In our case we take $S=T$ and
have $RK=G$, yielding $T \meet R \bs T / T \meet K$ on the right. \>{rem}

\<{remark} \lbl{Jtorsor} Let $O$ be an $R^*$-conjugacy class, defined over $\f$.  Given an element $c \in O$,
we defined above a torus $T_c$, a Jacobian $J_c$, and a map $z_c: O \to J_c$; the construction depended on $c$.  However the
triple $(T_c,J_c,z_c)$ descends to $(T,U,z)$ where $U$ is an $\f$-definable torsor
of an Abelian variety $J$ over $\f$, and $z: O \to U$.  For any $c \in O$ we have a $c$-definable isomorphism $(T_c,J_c,z_c) \to (T,U,z)$.  \>{remark}

\prf  
Let $O$ be an $R^*$-conjugacy class, defined over $\f$.  We define the canonical torus with adelic structure  $T$ associated to
$O \meet D(K)$ as follows.  For any $c \in O \meet D(K)$ let $T_c= \{a \in D^*: ac=ca \}$.
If $c,c' \in O \meet D(K)$ there exists $d \in D(K^{alg})$  with 
$dc d \inv = c'$; such a $d$ gives an isomorphism $ad_d: T_c \to T_{c'}$, preserving
adelic structure,
$ad_d(x) = d x d \inv$; but $ad_d$ does not depend on the choice of $d$, so we can
write $f_{c,c'} =ad_d$.  Clearly $f_{c',c''} f_{c,c'} = f_{c,c''}$.  We can factor our this system to obtain a torus $T$, with adelic structure; given any $c \in O$ we have an 
isomorphism $f_c: T \to T_c$; and $T$ is definable over any field of definition for $O$.
This induces an isomorphism $f_c: J \to J_c$, where $J = \TR \bs T(\Aa) / T(K)$.
 
For any $c,c' \in O \meet D(K)$ we obtain an element $j(c,c')$ 
of $J$:  pick $a \in R,  b \in D(K)^*$ with $ac=bc=c'$; let $t=a \inv b$; then $t \in T_c(\Aa)$; let $j(c,c') =f_c \inv (t)$.

If $c_1,c_2,c_3 \in O \meet D(K)$, let $a_i \in R$, $b_i \in K$, $a_1c_1=b_1c_1=c_2$,
$a_2c_2 = b_2 c_2 = c_3$.  Let $a_3 = a_2a_1, b_3 = b_2 b_1$; so
$a_3c_1=b_3 c_1 = c_3$.    Let $t_i = a_i \inv b_i$.  Then $j(c_1,c_3)=f_{c_1}(t_3)=
f_{c_1}( a_1 \inv a_2 \inv b_2 b_1) = f_{c_1}(  (a_1 \inv a_2 \inv b_2 a_1 ) (a_1 \inv b_1))$  Now $f_{c_1}(  (a_1 \inv a_2 \inv b_2 a_1 )= f_{c_2}(a_2 \inv b_2) = f_{c_2}(t_2)=j(c_2,c_3)$; and $f_{c_1}(a_1 \inv b_1) = f_{c_1}(t_1) = j(c_1,c_2)$.  Thus
$j(c_1,c_3)=j(c_2,c_3)  j(c_1,c_2)$.  It follows that there exists a $J$-torsor $U$,
defined over $\f$, and a map $z: O \to U$, with $z(c_1,c_2) + z(c_1)=z(c_2)$.  
\eprf

   \<{question}  Does the equation of \propref{E2}, descend to $\K$ ? 
\>{question}

\ssec{$\delta_K$ is geometric}
 
Consider  global functions $\phi$ given by a uniformly definable family $(\phi_v)$.
We assume  that $\phi_v$ has bounded support
for all $v$, contained in $R_v$ for almost all $v$; but not necessarily that $\phi_v$ is locally
constant.   Assume $\phi_v$ is $R_v^*$-invariant.  Recall that $\K=\K[Gr \inv]$.

\<{prop} \lbl{geometric}   Let  $\phi,\dot{\phi}$ be matching definable global functions  as above.  
Then $\CR(\phi) = \CR(\dot{\phi}) \in \K$.
  \>{prop}  
  
 \prf  The support of $\phi$ is a limited subset $X$ of $D(K)$; the equivalence
 relation $E$ of integral conjugacy is definable on $X$. Similarly, let $\dot{X}$ be the support of $\dot{\phi}$. and let $\dot{E}$ be integral conjugacy.  
 Since $\phi,\dot{\phi}$ match,
 we can identify the quotients $X/E, \dot{X}/E$; so we have quotient maps
 $\pi: X \to Y, \dot{\pi}: \dot{X} \to Y$.     By fibering over $Y$ (\ref{presheaf}) we can
 reduce to the case that $Y$ is a point, i.e. $X,\dot{X}$ form a single integral conjugacy class.  If this class is central, the statement is clear.  If it is not regular semi-simple,
 then $[X] = [\dot{X}] = 0 \in \K$ since $X(\f') = \emptyset$ whenever $L(\f')$
 is a field.  So we assume $X,\dot{X}$ are integral conjugacy classes of regular
 semi-simple elements.  
 
Since $\K[\mathcal{N} \inv] = \K[\mathcal{N} \inv][ \mathcal{N} \inv]$,
replacing $\K$ by $\K[Gr \inv]$ we may assume $\K=\K[Gr \inv]$.

 For any  $c \in X'$, let $\f'=\f(c)$.  
 
 \claim{}   $\CR(\phi) = \CR(\dot{\phi}) \in \K_{\f'}$. 
 
  Indeed
 there exists $\dot{c} \in \dot{X}(\f')$ (\lemref{local-conj}).  By strong invariance of $\phi,\dot{\phi}$
 we have $\CR(\phi) = \phi(c) [X]$, and $\CR(\dot{\phi}) =\dot{\phi}(\dot{c})[\dot{X}]$,
 and by strong matching we have $ \phi(c)= \dot{\phi}(\dot{c}) \in \K_{\f'}$.
By \lemref{equal}, $[X]=[\dot{X}] \in \K(\f')$.  So $\CR(\phi) = \CR(\dot{\phi}) \in \K_{\f'}$.

Let $U=\CR(\phi), \dot{U} = \CR(\dot{\phi}) $.  Since $[U]=[\dot{U}] \in \K(\f(c))$, 
 for any $c \in X$, summing over $c \in X$ we obtain:
 $$[\dot{U}]  [X] = [U][X]$$
Similarly, since  $[X]=[\dot{X}] \in \K(\f(c))$ for any $c \in X$, we have:
 $$ [X]^2 = [\dot{X}] [X]$$
 and by symmetry, $[\dot{X}]^2= [X] [\dot{X}]$.

Let $[A]= [L_1] [ \CT: T(K)]$.  Then $[A] \in Gr$.  
 
 If $X \meet CN = \emptyset$, then by \propref{E2} we have the relations of \ref{211}.  Thus $[X]=[\dot{X}] \in \K$.

Otherwise, by   \lemref{o-o},    $X \subseteq  CN$.  So $[X]=[\dot{X}] \in \K$. by \lemref{orb}.    \eprf 

\<{rem}  The proof of \lemref{geometric} does not require subtraction, and goes through
for the Grothendieck semiring.  \>{rem}

 \>{section}  

 \<{section}{An expression for the Fourier transform, and proof of \thmref{A}}
 
We can now express the Fourier transform at $\f((t))$ in terms of 
the Fourier transform at  $\f((t-1))$ (where $D$ splits) and  $\CR$.  
This will lead to a proof of \thmref{A}.
  \lbl{express}

For a finite set of places $w$, let  $R_w = \Pi_{v \in w} R_v$.    If $\phi=(\phi_v)_{v \in S}$ is
a family of local test functions on a set $S$ of places, and $\phi' =(\phi'_v)_{v \in S'}$
is a family of local test functions on a disjoint set $S'$ of places, we write
$\phi \cct \phi'$ for their conjunction on $S \union S'$.   

We will write $K_0,K_1$ for $K_{v_0},K_{v_1}$; here $v_1$ is the valuation with $v_1(t-1)>0$.

Let $O$ be an integral conjugacy class.   Since $O$ is $T^*$-limited, there exist   bounded definable sets $B_0,B_\infty$
such that if $\f' \models T^*$ then $O_v \subseteq B_v$ for $v(t) \neq 0$.  
Recall that $[O]$ is defined to be the class of $O' : = B_0 \meet O \meet B_\infty$.
More generally,   let
 $\CR_O(\phi) = \sum \{\phi(a): a \in O'(K) \} \in \K $; this clearly does not depend on the
 choice of $B_0,B_\infty$.     Also let $\CR_c=\CR_{O(c)}$, where
 $O(c)$ is the integral conjugacy class of $c$.    Since $O(c)$ is already a limited
 set, $\CR_c(\phi)$ is defined even if $\phi$ does not have bounded support.
  
 Let $\phi_{v}^{st}$ be the characteristic function of $R_{v}$.  We call this the standard test function at $v$. 
 Let $\phi_{v;n}(x) = \phi_{v}(e x)$ where $e$ is any element of $\k(t)$ with $v(e)=n$;
 for instance $c=(t-\alpha)^n$ if $v(t-\alpha)>0$.  This function clearly does not depend on the choice of $c$; we call such functions semi-standard.
 
 Let $\phi_1^{c,m}$ be the characteristic function of $Ad_{R_1}(c)+(t-1)^mR_1$. 
 
A collection $(\phi_v)_v$ is  semi-standard if $\phi_v$ is a semi-standard test function for each $v$, and standard almost everywhere.

 \<{lem} \lbl{lastplace}  For $v \neq 1$ let $Y_v \subset D$ be a bounded $ACVF_{\f(t)_v}$-definable set, with $Y_v = R_v$ for all $v$ outside some $ACF_\f$- definable finite subset of $\Pp^1$.  
  Let $c \in D(\f(t))$ with $c \in Y_v$ for $v \neq 1$.  Then there exists a bounded, 
  $R_1^*$-invariant $\f(t)_{v_1}=\f((t-1))$-definable neighborhood $Y_1$ of $c$ such that
if  $y \in Y= D(\k(t)) \meet \meet_v Y_v$    then $y,c$ are $R_1^*$-conjugate.  \>{lem}
  
 \prf  Let $U_m = c+(t-1)^mR_1$; it is  a  bounded, $\f(t)_{v_1}=\f((t-1))$-definable, neighborhood of $c$.  So is $V_m= ad_{R_1^*}( U_m)= \{x: (\exists y \in R_1^*)(yxy \inv \in U_m\}$; moreover $V_m$ is $R_1^*$-invariant.
     Let
 $Y^0 = D(\k(t)) \meet \meet_{v\neq 1} Y_v \meet V_0$.  Then $Y^0$ is
  a limited subset of $D(\k(t))$.  
  
Let ${R_1^*}_y = Ad_{R_1}(y)$ be the $R_1^*$-conjugacy class of $y$.  
Being the image of a bounded set defined by weak inequalities under a continuous definable map, 
$Ad_{R_1}(y)$ is a closed and bounded subset of $D(K_1)$ in the valuation topology.     

\claim{}  For any $y \in Y^0$ there exists $m=m(y)$  such that if $u \in U_m$ and $u,y$
are $R_1^*$-conjugate then $u,c$ are $R_1^*$-conjugate.

\prf  Fix $y \in Y^0$.  If $y \in Ad_{R_1}(c)$ then $m=0$ will do , since if $u,y$ are $R_1^*$-conjugate
then so are $u,c$. If $y \notin Ad_{R_1}(c)$, then since $Ad_{R_1}(y)$ is closed,
and $c \notin Ad_{R_1}(y)$, some neighborhood $U_m$ of $c$ is disjoint from $Ad_{R_1}(y)$.
In this case no $u \in U_m$ is $R_1^*$-conjugate to $y$.  \eprf
  
Now $Y^0$ is a limited (so $ACF_\f$-definable) set, and so by compactness,
for some $m$, for any $y \in Y^0$, if  $u \in U_m$ and $u,y$
are $R_1^*$-conjugate then $u,c$ are $R_1^*$-conjugate. 
  Let $Y_1 = ad_{R_1^*}( U_m)$. 
  This is a bounded, 
  $R_1$-invariant $\f((t-1))$-definable neighborhood of $c$.  Define $Y$ as above.
  If $y \in Y$ then $y \in Y^0$, and $y $ is $R_1^*$-conjugate to some $u \in U_m$.
  But then $u,c$ and hence $y,c$ are also $R_1^*$-conjugate, as required.    
 \eprf

 \<{lem}  \lbl{discrete2} Let $\phi_{0} $  
  be a local test function  on  $D$ over $K_0=K_{v_0}$, 
 and let $c \in D(F)$.     
 
   Let $(\theta_v: v  \neq 0,1\}$ be  a semi-standard collection,   with $\theta_v(c)=1$.   
    Let $\phi_1=\phi_1^{m,c}$,   $\theta'=\fF \inv \theta$, $\phi_1' = \fF \inv \phi_1$
   
      Then
 if $m$ is sufficiently large,  
  we have
   $$\CR_c \fF \phi_{0}   =  \CR (\phi_{0} \, \cct\phi_1' \cct\theta') $$
   
If  $\phi_0$ is  $R_0^*$-invariant, and $c \notin CN$, we have:
$$\Ff \phi_0 (c) =   (n  [L_1]  [\CT_c : T(K)] ) \inv \delta^K(\phi_0 \cct \phi_1' \cct \theta') )$$

 \>{lem} 
 
 \prf By \lemref{lastplace} we have  $\CR (\phi_{0} \, \cct\phi_1\cct\theta) = \CR_c(\phi_{0} \, \cct\phi_1\cct\theta) $.  Now $\phi_1, \theta $ take the value $1$ on $\CR_c$.  So
 $\CR (\phi_{0} \, \cct\phi_1\cct\theta)  = \CR_c (\phi_0)$.    This gives:
 
  $$   \CR_c (\phi_0) = \CR (\phi_{0} \, \cct\phi_1\cct\theta) $$

 Applying this formula  to $\phi_0'  :=   \fF \phi_{0}$, we find:
 $$\CR_c (\phi_{0}' )  = \CR (\phi_{0}' \, \cct\phi_1\cct\theta) $$
By Poisson summation, 
 $\CR( \phi_{0} \, \cct\phi_1' \cct \theta') =  \CR (\phi_{0}' \, \cct\phi_1\cct\theta) $.  
The first formula follows.  

   By \lemref{invariance}, $\Ff \phi_0$ is strongly invariant.  
So $\Ff \phi_0(y) = \Ff \phi_0 (c) \in \K_y$, and
  $\CR_c(\Ff \phi_0) = \sum_{y \in O(c)}\Ff \phi_0(y) = [O(c)] \Ff \phi_0(c)$.  
Hence,  $\Ff \phi_0 (c) = [O_c] \inv \CR_c(\Ff \phi_0)$,
and the lemma follows from the first formula and  \propref{E2}.
 \eprf

Note that $\theta'$ above is easily computed, and gives the same (absolute) value for $D, \dot{D}$.    Since at $1$ we have an isomorphism of $D,\dot{D}$ preserving
integral structure, the Fourier transform of $\phi_1$ can be computed with respect
to either ring, giving the same result $\phi_1'$.   Finally, note that the global 
term $[\CT_c : T(K)]$ is the same, via an explicit bijection, for adelically matching
$c,c'$ (\lemref{el-conj}.)

We can now deduce a {\em proof of \thmref{A}}.    Let $\phi_0,\dot{\phi_0}$ be   matching  $R_0^*$-invariant local
test functions at $0$.  We wish to show that   $\Ff \phi_0, \Ff \dot{\phi_0}$   also
  match.  By \lemref{invariance} (c), it suffices to consider rational 
points  $c, \dot{c}$ be   of matching  conjugacy classes $O, \dot{O}$ of $D,\dot{D}$; 
with $c  \notin CN$ (and so $\dot{c} \notin \dot{CN}$.) By \lemref{el-conj} there exists
$\dot{c}'$ adelically conjugate to $c$; by invariance we have $ \Ff \dot{\phi_0}(\dot{c})=
\Ff \dot{\phi_0}\dot{c}'$; so we may assume $c, \dot{c}$ match at every place.  
 In this case, \propref{geometric} and the explicit formula of \lemref{discrete2} shows that $\Ff \phi_0(c) = \Ff \dot{\phi_0} (\dot{c})$.

\>{section}

\section{Appendix 1:     Ind-definable sets}     \label{limited-sets-s}

We include here some standard definitions, largely lifted from the exposition in \cite{CH}.
 

 A structure $N$ for a finite relational language
$L$ is  {\em piecewise-definable} over another structure $k$ if there exist definable $L$-structures
$N_i$ and definable $L$-embeddings $N_i \to N_{i+1}$ such that $\lim_i N_i$ is isomorphic to $N$.   A {\em definable subset} of $\lim_i N_i$ is just a definable subset of some $N_i$.  
A {\em morphism} $f: N \to N'$ is an $L$-embedding such that for any definable $S \subseteq N$, $f(S)$ is a definable subset of $N'$, and $f|S$ is a definable map $S \to f(S)$.

A piecewise definable set is an Ind-object over the category of definable sets with {\em injective}
definable maps; we will not consider other Ind-objects in this paper, so will
use the term "Ind-definable" synonymously with ``piecewise definable''.

Let $C$ be the category of $L$-structures interpretable in $k$, with definable $L$-embeddings
between them.  Since all maps in $C$ are injective, every object of $Ind C$ is strict.
If $A \in Ind C$ is represented by a system $(A_i)_{i \in I}$, 
let $\phi(A) = lim_i A_i$ (inductive limit of $L$-structures.)  For a map $f: A \to B$
in $Ind C$, $\phi(f)$ is defined in the obvious way.  Then $\phi$ is an equivalence of
categories.  Unlike the case of $Pro C$, there is no saturation requirement on $k$.

\<{lem} \label{limited-1}
(1)  If $N$ is quantifier-free definable over $L$, and $L$ is piecewise-definable over $k$, then $N$ is piecewise-definable over $k$. 

 (2)  Let $k$ be a  field, and let $L=k(b_1,\ldots,b_n)$ be a finitely generated field extension of $k$.
Then $(L,+,\cdot,b_1,\ldots,b_n,k)$ is piecewise definable over $k$.    (More precisely there exists a piecewise definable
$k$-algebra $L'$ and an isomorphism $\psi: L  \to L'$ of $k$-algebras.)

(3)  For any variety $V$ over $L$, $V(L)$ can be viewed as   piecewise-definable
over $k$.  (I.e. $\psi (V(L))= V^{\psi}(L')$ is piecewise-definable over $k$.)

 \>{lem}

\prf   (1) is clear.  For (2), $L$ is a finite extension of  a purely transcendental extension 
$k(t)=k(t_1,\ldots,t_n)$ of $k$. Clearly $L$ is  quantifier-free definable over $k(t)$.
Hence by (1) it suffices to show that  $k(t)$  is piecewise-definable
over $k$.  Indeed let  $S_n$
be the set of rational functions $f(t)/g(t)$ with $\deg(f),\deg(g) \leq n$, and let 
$+,\cdot $ be the graphs of addition and multiplication restricted to $S_n^3$.  Then 
$\lim_n S_n = k(t)$.   

(3) Note that the $k$-algebra isomorphism $\psi$ induces a map $V(L) \to V^{\psi}(L')$,
also denoted $\psi$.  (3) follows from (1) and (2).  \eprf

\<{defn}  Let $L,V$ be as in \ref{limited-1}.  A subset $Y$ of $V(L)$ is
called {\em limited} if for some isomorphism $\psi: L \to L'$
to a piecewise-definable field, $\psi(Y)$ is contained in a definable subset of 
the piecewise-definable set $V^\psi (L')$.  \>{defn}

Let us mention two further equivalent formulations , one geometric and one model-theoretic.  

(1)  When $tr. deg._k L = 1$, and when $V$ comes with a projective embedding,
one has the notion of a Weil height of a point of $V(L)$.  Then a limited subset of $V(L)$ is a set of bounded height.    (2)  Let $T$ be the $\omega$-stable theory of pairs $(k,K)$ of algebraically closed fields, with $k < K$.  Assume $(k,K) \models T$.  A subset $Y$ of $V(K)$ is {\em limited} if it is $k$-internal,
i.e. $Y \subseteq \dcl(b,k)$ for some finite $b$.  (In this case $Y  \subseteq V(L)$
 for some subfield $L$ of $K$,
   finitely generated over $k$.)

\<{section}{Appendix 2:  forms}

\<{defn}  Let $T$ be a   theory, $D$ a definable set, and $R_i$ a definable subset of $D^{m_i}$ ($i=1,\ldots,m$).  By a {\em form} of $(D,R_i)_i$ we mean a structure  $(D',\dot{R}_i)_{i}$,
with $D'$ definable in $T$ and $\dot{R}_i$ a $T$-definable subset of $D^{m_i}$, such that
for any $M \models T$ there exists a $T_M$-definable isomorphism $(D,R_i)_i \to (D',\dot{R}_i)_i$,
i.e. a $T_M$-definable bijection $D \to D'$ carrying $R_i$ to $\dot{R}_i$. \>{defn}

For instance, a torus is by definition a form of $G_m^n$, with respect to the theory $ACF$.  

For the rest of the section we discuss forms for $ACVF$ or $ACVF_F$, where 
  $F$ is a valued field, with residue field $\f$.     Let $\Oo$ denote the valuation ring, $\Mm$ the maximal ideal.  
  
Let $M_n$ (respectively $ M_n(\Oo)$) denote the ring of $n \times n$ (integral) matrices.
Thus $D$ is a form of $(M_n,M_n(\Oo))$ iff $D$ is a form  definable finite dimensional central simple algebra. 

A  form  of $(M_n,M_n(\Oo))$ over $F$ is  a pair $(D,R)$, with $D$ an $ACF_F$-definable algebra and $R$ a definable subring, 
such that if $K \models ACF_F$ then there exists an $ACF_K$-definable isomorphism
$D \to M_n$ carrying $M_n$ to $M_n(\Oo)$.    

If $V$ is a definable vector space, by a {\em lattice}  we mean a definable $\Oo$-submodule $\Lambda$ of $V$, such $(V,\Lambda)$ is a form of $(K^n,\Oo^n)$ (for
$n=\dim(V)$.) 

 Thus    $(D,R)$ is a form of $(M_n,M_n(\Oo))$ iff there exists a 
   an $ACF_K$-definable $D$-module $A$ of dimension $n$, and a definable lattice $\Lambda \leq A$, such that $R = \{r \in D: r \Lambda = \Lambda\}$.  If $A,\Lambda$
   can be found over an unramified extension of $F$, we say that $(D,R)$ is
   an {\em unramified} form.  While it is mostly   unramified forms that are of interest for us, much of the discussion can be carried out more generally.
   
     Since $GL_n(\Oo)$ leaves invariant the ideal $\Mm M_n(\Oo)$, any
definable integral form   $R$ has a unique definable ideal $M$, such that $(D,R,M)$ is a form
of $(M_n,M_n(\Oo), \Mm M_n(\Oo))$.

The trace map $tr \circ \phi$ does not depend on the choice of $\phi$, so it is defined over $F$,
and we denote it by $tr$.  Similarly for $\det$.   In particular we have a bilinear form $tr(xy)$ defined over $F$.

\ssec{Characterizations of integral forms}  

For an $ACVF_F$- definable ring $R$, we will say ``definably compact'' for ``$(R,+)$ is generically metastable", i.e. for: ``$(R,+)$ admits a stably dominated translation invariant type" (\cite{HHM}.)
 For a subring of an algebra $D$, this is equivalent to:
$R$ is bounded, and definable by weak valuation inequalities.  
 In this case,
for some $\phi, \theta$, $R$ is defined by a formula $\phi(x,a)$, with $a \in \theta(F)$, and for almost all local fields
$F'$ and $a' \in \theta(F')$, $\phi(x,a')$ defines a compact ring.   Say $R$ is ``maximally definably compact"  if it is definably compact and 
  is contained in no bigger definably compact ring, even over $K$.  

Given a definable lattice $\Lambda \leq D$, let $\Lambda^\perp = \{x: (\forall y \in \Lambda) (tr(xy) \in \Oo)$.   This is another definable lattice, freely generated as an $\Oo$-module by the dual basis
to a basis for $\Lambda$.  Say $\Lambda$ is self-dual if $\Lambda^\perp = \Lambda$.

Let $R$ be definably compact, $K \models ACVF_F$.
Then $R$ is contained in a conjugate $\dot{R}$ of $GL_n(\Oo)$ (defined over $K$). 
 Any such conjugate is self-dual.  If 
  $R$ is also self-dual, $\dot{R} \subseteq R$, so $R=\dot{R}$.  Conversely, 
if $R$ is maximally definably compact then $R=\dot{R}$, so $R$ is self-dual.  

Thus the following conditions on a definably compact subring $R$ are equivalent: $R$ is  self-dual, $R$ is  a maximal definably compact,  $R$ is (eventually, i.e. in a model) conjugate to $M_n(\Oo)$.

For any definable subring $R$ of $GL_n$, let $N(R) =  \{a \in GL_n: (\forall b \in R) (a \inv b a \in R ) \}$
be the normalizer.  This is a definable subgroup of $GL_n$ containing the center $Z$.  
We call $R$   {\em self-normalizing} if $N(R) = Z R^*$.  

\<{lem} \lbl{match}  Let $R$ be a self-normalizing subring of $D$, a form of $M_n$.    Let $(\dot{D},\dot{R})$
be a form of $(D,R)$.  Let $E$ be the definable equivalence
relation of $R^*$-conjugacy on $D$, and similarly $\dot{E}$.  Let $D/E, \dot{D}/\dot{E}$
be interpreted in ACVF.  There exists a definable bijection $f: D/E \to \dot{D} / \dot{E}$.  
\>{lem}

\prf  Let $Hom^*(A,B)$ denote the set of $k$-algebra isomorphisms $A \to B$. 
The $Hom^*(D,\dot{D})$ is a form of $Hom^*(M_n,M_n) = PGL_n$; in particular it is a definable
set.  Let $H = \{h \in Hom^*(D,\dot{D}): h(R) = \dot{R} \}$.  Then $H$ is also a definable
set,  a torsor for $N(R)/Z$.  Any $h \in H$ induces a bijection $D/E \to  \dot{D} / \dot{E}$,
which is $h$-definable.  However this bijection does not depend on the choice of $h$,
so it is definable.
\eprf

Though we will use only forms of $M_n(\Oo)$, we note in passing   another, non-maximal, compact definable ring.

\<{example} \lbl{iwahori}  \rm   
Let $I$ be the Iwahori algebra of $n \times n$ matrices from $\Oo$ 
with superdiagonal entries in $\Mm$.   $I$ can be viewed as a definable subring of the matrix ring $M_n$.
Then $I$ is self-normalizing as an ACVF-definable ring, though $I(\Qq_p)$ or $I(\Cc((t))$ are not.
  First, 
if $a \in GL_n$ normalizes $I$, then it must normalize $M_n(\Oo)$; the reason is that
by considering elements of the form $t^\alpha$ for $\alpha \to 0$, one can approximate
elements of $M_n(\Oo)$ by elements of $I$.  Since in $M_n(\Oo / \Mm)$, the algebra
of lower triangular matrices is self-normalizing, it follows that $I$ is self-normalizing in $GL_n$, modulo the center.     \>{example}



 \<{lem}  \lbl{r1}
Let $D$ be a form  of $M_n$ over a nontrivially valued field $F$, and let $K \models ACVF_F$.  
Then there exists an $ACVF_F$-definable subring $R$ such that 
there exists an $ACF_K$-definable isomorphism $h: D \to M_n$ with $h(R) = M_n(\Oo)$.  
    \>{lem} 
    
\prf   Since $F^a \models ACVF_F$, there exists a finite $ACVF_F$ definable set $W$ 
and   for $w \in W$ an    $F(w)$-definable
representation $V(w)$ of $D$, of dimension $n$.  Given $w,w' \in W$, there exists
a finite $F(w,w')$-definable set $Y_w$ and for $y \in Y$ an $F(w,w',y)$- definable isomorphism $g_y$ between the two representations. 
Moreover two 
isomorphisms $g_y,g_{y'}: V_w \to V_{w'}$ differ by a scalar $c(y,y')$.   We may assume all $y \in Y$ have the same type over $F(w,w')$.  But we can define a partial ordering on $Y$,
$y \leq y'$ iff $\val( c(y,y')) \geq 0$.  Since $Y$ is finite the partial ordering must be trivial,
i.e. $\val c(y,y') = 0$ for all $y,y'$.   It follows that the maps $g_y$ induce a unique
isomorphism $g_{w,w'}: \Oo^* \bs V_w \to \Oo^* \bs V_{w'}$.    (Alternatively
by Hilbert 90, since $Hom_{D}(V_w,V_{w'})$ is a 1-dimensional vector space,
there exists $G_{w,w'}: V_w \to V_{w'}$, defined over $F(w,w')$; let
$g_{w,w'}:  \Oo^* \bs V_w \to \Oo^* \bs V_{w'}$ be the induced map. )

Now we need some Galois cohomology, which is easiest to do from first principles.  
Let $C^m(W,\G)$ be the set of definable maps $W^m \to \G$.  Define a coboundary 
map $d: C^m(W,\G) \to C^{m+1}(W,\G)$ in the usual way.  Namely given
$f \in C^m(W,\G)$, let $d f = F$ where $F(w_0,\ldots,w_m) = \sum (-1)^i f({\rm omit \ } w_i)$.
Since $\G$ is uniquely divisible,  the cohomology groups $H^m(W,\G)$, $m \geq 1$ are trivial:
let $f \in C^m(W,\G)$ and assume $df=0$.   Given $w \in W$, let $f_w  \in C^{m-1}(W,\G)$ be defined
by $f_w(w_1,\ldots,w_{m-1})= -f(w,w_1,\ldots,w_{m-1})$. Then formally the relation
$df=0$ gives $d f_w = f$.  Now $f_w$ is not definable, but let $F$ be
the average of all $f_w$; then $F$ is definable and $dF = f$.

In particular, let $f(w,w',w'') = g_{w,w''}\inv g_{w',w'' } g_{w,w'}$.  Then $f(w,w',w'')$ is a scalar
endomorphism of $V_w$, modulo $\Oo^*$, so it can be viewed as an element of $\G$.  We have $df=0$ (this can be verified by fixing some $w_0$, identifying $\Oo^* \bs V_w$
with $\Oo^* \bs V_{w_0}$ via $g_{w_0,w}$, and using the commutativity of $\G$.)
    So $f=d F$ for
some $F \in C^2(W,\G)$. Replacing $g_{w,w'}$ by $F(w,w') \inv g_{w,w'} $
(i.e. by $g_{w,w'}$ composed by the endomorphism of division by a scalar with value
$F(w,w')$),  we may assume 
$f=0$, i.e. $g_{w,w''} = g_{w',w''} g_{w,w'} $.  So we have a commuting  system of isomorphisms between 
the integrally projectivized representations $\Oo^* \bs V_w$.

Now we can find an $F(w)$-definable lattice $\Lambda_w \in L(V_w)$.   
Let $\Lambda^*_w = \meet_{w' \in W} g_{w,w'} \inv L(V_{w'})$.  Then 
$\Lambda^*_w$ is also an  $F(w)$-definable lattice
\footnote{By \cite{acvf1}, a definable $\Oo$-submodule of $n$-space is a lattice iff
the intersection with each one-dimensional subspace is a closed ball; this property is
evidently preserved under finite intersections.}
; and $g_{w,w'} \inv \Lambda^*_{w'} = \Lambda^*_w$.

Let $R_w = \{r \in D:  r \Lambda^*_w \subseteq  \Lambda^*_w \} $.  Then $R_w$
does not depend on $w$ as one sees using $g_{w,w'}$.  Let $R=R_w$.  Then $R$ clearly
satisfies the requirements.  

\eprf

Let $T_n$ be the diagonal subalgebra of $M_n$.  
Consider diagonalizable algebraic subrings of $D$, i.e. definable subrings $T$ such that 
$(D,T)$ is a form of $(M_n,{T_n})$.  Then there exists a unique definable subring $ \Oo_T$
of $T$, such that $(D,T, \Oo_T)$ is a form of $(M_n,{T_n}, {T_n}(\Oo))$.  Indeed after base change, there exists an isomorphism $(D,T) \to (M_n,{T_n})$ of pairs of rings; it is well-defined up to
composition with an element of the Weyl group $\Sym(n)$; since $\Sym(n)$ respects ${T_n}(\Oo)$, the pullback of ${T_n}(\Oo)$ does not depend on the choice of isomorphism, and is definable.   
 It is not necessarily the case that $\Oo_T $ is the $\Oo$-module generated by
$\Oo_{T(F)}$.  Note that 
by Hilbert 90, $T$ is determined by $T(F)$.

 We call $R$ a 
 {\em definable integral form} for $(D,T)$ if $(D,T,R)$ is a   form of $(M_n,{T_n},M_n(\Oo))$.    So $R \meet T = \Oo_T$.

 \<{lem}  \lbl{r2}
Let $D$ be a form  of $M_n$ over a nontrivially valued field $F$, and let $K \models ACVF_F$.  Let $T$ be   defined over $F$, with $(D,T)$ a  form of $(M_n,{T_n})$.
Then there exists definable integral form for $(D,T)$.       \>{lem}

\prf Same as \lemref{r1}.  We take $V(w)$ to be graded by one-dimensional eigenspaces
of $T$, and we choose $\Lambda_w$ to be generated by $T$-eigenvectors.  \eprf

In case $D \cong M_n$, we have $D \cong End(T)$; let $R_T = End_\Oo \Oo_T$; then
 $(D,T,R_T)$ is a  definable integral form for $(D,T)$.

Let  $T$ be a diagonalizable subring, with normalizer $N$.
  Let $\mD=\mD(D,T,F)$ be the set of definable integral forms of $(D,T)$  over $F$.  
 Let  $Z$ be the center of $D^*$.
  Let $X =Hom(G_m,T^*/Z)$
and $\Delta = X  \tensor \G$ (so $\Delta$ is parameterically isomorphic to $\G^{\dim(T)-1}$.)  Let $\mD_T^F$ be the  definable integral forms for $(D,T)$,
and $\mDTc ^F= \mD_T^F /  {N(F)^*} $
be $\mD_T(F)$ 
 up to $D(F)^*$-conjugacy.   If the identity of $F$ is clear we will omit the superscript.
 
If $K \models ACVF$ and $R,\dot{R}$ are two definable integral forms for $D$ (or for $(D,T)$),
then $R,\dot{R}$ are  conjugate in $D(K)$ (respectively, in $R(K) Z(K)$.)   We will use this in (1),(2) below. 
 
In (4),(5) below we use the fact that $D$ splits over the maximal unramified algebraic extension $F^{unr}$ of $F$; see \cite{serre} II 3.2, Corrolaire, and 3.3(c) ($H^1(F^{unr}, PGL^n) = 0$.)
 
 In (2) below we assume residue char. 0; in each case what we really
use is that $H^1(F,A)=0$ for certain definable unipotent groups $A$.   For ACVF, unlike
ACF, this is not automatic even over perfect fields; but   the instances we require may be true in char. p too.

\<{lem} \lbl{r3}  \<{enumerate}
\item   $\mDTc ^F$ has  a canonical structure  of torsor over 
$\Delta_{def}(F) / \Delta(F)$; where $\Delta_{def}(F)$ is
the set of points of $\Delta(F^{alg})$  invariant under $Aut(F^{alg}/F)$, 
and $\Delta(F)$ is the set of points of $\Delta$ represented in $F$.

\item  (char. 0.)    Assume  there exists a definable $B$  such that $(D,T,B)$ is a form
of $(M_n,{T_n},B_n)$, with $B_n$ the upper triangular matrices.   Then   there   is a canonical retraction $\rho:\mD \to \mD_T$.     In residue char. 0, 
 any $R \in \mD$, $R$ and $\rho(R)$ are $D^*(F)$-conjugate.  Hence  every definable integral form for $D$ is $D(F)$-conjugate
to a definable integral form for $(D,T)$.
\item Let $F'/F$ be an unramified field extension.  Then the natural map $\mDTc^F  \to \mDTc^{F'}$ is injective.  
 

\>{enumerate}
  \>{lem}

\prf   
\<{enumerate}
\item   
Let $R \in \mD_T$.  Let $N \leq D^*$ be the normalizer of $T$ (a definable group.)
As observed above, any   element   of $\mD_T$ has the form $a \inv R a$ for some $a \in N(K)$.
   But $N R = T^* R$, since the Weyl group is represented in $R(K)^*$.
So we can take $a \in T(K)$.  Now $a \inv R a$ is by assumption $ACVF_F$-definable,
and the normalizer of $R$ is $R^*Z$, so 
$N_T(R) := T \meet R = \Oo_T Z$, using the observation above that $T \meet R = \Oo_T$.  hence $a \Oo_T Z$ is definable.  Now $T/ (\Oo_TZ) = \Delta$.
The definable points of $\Delta$ are $\Delta(F^{alg})$.  This gives a surjection $\Delta_{def}(F)    \to \mD_T$,
and hence   $\Delta_{def}(F)  / \Delta(F) \to \mDTc$.  It is easy to check injectivity.

\item    Any   element   of $\mD$ has the form $a  GL_n(\Oo) a \inv$ for some $a \in GL_n(K)$.
 Since   $GL_n= B_n GL_n(\Oo)$,   we can take $a \in B_n(K) $.  Let $ss(x)$ be the semi-simple part of $x$.  Then $ss$ commutes with conjugation, hence gives a well-defined map on $D$;
 it induces a homomorphism $B_n \to {T_n}$.  This in turn induces a   map
 $s: B_n / B_n(\Oo) Z \to {T_n} / {T_n}(\Oo) Z$.

  It remains to show that $\rho(R),R$ are $D^*(F)$-conjugate.  Say 
  $R=  a  GL_n(\Oo) a \inv$  with $a \in B_n$; write   $a=a_s a_u$ with $a_s \in {T_n}$
 and $a_u \in U_n$, this being the strictly upper-triangular matrices.   So
  $\rho(R) = a_s GL_n(\Oo) a_s \inv $, hence $R= a_u \rho(R) a_u \inv$.
 Now $S=\{u \in U_n: u \rho(R) u \inv = \rho(R)\}$ is an $ACVF_F$-definable subgroup of $U_n$.
This group is geometrically connected.  In characteristic $0$ it is clear that $H^1(Aut(F^a/F),S)=0$, so there exists a definable point.

\item   Pick $R \in  \mD_T$.  Then by change of scalars we can view $R$ as an element of $\mD_T$
  over $F'$.    By (1), $\mDTc ^F$ corresponds bijectively to a subgroup of $\Delta(F^{alg})/\Delta(F)$, 
while  $\mDTc ^{F'}$ corresponds bijectively to $\Delta(F^{alg})/\Delta(F')$.  However these two groups are the same.

%

 \>{enumerate}
\eprf

\<{rem}  \rm 
Let $G_{ch}$ be a Chevalley group.  Let $G$ be a form of $G_{ch}$ over $F$.
Define an   {\em integral form} of $G$   to be
 an $ACVF_F$-definable subgroup $H$ of $G$, such that $(G,H)$ is a form of
 $G_{ch}, G_{ch} (\Oo)$.   It seems that analogues of the above results should be true. 
 
Above we used the fact that  $M_n$ has no outer automorphisms, which is not true for $G_{ch}$.  However 
every outer automorphism of $G_{ch}(K)$ is an inner automorphism composed with a graph automorphism, and the graph automorphisms preserve $G_{ch}(\Oo)$.  This makes it
possible to consider definable integral forms of $G$ for a form $G$ of $G_{ch}$, so that
two forms are $G$-conjugate. 
   \>{rem} 

\>{section}

 \<{section}{Appendix 3:  Multiplicative convolution.}  \lbl{convolution}
 
Our results on the stability of the Fourier transform have an immediate
consequence for additive convolution:  given two pairs of   matching
local test functions on $D,\dot{D}$, their additive convolutions also match.

This statement can be phrased without the intervention of   additive characters,
and may be valid for the Grothendieck ring $\K[Gr \inv]$; our proof however
shows it in $\K_e[Gr \inv]$. Indeed the Fourier transform transposes the problem
into a similar one using pointwise products, which is obvious.  

Here we assume characteristic $0$ in order to point out a relation between
 this additive result and the analogous multiplicative statement, as in  \cite{DKV}.
 
We note that the ``orbit method" isomorphism between convolution algebras of 
nilpotent groups and their Lie algebras (\cite{BD}  Prop. 2.4) goes through for motivic convolution algebras.   

In our setting, we have the division algebra $D$, with subring $R$, ideals $M_n$ of elements of determinental valuation $\geq n$.   The exponential map defines a bijection 
$x \mapsto 1 + x + \ldots$, $A_n:= M_1 / M_n \to (1+M_1) / (1+M_n) =: G_n$.  This induces a bijection $\exp$ between $D^*$-conjugacy classes on the algebraic groups $A_n$ and on $G_n$.  

\<{lem}  $\exp$    induces an isomorphism of motivic rings $Fn(G_n)^{D^*} \to Fn(A_n)^{D^*}$.  \>{lem}

\prf  
Let $C_1,C_2$ be two conjugacy classes, and let $c \in A_n$.
Then 
we have to show that 
$A=\{(x_1,x_2) \in C_1 \times C_2: x_1+x_2 = c \}$
has the same class in the Grothendieck group as 
 $B=\{(x_1,x_2) \in C_1 \times C_2: exp(x_1) exp(x_2) = exp(c) \}$.
 
 Let $H(x_1,x_2) = (\exp(ad \phi(x,y) ) (x) , \exp(ad \psi(x,y) ) (y))$, where
 $\phi(X,Y),\psi(X,Y)$ are the Lie polynomials from   \cite{BD}  Lemma. 2.5.  
 By this lemma, $B = H \inv (A)$.   We are thus done given:  

\Claim{}   Let $\phi(X,Y),\psi(X,Y)$ be Lie polynomials.  Then the function $A_n^2 \to A_n^2$ defined by:
$$(x,y) \mapsto (\exp(ad \phi(x,y) ) (x) , \exp(ad \psi(x,y) ) (y))$$
is bijective.

\prf of Claim:  if $H(u) = H(u')$, we show by induction on $k \leq n$ that
$u \equiv u' \mod M_k$.  Given that $u=u' \mod M_k$ we have
$\phi(u) = \phi(u'),  \psi(u)=\psi(u')   \mod M_{k+1}$, so $u=u' \mod M_{k+1}$.  \eprf
\eprf

In the classical case, given the isomorphism on $U$,  the full multplicative isomorphism can be obtained using  character-theoretic methods.  As the characters involved are uniformly definable, it seems likely that this can be done motivically too.

\>{section} 

\<{thebibliography}{delon}

\bibitem{baur}
Baur, Walter On the elementary theory of pairs of real closed fields. II.  J. Symbolic Logic  47  (1982), no. 3, 669--679

\bibitem{BD}  Mitya Boyarchenko, Vladimir Drinfeld.   A motivated introduction to character sheaves and the orbit method for unipotent groups in positive characteristic. 
math.RT/0609769 

\bibitem{CK} 
Chang, C. C.; Keisler, H. J. Model theory. Third edition. Studies in Logic and the Foundations of Mathematics, 73. North-Holland Publishing Co., Amsterdam, 1990.  

\bibitem{CH}   	Difference fields and descent in algebraic dynamics - I. Zo\'e Chatzidakis, Ehud Hrushovski. math.LO arXiv:0711.3864 

\bibitem{delon}
 Delon, Fran\c{c}oise Ind\'ecidabilit\'e de la th\'eorie des paires imm\'ediates de corps valu\'es henseliens
 J. Symbolic Logic 56 (1991), no. 4, 1236--1242
 
\bibitem{acvf1}   Haskell, Deirdre; Hrushovski, Ehud; Macpherson, Dugald Definable sets in algebraically closed valued fields: elimination of imaginaries.  J. Reine Angew. Math.  597  (2006), 175--236. 

\bibitem{HHM} Haskell, Deirdre; Hrushovski, Ehud; Macpherson, Dugald Stable domination and independence in algebraically closed valued fields. Lecture Notes in Logic, 30. Association for Symbolic Logic, Chicago, IL; Cambridge University Press, Cambridge, 2008. xii+182 pp

\bibitem{jacobson} Nathan Jacobson, Basic Algebra,  San Francisco, W. H. Freeman [1974-80], 2 v. illus. 25 cm.

 \bibitem{leloup1}
 Leloup, G\'erard, Th\'eories compl\`etes de paires de corps valu\'es henseliens.   J. Symbolic Logic 55 (1990), no. 1, 323--339.


 \bibitem{GK}
  Gaitsgory, D.; Kazhdan, D. Representations of algebraic groups over a 2-dimensional local field. Geom. Funct. Anal. 14 (2004), no. 3, 535--574.

\bibitem{HK} Hrushovski, Ehud; Kazhdan, David
Integration in valued fields.  Algebraic geometry and number theory, 261--405,
Progr. Math., 253, Birkhauser Boston, Boston, MA, 2006. 

\bibitem{Kam}  Moshe Kamensky,   Ind- and Pro- definable sets, math.LO/0608163

\bibitem{DKV}
Deligne, P.; Kazhdan, D.; Vign\'eras, M.-F. Repr\'esentations des alg\`ebres centrales simples $p$-adiques.   Representations of reductive groups over a local field,  33--117, Travaux en Cours, Hermann, Paris, 1984. 

\bibitem{DL} Denef, J.; Loeser, F. Motivic integration and the Grothendieck group of pseudo-finite fields.  Proceedings of the International Congress of Mathematicians, Vol. II (Beijing, 2002),  13--23, Higher Ed. Press, Beijing, 2002

\bibitem{lichtenbaum} Stephen Lichtenbaum,  Duality Theorems for Curves over P-adic Fields,  Inventiones math. 7, 120-136 (1969)
 
\bibitem{MvD} van den Dries, Lou; Macintyre, Angus The logic of Rumely's local-global principle.  J. Reine Angew. Math.  407  (1990), 33--56 
 
\bibitem{sebag}
Sebag, Julien IntŽgration motivique sur les schŽmas formels. (French) [Motivic integration on formal schemes]  Bull. Soc. Math. France  132  (2004),  no. 1, 1--54.

     \bibitem{serre}   Serre, Jean-Pierre Lectures on the Mordell-Weil theorem. Translated from the French and edited by Martin Brown from notes by Michel Waldschmidt. Aspects of Mathematics, E15.   Vieweg  Braunschweig, 1989. 

\bibitem{serre-loc}
Serre, Jean-Pierre Local fields. Translated from the French by Marvin Jay Greenberg. Graduate Texts in Mathematics, 67. Springer-Verlag, New York-Berlin, 1979. viii+241 pp.

\>{thebibliography}

\end{document}